\documentclass[american,english]{article}
\usepackage[T1]{fontenc}
\usepackage[latin9]{inputenc}
\usepackage{color}
\usepackage{units}
\usepackage{mathrsfs}
\usepackage{amsmath}
\usepackage{amsthm}
\usepackage{amssymb}
\usepackage{setspace}
\usepackage{wasysym}
\makeatletter
\theoremstyle{plain}
\newtheorem{thm}{\protect\theoremname}[section]
\theoremstyle{remark}
\newtheorem*{acknowledgement*}{\protect\acknowledgementname}
\theoremstyle{definition}
\newtheorem{defn}[thm]{\protect\definitionname}
\theoremstyle{plain}
\newtheorem{lem}[thm]{\protect\lemmaname}
\theoremstyle{plain}
\newtheorem{cor}[thm]{\protect\corollaryname}
\theoremstyle{plain}
\newtheorem{prop}[thm]{\protect\propositionname}
\theoremstyle{remark}
\newtheorem{claim}[thm]{\protect\claimname}
\newcommand{\lyxaddress}[1]{
	\par {\raggedright #1
	\vspace{1.4em}
	\noindent\par}
}

\usepackage{enumitem}
\setenumerate{label=(\roman{enumi})}
\usepackage{tikz}

\makeatother

\usepackage{babel}
\addto\captionsamerican{\renewcommand{\acknowledgementname}{Acknowledgement}}
\addto\captionsamerican{\renewcommand{\claimname}{Claim}}
\addto\captionsamerican{\renewcommand{\corollaryname}{Corollary}}
\addto\captionsamerican{\renewcommand{\definitionname}{Definition}}
\addto\captionsamerican{\renewcommand{\lemmaname}{Lemma}}
\addto\captionsamerican{\renewcommand{\propositionname}{Proposition}}
\addto\captionsamerican{\renewcommand{\theoremname}{Theorem}}
\addto\captionsenglish{\renewcommand{\acknowledgementname}{Acknowledgement}}
\addto\captionsenglish{\renewcommand{\claimname}{Claim}}
\addto\captionsenglish{\renewcommand{\corollaryname}{Corollary}}
\addto\captionsenglish{\renewcommand{\definitionname}{Definition}}
\addto\captionsenglish{\renewcommand{\lemmaname}{Lemma}}
\addto\captionsenglish{\renewcommand{\propositionname}{Proposition}}
\addto\captionsenglish{\renewcommand{\theoremname}{Theorem}}
\providecommand{\acknowledgementname}{Acknowledgement}
\providecommand{\claimname}{Claim}
\providecommand{\corollaryname}{Corollary}
\providecommand{\definitionname}{Definition}
\providecommand{\lemmaname}{Lemma}
\providecommand{\propositionname}{Proposition}
\providecommand{\theoremname}{Theorem}

\begin{document}
\title{Moment-optimal finitary isomorphism for i.i.d. processes of equal
entropy}
\author{Uri Gabor\thanks{supported by ISF research grants 1702/17 and 3056/21.}}
\maketitle
\begin{abstract}
The finitary isomorphism theorem, due to Keane and Smorodinsky, raised
the natural question of how ``finite'' the isomorphism can be, in
terms of moments of the coding radius. More precisely, for which values
$\theta>0$ does there exist an isomorphism $\phi$ between any two
i.i.d. processes of equal entropy, with coding radii $R_{\phi},R_{\phi^{-1}}$
exhibiting finite $\theta$-moments? \cite{key-19,key-1}. Parry \cite{key-6}
and Krieger \cite{key-8} showed that those $\theta$ must be lesser
than 1 in general, and Harvey and Peres \cite{key-9} showed that
$\theta<\frac{1}{2}$ in general. However, the question for $\theta\in(0,\frac{1}{2})$
remained open, and in fact no general construction of an isomorphism
was shown to exhibit any non trivial finite moments.

In the present work we settle this problem, showing that between any
two aperiodic Markov processes (and i.i.d. processes in particular)
of the same entropy, there exists an isomorphism $\phi$ with coding
radii exhibiting finite $\theta$-moments for all $\theta\in(0,\frac{1}{2})$.
The isomorphism is constructed explicitly, and the tails $\mathbf{P}[R_{\phi}>n],\mathbf{P}[R_{\phi^{-1}}>n]$
are shown to be of order $O(n^{-\frac{1}{2}}\log^{O(1)}n)$, which
is optimal up to a poly-logarithmic factor.
\end{abstract}
The coding radius of a function $f$ defined on a stationary process
$X\in A^{\mathbb{Z}}$ is the stopping time function $R:A^{\mathbb{Z}}\rightarrow\mathbb{N}\cup\{\infty\}$
defined by
\[
R(x)=\min\left\{ r\geq0:\,(X_{i}=x_{i})_{|i|\leq r}\quad\Longrightarrow\quad f(X)=f(x)\;a.s.\right\} .
\]
A stationary map $\phi$ between two processes is \textit{finitary}
if its zero-coordinate function $\phi_{0}$ has an a.s. finite coding
radius. That coding radius is denoted by $R_{\phi}$. $\phi$ is an
\textit{isomorphism} if it has a stationary inverse map, and if both
$\phi$ and $\phi^{-1}$ are finitary, we say that $\phi$ is a \textit{finitary
isomorphism}, and the two processes are said to be\textit{ finitarily
isomorphic}.

In ergodic theory, finitary (and stationary) maps provide an alternative,
stronger notion of morphism than the more common notions of measurable
factors and isomorphisms.

In the same year Kolmogorov introduced entropy as an invariant for
isomorphism \cite{key-10}, one of the first examples of finitary
isomorphism between i.i.d. processes was given by Mešalkin \cite{key-4},
showing that for a certain class of i.i.d. processes, entropy is a
full invariant for isomorphism. The celebrated Isomorphism Theorem
due to Ornstein \cite{key-5} states that this is the case for the
general class of i.i.d. processes: two i.i.d. processes with the same
entropy are isomorphic. While Ornstein\textquoteright s result implies
that of Mešalkin, the construction of Mešalkin has a few advantages
that are missing in Ornstein\textquoteright s work: The isomorphism
in Ornstein\textquoteright s result is not given explicitly but rather
shown to exist. In contrast, not only is Mešalkin\textquoteright s
code given explicitly, but it has the additional property of being
finitary. From the topological point of view, this is equivalent to
saying that on a set of full measure, the map is continuous. From
the combinatorial point of view, a map is finitary if it can be implemented
by a machine, as it requires a finite set of data for encoding each
of the image\textquoteright s symbols. 

The gap between Mešalkin\textquoteright s code and the Isomorphism
Theorem has narrowed a decade later, when Keane and Smorodinsky proved
that any two finite alphabet i.i.d. processes of the same entropy
are finitarily isomorphic \cite{key-2}. In other words, on the class
of finite alphabet i.i.d. processes, entropy is a full invariant even
for the stricter notion of finitary isomorphism.

While finitary coding has the advantage of the coding radius $R$
being finite a.s., for that to be applicable one has to know how small
$R$ can be made. It is clear that in general $R$ cannot be bounded
by any constant, but it was not clear whether one can require finite
expectation for $R$, i.e. $\mathbf{E}[R]<\infty$, or any finite
moment, $\mathbf{E}[R^{\theta}]<\infty$ for some $\theta>0$. 

For factoring an i.i.d. process on a strictly lower entropy process,
Keane and Smorodinsky \cite{key-3} provided in any such situation
a factor map which is finitary and has finite expected coding radius,
as demonstrated later by Serafin \cite{key-11}. That was significantly
improved by Harvey, Holroyd, Peres and Romik \cite{key-1}, showing
that in that case of strictly lesser entropy, there is always a finitary
factor map with coding radius $R$ having exponential tails: 
\[
\mathbf{P}[R>n]=O(e^{-cn})\text{ for some constant }c>0,
\]
which in particular implies that 
\[
\forall\theta>0,\;\mathbf{E}[R^{\theta}]<\infty.
\]

While factoring one i.i.d. process on another with strictly lower
entropy appears to enjoy good bounds on the coding radius tail, for
isomorphisms negative results have piled up: It was shown independently
by Parry \cite{key-6} and Kreiger \cite{key-8} that finitary isomorphism
with finite expected coding radius imposes certain restrictions on
the relation between the i.i.d. processes in question. Based on these
works together with \cite{key-12}, Schmidt \cite{key-7} showed that
two i.i.d. processes can only be finitarily isomorphic with finite
expected coding radius if their marginal probabilities are permutations
of each other. In 2003 it was shown by Harvey and Peres \cite{key-9}
that if $X$ and $Y$ are i.i.d. processes of the same entropy, and
there is a finitary factor map $\phi(X)=Y$ with coding radius $R_{\phi}$
exhibiting finite half moment, i.e. 
\[
\mathbf{E}[R_{\phi}^{\frac{1}{2}}]<\infty,
\]
then the information functions of the zero marginals in each process
share the same variance:
\[
\mathrm{Var}(I(X_{0}))=\mathrm{Var}(I(Y_{0}))
\]
In light of this result, Mešalkin's isomorphism is optimal, in the
sense that its coding radii achieve all moments less than half: the
code applies to processes that do not share the same marginal informational
variance, thus by Harvey and Peres result the coding radii of any
isomorphism between them do not attain finite half moment. However,
both coding radii of Mešalkin's code satisfy
\[
\mathbf{P}[R>n]=\Theta(n^{-\frac{1}{2}})
\]
thus for any $\theta<\frac{1}{2}$,
\[
\mathbf{E}[R_{\phi}^{\theta}]<\infty.
\]

The main result of this work is that in the class of i.i.d. processes,
as well as for the more general class of aperiodic Markov shifts (irreducible
Markov chains of period one), entropy is a full invariant for finitary
isomorphism with finite less-than-half moments. 
\begin{thm}
\label{thm:Main=000020Thm}For any two finite alphabet aperiodic Markov
shifts $X\sim(A^{\mathbb{Z}},\mu),Y\sim(B^{\mathbb{Z}},\nu)$ of the
same entropy, there exists a finitary isomorphism $\phi:A^{\mathbb{Z}}\rightarrow B^{\mathbb{Z}}$
with coding radii $R_{\phi},R_{\phi^{-1}}$ satisfying:
\begin{equation}
\mathbf{P}[R_{\phi}>n]=O(n^{-\frac{1}{2}}\log^{O(1)}n),\qquad\mathbf{P}[R_{\phi^{-1}}>n]=O(n^{-\frac{1}{2}}\log^{O(1)}n).\label{eq:main=000020theorem}
\end{equation}
Consequently, both coding radii achieve finite $\theta$-moment for
all $\theta<\frac{1}{2}$.
\end{thm}

Using $\tilde{O}$-notation (see subsection \ref{subsec:Definitions-and-a}),
the estimate (\ref{eq:main=000020theorem}) can be written as
\begin{equation}
\mathbf{P}[R_{\phi}>n]=\tilde{O}(n^{-\frac{1}{2}}),\qquad\mathbf{P}[R_{\phi^{-1}}>n]=\tilde{O}(n^{-\frac{1}{2}}).\label{eq:main=000020theorem-1}
\end{equation}

\subsubsection*{Proof outline of Theorem \ref{thm:Main=000020Thm}}

We begin by describing the method of Keane-Smorodinsky \cite{key-2}
for constructing a finitary isomorphism $\phi$ between two i.i.d.
processes $X\sim(A^{\mathbb{Z}},\mathbf{p}^{\mathbb{Z}})$ and $Y\sim(B^{\mathbb{Z}},\mathbf{q}^{\mathbb{Z}})$.
Consider first the case where $\mathbf{p}$ and $\mathbf{q}$ share
an entry with the same probability, say, $\mathbf{p}(a)=\mathbf{q}(b)$.
Using this assumption, the first step is to couple $(X,Y)$ along
the process 
\[
(\mathbf{1}_{\{X_{n}=a\}})_{n\in\mathbb{Z}}\overset{D}{=}(\mathbf{1}_{\{Y_{n}=b\}})_{n\in\mathbb{Z}}.
\]
A factor $\xi=\xi\left((\mathbf{1}_{\{X_{n}=a\}})_{n\in\mathbb{Z}}\right)$
uses the data in $(\mathbf{1}_{\{X_{n}=a\}})_{n\in\mathbb{Z}}$ to
split simultaneously the sample $X$ and the intended sample $\phi(X)$
into finite words in a stationary manner, $X=\cdots X_{S^{-1}}X_{S^{0}}X_{S^{1}}\cdots$
(here $S^{i}$ are intervals and $X_{S^{i}}=(X_{k})_{k\in S^{i}}$),
and each word $X_{S^{i}}$ is used to encode, at least partially,
the respective word $\phi(X)_{S^{i}}$ of the intended image. The
factor $\xi$ is called marker process, as it marks the places of
splitting the sample. Conditioned on the sample $(\mathbf{1}_{\{X_{n}=a\}})_{n\in\mathbb{Z}}$
and ignoring 1-symboled coordinates, these words can be viewed as
having the laws
\[
\begin{array}{lr}
X_{S^{i}}\sim((A\backslash\{a\})^{\ell},\mathbf{p}_{a}^{\ell}); & \phi(X)_{S^{i}}\sim(B\backslash\{b\})^{\ell}\mathbf{q}_{b}^{\ell})\end{array}
\]
where $\mathbf{p}_{a}(\cdot)=\frac{\mathbf{p}(\cdot)}{1-\mathbf{p}(a)}$,
$\mathbf{q}_{b}(\cdot)=\frac{\mathbf{q}(\cdot)}{1-\mathbf{q}(b)}$,
and $\ell=\ell(S^{i})$ is the number of non-$a$-symboled coordinates
of $X_{S^{i}}$.

In order to encode a word in $(A\backslash\{a\})^{\ell}$ to a word
in $(B\backslash\{b\})^{\ell}$, it is necessary that the image word
would have measure greater or equal that of the source word. As this
cannot be guaranteed for words of two general processes, the spaces
of words are partitioned into classes of words, so that most of the
classes of $(A\backslash\{a\})^{\ell}$-words will have measure lesser
than those of the corresponding $(B\backslash\{b\})^{\ell}$-classes.
The collection of source-word classes of \textquotedblleft small\textquotedblright{}
measure, which altogether consist of most of the space $(A\backslash\{a\})^{\ell}$,
are then successfully encoded into classes in the image. In these
partitions, each class consists of all words that are equal on some
specified subset of indices $\mathcal{J}\subset\{1,...,\ell\}$ to
some specified word $(w_{j})_{j\in\mathcal{J}}$. The room in choosing
for each class the subset $\mathcal{J}$ and the symbols $(w_{j})_{j\in\mathcal{J}}$
it is described by, allows one to control the measure of each of the
classes: if $U\subset(A\backslash\{a\})^{\ell}$ is a class described
by $(w_{j})_{j\in\mathcal{J}}$ as above, then
\[
\mathbf{p}_{a}^{\ell}(U)=\prod_{j\in\mathcal{J}}\mathbf{p}_{a}(w_{j}).
\]
Despite the need of the classes in $(B\backslash\{b\})^{\ell}$ to
be of large enough measure, it is likely that the smaller the measure
of a $(B\backslash\{b\})^{\ell}$-words class is, the larger is the
number $|\mathcal{J}|$ of coordinates the class is specified by,
which, in turn, is the number of coordinates determined by encoding
the $(A\backslash\{a\})^{\ell}$-word $X_{S^{i}}$. To determine the
remained unknown coordinates, this procedure of splitting the sample
and encoding classes of words is done recursively again and again:
As for the splitting, given a sample $(\mathbf{1}_{\{X_{n}=a\}})_{n\in\mathbb{Z}}$,
the marker process $\xi$ actually gives a collection of splits $\left\{ X=\cdots X_{S_{r}^{-1}}X_{S_{r}^{0}}X_{S_{r}^{1}}\cdots:r\in\mathbb{N}\right\} $
. Denoting by $S_{r}$ (for each $r$) the interval drawn from the
$r$th-splitting of $X$ that contains the 0-coordinate, the collection
of splits satisfies $S_{r}\subset S_{r+1}$ and $\ell(S_{r})\longrightarrow\infty$.
In the $r$th step, classes of source words drawn from the $r$th-splitting
are associated to classes of image words in a way compatible with
preceding steps. The words length being larger, allows the classes
of $(B\backslash\{b\})^{\ell}$ to be specified on a larger fraction
of the coordinates, while at the same time, more of the space $(A\backslash\{a\})^{\ell}$
can be grouped into small measure classes, thus in later stages more
samples are encoded and more indices of the image are determined.
This gives eventually a homomorphism $\phi$ from $A^{\mathbb{Z}}$
to $B^{\mathbb{Z}}$. Doing simultaneously similar associations in
the other direction, from $B^{\mathbb{Z}}$ to $A^{\mathbb{Z}}$,
gives an isomorphism. 

For the following discussion we need some more notation: Recall that
for each $r$, $S_{r}$ denotes the interval drawn from the $r$th-splitting
of $X$ that contains the 0-coordinate. For any $n$, let $S(n)$
be the maximal $S_{r}$ contained in $[-n,n]$. The probability that
$\phi(X)_{0}$ is determined by $X_{S_{r}}$ can be estimated by
\[
\mathbf{P}[X_{S_{r}}\text{ determines }\phi(X)_{0}]\geq1-\mathbf{P}[E_{r}]-\mathbf{P}[F_{r}]
\]
where $E_{r}$ denotes the event that the class of $X_{S_{r}}$ was
not successfully encoded, and $F_{r}$ denotes the event that the
class of $X_{S_{r}}$ has been encoded to the class of $\phi(X)_{S_{r}}$,
but the zero coordinate has not been determined by the image class.
As $E_{r}$ mainly consists of those words in the source whose classes
are of too large measure to be encoded, by what we claimed above,
the probability of $E_{r}$ tends to $0$ as $r\rightarrow\infty$.
It also can be easily shown that the probability of $F_{r}$ is bounded
by the expected fraction of coordinates $i\in S_{r}$ so that $\phi(X)_{i}$
is not determined by the class of $\phi(X)_{S_{r}}$. By what was
claimed above, this expectation also tends to 0 as $r\rightarrow\infty$.
Thus the coding is finitary. 

Conditioned on the length $\ell=\ell(S_{r})$ of $S_{r}$, an application
of the Azuma-Hoeffding inequality to the information functions $\left(\log\mathbf{p}_{a}(w_{i})\right)_{i=1}^{\ell}$
enables us to choose the classes in $(A\backslash\{a\})^{\ell}$ and
$(B\backslash\{b\})^{\ell}$ so that both $\mathbf{P}[E_{r}]$ and
$\mathbf{P}[F_{r}]$ will be of order $\tilde{O}(\ell^{-\frac{1}{2}})$
(See Theorem \ref{thm:Azuma-Hoeffding-inequality} of section \ref{sec:Common=000020entry}
below). Thus, the conditional probability to determine $\phi(X)_{0}$
by reading $X_{[-\ell,\ell]}$ is 
\begin{equation}
\mathbf{P}[R_{\phi}\leq\ell\boldsymbol{|}\ell(S_{r})=\ell]=1-\tilde{O}(\ell^{-\frac{1}{2}}).\label{eq:conditional=000020R}
\end{equation}
This is the very basic idea of why a coding as in Theorem \ref{thm:Main=000020Thm}
should exists. In order to make this idea into a real coding satisfying
the conclusion of Theorem \ref{thm:Main=000020Thm}, the first issue
to handle is how to make (\ref{eq:conditional=000020R}) to hold unconditionally,
and for all $\ell$. For that, we construct a marker process with
good control on its interval lengths $\ell(S_{r})$ for each $r$:
with this marker process, for any $n$ there is some $r$ s.t. the
likely length of $S_{r}$ is $\text{\ensuremath{\tilde{\Theta}(n)}. }$Thus,
for any $n$, for some $r$, the probability $\mathbf{P}[R_{\phi}\leq n]$
can be estimated by the left hand side of (\ref{eq:conditional=000020R}),
so that
\begin{align*}
\mathbf{P}[R_{\phi}\leq n] & \apprge\mathbf{P}[R_{\phi}\leq\ell(S_{r})\boldsymbol{|}n\geq\ell(S_{r})=\tilde{\Theta}(n)]=1-\tilde{O}(n^{-\frac{1}{2}}).
\end{align*}

The second issue to deal with is lifting this result to the general
case, where $\mathbf{p}$ and $\mathbf{q}$ do not share a common
entry. This assumption was used to split the two processes simultaneously
by the same marker process. One advantage of the simultaneous splitting
is that the splitting procedure costs the same amount of information
in the two processes, and leaves behind two independent processes
of the same entropy. This by itself is not crucial and could have
been overcome. For example, if we would consider homomorphisms, it
seems that one could split the source alone into words, while the
waste of information and the effect on the independence of the process
is negligible. The main advantage of the simultaneous splitting is
that it enables to restore back from the image the source sample in
a finitary manner, which is crucial for a map to be an isomorphism.
In the original work of Keane and Smorodinsky, the same problem was
resolved by taking an intermediate i.i.d. process $\left(C^{\mathbb{Z}},\mathbf{r}^{\mathbb{Z}}\right)$
with marginal $\mathbf{r}$ that share a common entry with both $\mathbf{p}$
and $\mathbf{q}$. Hence, once there is a finitary isomorphism between
the intermediate process and each one of the processes, a composition
of the isomorphisms gives a finitary isomorphism between the two processes
in question. We use this idea of an intermediate process, but that
alone would not resolve the problem, as in contrast to the finitary
property, the coding radii tail law is not preserved under composition
in general. The intermediate process is endowed with two collections
of ``words classes''. While each collection is used to construct
a different isomorphism with one of the processes in question, they
are heavily related to each other, in the sense that for a given sample
of the intermediate process, the knowledge of its class on some split
interval in one collection, determines its class on some split interval
of large length (w.r.t. the former interval length). This enables
us to compose the two isomorphisms without losing the coding radius
tail law. Although it may seem like a technical issue, this step of
lifting the common-entry-result to the general case was our major
challenge in this work.

\subsubsection*{Organization of the paper}

In the first section, we provide a marker process with good control
of the splitting lengths. The second section provides the main machinery
of the coding, and the main theorem is proved for the special case
of two aperiodic Markov processes sharing a state with equal distribution.
The third section deals with the general case, there the main theorem
is proved. 
\begin{acknowledgement*}
This paper is part of the author\textquoteright s Ph.D. thesis, conducted
under the guidance of Professor Michael Hochman, whom I would like
to thank for all his support and advice.
\end{acknowledgement*}

\section{\protect\label{sec:Marker-processes}Marker processes}

\subsection{\protect\label{subsec:Definitions-and-a}Definitions and a lemma}

During the paper we will use extensively expressions of the kind $n^{b}\log^{a}n$.
For cleaner formulations, we make the following convention: for any
$a,b\in\mathbb{R}$ define the function $\Lambda_{a,b}:[0,\infty)\rightarrow[0,\infty)$
by
\begin{equation}
\Lambda_{a,b}(x)=x^{b}\log^{a}x\label{eq:Lambda_c}
\end{equation}
and in case that $b=1$ we abbreviate $\Lambda_{a}(x)=\Lambda_{a,1}(x)$.
In practice, we will make use of these functions only with $b\in\{\pm1,\pm\frac{1}{2}\}$.

We use the notation $\tilde{O}$ (``soft $O$ notation'') defined
for functions $f,g$ on $[0,\infty)$ by
\[
f(x)=\tilde{O}(g(x))\iff\exists C,k>0\text{ s.t. }\forall x,\;|f(x)|\leq C|\Lambda_{k}(g(x))|
\]
and similarly $\tilde{\Omega}$ is defined by:
\[
f(x)=\tilde{\Omega}(g(x))\iff\exists C>0,k\in\mathbb{R}\text{ s.t. }\forall x,\;f(x)\geq C|\Lambda_{k}(g(x))|
\]
By $f(x)=\tilde{\Theta}(g(x))$ we mean that both $f(x)=\tilde{O}(g(x))$
and $f(x)=\tilde{\Omega}(g(x))$.

We devote a notion $\Delta(n)$ for superpolynomial decay:
\[
f(n)=\Delta(n)\iff\forall k\geq1,f(n)=O(n^{-k}).
\]

Occasionally, we will use Landau notation for both describing an event
and estimating its probability at the same time. To avoid inaccuracy,
we stress here the meaning of such statements: Suppose $X(n)$ and
$Y(n)$ are random functions ($Y(n)$ might be deterministic). The
``event''
\[
\{X(n)=\Theta(Y(n))\}
\]
refers to any choice of functions $F_{1},F_{2}:[0,\infty)\rightarrow[0,\infty)$
of order $\Theta(n)$, and to the corresponding sequence of events
\[
E_{n}=\{F_{1}(f(Y(n)))\leq X(n)\leq F_{2}(f(Y(n)))\}.
\]
Accordingly, the assertion
\[
\mathbf{P}[X(n)=\Theta(f(Y(n)))]=O(g(n))
\]
means there exist functions $F_{1},F_{2}$ on $[0,\infty)$ both satisfying
$F_{i}(n)=Q_{i}(f(n))$ for some $Q_{1},Q_{2}>0$, s.t. 
\[
\mathbf{P}[F_{1}(Y(n))\leq X(n)\leq F_{2}(Y(n))]=Q'(g(n)).
\]

For a set of symbols $A$ and a set of integers $I\subset\mathbb{Z}$,
we write $A^{I}$ for the space of sequences with symbols in $A$
labeled by $I$:
\[
A^{I}=\left\{ \left(a_{i}\right)_{i\in I}:a_{i}\in A\right\} 
\]

\begin{defn}
\label{def:Let--be}Let $X,Y$ be two random variables on the same
probability space $\Omega$. The event \textit{``$X$ determines
$Y$''} (denoted by $\{X\Rrightarrow Y\}$) refers to the set of
all $\omega\in\Omega$ for whom the conditioned random variable $Y|(X=X(\omega))$
is deterministic:
\[
\{X\Rrightarrow Y\}=\{\omega:\,\forall\omega',\,X(\omega')=X(\omega)\:\Rightarrow\:Y(\omega')=Y(\omega)\}
\]
\end{defn}

Two elementary properties of this notion will be of use during the
proofs of this paper:

1. Measurability: $\{X\Rrightarrow Y\}\cap\varSigma(Y)\subset\varSigma(X)$
and in particular $\{X\Rrightarrow Y\}\in\varSigma(X).$ (here $\varSigma(X)$
(and $\varSigma(Y)$) denotes the minimal sigma-algebra under which
$X$ (resp. $Y$) is measurable.)

2. Transitivity: $\{X\Rrightarrow Y\}\cap\{Y\Rrightarrow Z\}\subset\{X\Rrightarrow Z\}$.

With this notation, given a stationary map $\phi$ on a process $X$,
we can reformulate its coding radius $R_{\phi}$ as
\[
R_{\phi}(x)=\min\{n:x\in\{X_{[-n,n]}\Rrightarrow\phi(X)_{0}\}\}
\]

Let $X\sim(A^{\mathbb{Z}},\mu)$. A state $a\in A$ is called a \textit{renewal
state} if conditioned on $X_{0}=a$, the random samples $X_{[-n,-1]}$
and $X_{[1,n]}$ are independent for all $n$. 
\begin{defn}
\label{def:The-distribution-of}The \textbf{\textit{distribution of
}}\textit{$\boldsymbol{a}$} is the process $\widehat{X}^{a}=\left(\mathbf{1}\{X_{n}=a\}\right)_{n\in\mathbb{Z}}$.
When there is no ambiguity about the choice of $a$, we abbreviate
$\widehat{X}=\widehat{X}^{a}$.
\begin{defn}
A \textbf{\textit{marker process}} (or a \textit{marker factor}) is
a stationary map of the distribution of $a$, $\xi=\xi(\widehat{X})\in(\mathbb{N}\cup\{0\})^{\mathbb{Z}}$
that outputs a sequence of non-negative integers, so that
\end{defn}

\end{defn}

\begin{itemize}
\item $\xi(\widehat{X})_{0}>0\Rightarrow\widehat{X}_{0}=1$.
\item For any integer $r>0$, $\mathbf{P}[\xi_{0}\geq r]>0$.
\end{itemize}
\begin{defn}
A \textbf{\textit{skeleton}} for $\xi(\widehat{X})$ is any tuple
$S=(\widehat{x}_{[\alpha,\beta]},\xi(\widehat{x})_{[\alpha,\beta]})$
that its indices $\alpha<\beta$ satisfy 
\[
\xi_{\alpha},\xi_{\beta+1}>\max\left\{ \xi_{\alpha+1},...,\xi_{\beta}\right\} .
\]
(the evaluated point $\widehat{x}$ is omitted for brevity.) In that
case, if $I\subset[\alpha,\beta]$, we say that \textit{$S$ covers
$x$ at $I$}. Often we treat the skeleton $S$ as the set of coordinates
it covers, $[\alpha,\beta]$. Thus for example, $i\in S$ means that
$S$ covers $i$, $x_{S}$ refers to the name $x_{[\alpha,\beta]}$,
and $|S|$ refers to the size $\beta-\alpha+1$ of $[\alpha,\beta]$.
$S$ is a \textit{skeleton of rank $r$}, or just an \textit{$r$-skeleton},
if
\[
\xi_{\alpha},\xi_{\beta+1}\geq r>\max\left\{ \xi_{\alpha+1},...,\xi_{\beta}\right\} .
\]
the minimal $r$ for which the above holds is the minimal rank of
the skeleton, denoted by $\text{mr}(S)$. Any $r$-skeleton $S$ can
be uniquely decomposed as the concatenation of its $(r-1)$-sub-skeletons:
\[
S=S^{1}\cdots S^{m}.
\]
If $\text{mr}(S)=r>0$ then all skeletons $S^{1},...,S^{m}$ are \textit{proper}
sub-skeletons, i.e. non of them equals $S$.
\end{defn}

The marker process $\xi$ serves as representing a sequence $(\mathcal{P}_{r}(X))_{r=1}^{\infty}$
of coarser and coarser partitions of the infinite word $X=(X_{n})_{n\in\mathbb{Z}}$
into finite words - the $r$th partition $\mathcal{P}_{r}(X)$ consists
of the words $\{X_{S}\}_{S}$, where $S$ ranges over all skeletons
of minimal rank $r$ that cover $X$. Accordingly, in order to represent
only a subsequence of the partitions, say, $(\mathcal{P}_{r_{n}}(X))_{n=1}^{\infty}$
for some increasing sequence $(r_{n})\subset\mathbb{N}$, one can
\textit{thin out $\xi$} to get another marker process $\bar{\xi}$
by the following rule: For any $o\in\mathbb{Z}$, let
\begin{equation}
\bar{\xi}_{o}=\begin{cases}
0 & \xi_{o}=0\\
\max\{n:r_{n}\leq\xi_{o}\} & \xi_{o}\neq0
\end{cases}\label{eq:-9}
\end{equation}

Given a sample $x\in A^{\mathbb{Z}}$ and an integer $o\in\mathbb{Z}$,
we denote by $S_{r}(o,x)$ the skeleton with minimal rank $r$ that
covers $x$ at $o$. When the integer $o$ is zero, we simply write
$S_{r}(x):=S_{r}(0,x)$, and when the sample is the random variable
$X$, we abbreviate $S_{r}:=S_{r}(0,X)$.

Let $X\sim(A^{\mathbb{Z}},\mu)$ be an aperiodic Markov process of
memory 1 on a finite state space $A$. That means:
\begin{itemize}
\item The state space $A$ is finite, $|A|<\infty$.
\item $X$ is of memory 1, that is, for any sample $x\in A^{(-\infty,0]}$,
$\mathbf{P}[X_{0}=x_{0}|X_{(-\infty,-1]}=x_{(-\infty,-1]}]=\mathbf{P}[X_{0}=x_{0}|X_{-1}=x_{-1}]$
(all its states are renewal).
\item For any $a',a''\in A$ there is some $n$ s.t. $\mathbf{P}[X_{0}=a',X_{n}=a'']>0$
($X$ is irreducible).
\item For any $a',a''\in A$, the greatest common divisor of the set $\{n:\mathbf{P}[X_{0}=a',X_{n}=a'']>0\}$
equals 1 ($X$ has period 1).
\end{itemize}
The aim of this section is to provide, for any such $X$ with any
of its states $a\in A$, a marker factor $\xi=\xi(\widehat{X}^{a})$
that exhibits a good control on its skeletons sizes, as well as on
the radius needed to encode them. This is stated precisely in the
following lemma, referred to as the Marker process Lemma:
\begin{lem}
\label{lem:Marker=000020process}(Marker process Lemma) Let $X$ be
an aperiodic Markov process on finite state space, and let $a\in A$
be any of its states. There exists a marker factor $\xi=\xi(\widehat{X}^{a})$
with the property that for some $c>0$, for any $n$ there is some
$r=r(n)$ so that
\[
\begin{gathered}\mathbf{P}[X_{[-n,n]}\Rrightarrow S_{r}(X)]=1-\Delta(n)\\
and\\
\mathbf{P}[\frac{n}{\log^{c}n}\leq|S_{r}(X)|\leq n]=1-\Delta(n).
\end{gathered}
\]
\end{lem}

By thinning out the marker process provided by the above lemma (as
in (\ref{eq:-9})), one gets the following result:
\begin{cor}
\label{cor:Marker=000020process}For any large enough $L_{1},c>0$
and any sequence $(L_{r})_{r=2}^{\infty}$ satisfying $\frac{L_{r}}{\log^{c}L_{r}}>L_{r-1}$
for each $r$, there exists a marker factor $\xi=\xi(\widehat{X}^{a})$
satisfying:
\begin{align}
\mathbf{P}[X_{[-L_{r},L_{r}]} & \Rrightarrow S_{r}(X)]=1-\Delta(L_{r})\nonumber \\
\mathbf{P}[\frac{L_{r}}{\log^{c}L_{r}}\leq & |S_{r}(X)|\leq L_{r}]=1-\Delta(L_{r})\label{eq:lengths=000020sequence=0000201}
\end{align}
\end{cor}

\begin{defn}
\label{def:nice=000020marker}A marker factor $\xi$ is said to be
\textit{nice} and a sequence $(L_{r})_{r\geq0}$ is called a \textit{lengths
sequence }for $\xi$, if $\xi$ satisfies the conclusion of the Marker
process lemma (Corollary \ref{cor:Marker=000020process}) w.r.t. $(L_{r})_{r\geq0}$. 

In the next subsection we present a construction of a marker process,
which is shown in Section \ref{subsec:Proof-of-lemma} to satisfy
the conclusions claimed in the above lemma.
\end{defn}

\subsection{\protect\label{subsec:Construction-of-alpha}Construction of the
marker process $\xi$}

\subsubsection*{Outline of the construction of $\xi$}

Fix some $c>0$ and increasing ``lengths'' $(L_{r})_{r\geq0}$ with
$L_{r}\leq\Lambda_{-c}(L_{r+1})$ (recall the notation $\Lambda_{-c}$
as in (\ref{eq:Lambda_c})). We first find a sequence of events $(E_{r})_{r\geq0}$
, so that for any $r$, the time 
\[
\tau_{E_{r}}=\min\left\{ i>0:\sigma^{i}X\in E_{r}\right\} -\max\left\{ i\leq0:\sigma^{i}X\in E_{r}\right\} 
\]
 that passed since the last occurrence of translations of $X$ in
$E_{r}$ up to $0$, till its first occurrence in is close to $L_{r}$,
with high probability:
\[
\mathbf{P}[\Lambda_{-c}L_{r}\leq\tau_{E_{r}}\leq\Lambda_{c}L_{r}]=1-o_{r\rightarrow\infty}(1).
\]
Naively, we would like to treat the occurrence times of translations
of $X$ in $E_{r}$ as place-marks where the sequence $X$ breaks
into $r$-skeletons, as this provides skeletons that are likely to
have sizes as in the Marker process lemma. However, unless $E_{r}\subset E_{r-1}$
for all $r$, this procedure will yield $r$-skeletons that are not
always being decomposed into $r-1$-skeletons. Thus, we slightly modify
this partition of $X$ by adding to each $r$-skeleton the $r-1$-skeleton
that overlaps with its left boundary, while removing the $r-1$-skeleton
that overlaps with its right boundary (this can be viewed as replacing
inductively the event $E_{r}$ by a modified event $\bar{E}_{r}$
so that the sequence $(\bar{E_{r}})_{r}$ is decreasing). In order
to do this without causing a significant change to the size of the
$r$-skeletons, the events $(E_{r})_{r}$ are chosen in advance to
have return times concentrated far enough away from each other, so
that adding to (or subtracting from) a candidate $r$-skeleton some
candidate $r-1$-skeleton, causes negligible effect on the former
skeleton size. Also, in case an occurrence time of $E_{r}$ is much
larger than its expected value, then instead of defining a corresponding
$r$-skeleton of that size, we define all ($r-1)$-skeletons there
to be $r$-skeletons too. This guarantees that adding overlapping
skeletons of lower rank to a skeleton of high rank doesn't significantly
change the size of the higher-rank skeleton, regardless of whether
the lower rank skeletons are ``likely'' ones or not.$\hfill\square$

Given $X$ and $a$, recall the distribution of $a$ $\widehat{X}=\left(\mathbf{1}\{X_{n}=a\}\right)_{n\in\mathbb{Z}}$
(Def. \ref{def:The-distribution-of}). Let 
\[
k=\min\left\{ k'\geq0:\mathbf{P}\left[X_{0}=a,X_{k'}=a\right]>0\right\} 
\]
and let $w=0^{k-1}1$. Because 1 is a renewal state for $\widehat{X}$,
$w^{\ell}$ appears in $\widehat{X}_{[1,\ell|w|]}$ with positive
probability for every $\ell$, and a generic sample of $\widehat{X}$
can be uniquely written as
\[
\begin{array}{cr}
\hspace{2cm}\cdots w^{\ell_{1}}0^{k_{1}}w^{\ell_{2}}0^{k_{2}}\cdots & \hspace{1cm}(k_{i}\geq k)\end{array}
\]
(by aperiodicity of $X$ and the choice of $w$, $0w$ must also appear
in $\widehat{X}_{[0,|w|]}$ with positive probability.)

Fixing an increasing sequence of positive integers $(N_{r})_{r\geq0}$
that will be specified later (see (\ref{eq:-3}) below), we define
for each $r$ the event
\begin{align}
D_{r} & =\left\{ \widehat{X}_{[-|w|,(N_{r}-1)|w|]}=0w^{N_{r}}\right\} .\label{eq:-58}
\end{align}

Define the indicator process $\mathbf{1}_{D_{r}}(X)$ by
\[
\mathbf{1}_{D_{r}}(X)_{n}=\mathbf{1}\{\sigma^{n}X\in D_{r}\}.
\]
(here $\sigma:A^{\mathbb{Z}}\rightarrow A^{\mathbb{Z}}$ is the left
shift operator, defined by $(x_{n})_{n\in\mathbb{Z}}\overset{{\scriptscriptstyle \sigma}}{\longmapsto}(x_{n+1})_{n\in\mathbb{Z}}$).
Accordingly, define the ``lengths'' sequence 
\begin{equation}
L_{r}=\mathbf{P}\left[D_{r}\right]^{-1}\label{eq:Def.=000020of=000020L_r}
\end{equation}
which can be viewed as the expected return times of the words $0w^{N_{r}}$
in $\widehat{X}$. For each $r\geq0$ define the event
\begin{equation}
E_{r}=\left\{ \mathbf{1}_{D_{r}}(X)_{[0,\Lambda_{-1}(L_{r})]}=10^{\Lambda_{-1}(L_{r})}\right\} .\label{eq:-1}
\end{equation}
These events will serve as mark-places to define the skeletons of
$\xi$. We first show that the length of the word between the first
positive occurrence and the first non positive occurrence of $E_{r}$
is concentrated in the interval $[\Lambda_{-1}L_{r},\Lambda_{2}L_{r}]$.
For an event $E$, an element $x\in A^{\mathbb{Z}}$ and an index
$n\in\mathbb{Z}$, write 
\begin{align*}
\tau_{E}^{>n}(x) & =\min\{i>n:\sigma^{i}x\in E\}\\
\tau_{E}^{<n}(x) & =\max\{i<n:\sigma^{i}x\in E\}
\end{align*}
 and similarly, define $\tau_{E}^{\geq n}(x):=\tau_{E}^{>n-1}(x)$
and $\tau_{E}^{\leq n}(x):=\tau_{E}^{<n+1}(x)$. 
\begin{prop}
\label{prop:return=000020times=000020concentration}For any $d>2$,
\[
\mathbf{P}\left[\tau_{E_{r}}^{>0}-\tau_{E_{r}}^{\leq0}\in[\Lambda_{-1}L_{r},\Lambda_{d}L_{r}]\right]=1-\Delta(L_{r}).
\]
\end{prop}

\begin{proof}
By the definition of $E_{r}$ (\ref{eq:-1}), the lower bound $\tau_{E_{r}}^{>0}-\tau_{E_{r}}^{\leq0}\geq\Lambda_{-1}L_{r}$
always holds. Let $\tau_{E_{r}}^{\leq0}=k_{0}<k_{1}<\cdots<k_{t+1}=\tau_{E_{r}}^{>0}$
be all indices in the range $[\tau_{E_{r}}^{\leq0},\tau_{E_{r}}^{>0}]$
for which $\mathbf{1}_{D_{r}}(X)$ gives 1 (here $t$ stands for the
number occurrences of $D_{r}$ after $\tau_{E_{r}}^{\leq0}$ and before
$\tau_{E_{r}}^{>0}$). We have 
\[
k_{1}-k_{0}>\Lambda_{-1}L_{r}\geq\max_{1\leq i\leq t}\{k_{i+1}-k_{i}\}.
\]
Notice that
\begin{align}
\{\tau_{E_{r}}^{>0}-\tau_{E_{r}}^{\leq0}>\Lambda_{d}L_{r}\} & \subset\{t>\log L_{r}\}\,\cup\,\{k_{1}-k_{0}>\frac{1}{2}\Lambda_{d}L_{r}\}\backslash\{t>\log L_{r}\}\label{eq:-16}
\end{align}
thus bounding both $\mathbf{P}[t>\log L_{r}]$ and $\mathbf{P}[k_{1}-k_{0}>\Lambda_{d}L_{r}-L_{r}]$
by $\Delta(L_{r})$ will imply the claim of the proposition. We split
the proof of these bounds, and they will be proven in Claims \ref{prop:asdfasdf}
and \ref{claim:In-the-same} below.
\end{proof}
\begin{claim}
\label{prop:asdfasdf}In the same circumstance as in (\ref{eq:-16}),
one has
\[
\mathbf{P}[t>\log L_{r}]=\Delta(L_{r}).
\]
\end{claim}

We give first a rough idea of the claim: From the definition of $E_{r}$,
for any $1\leq i<t$ we have that 
\begin{equation}
k_{i+1}-k_{i}<\Lambda_{-1}L_{r}.\label{eq:-27}
\end{equation}
By the definition of $D_{r}$, having a renewal state at 0, the event
(\ref{eq:-27}) is independent of $X_{(-\infty,k_{i}-1]}$, thus a
union bound on it will give 
\begin{align*}
\mathbf{P}[t>\log L_{r}] & \apprle\mathbf{P}[k_{i+1}-k_{i}<\Lambda_{-1}L_{r}.]^{\log L_{r}}\\
 & \lesssim\left(\log L_{r}\cdot\Lambda_{-1}L_{r}\cdot\mathbf{P}[X\in D_{r}]\right)^{\log L_{r}}\\
 & =\left(\log^{-1}L_{r}\right)^{\log L_{r}}.
\end{align*}

We now give the detailed proof for the bound:

\textit{Proof of Claim \ref{prop:asdfasdf}.} We have 
\begin{equation}
\mathbf{P}[t>\log L_{r}]\leq\mathbf{P}[\sum_{n=\tau_{D_{r}}^{>0}}^{\tau_{E_{r}}^{>0}}\mathbf{1}_{D_{r}}(X)_{n}>\frac{1}{2}\log L_{r}]+\mathbf{P}[\sum_{n=\tau_{E_{r}}^{\leq0}}^{\tau_{D_{r}}^{\leq0}}\mathbf{1}_{D_{r}}(X)_{n}>\frac{1}{2}\log L_{r}]\label{eq:-24}
\end{equation}
The first summand can be written as
\begin{multline}
\mathbf{P}[\sum_{n=\tau_{D_{r}}^{>0}}^{\tau_{E_{r}}^{>0}}\mathbf{1}_{D_{r}}(X)_{n}>\frac{1}{2}\log L_{r}]=\\
=\mathbf{P}[\sigma^{\tau_{D_{r}}^{>0}}X\notin E_{r}]\cdot\mathbf{P}[\sum_{n=\tau_{D_{r}}^{>{\scriptscriptstyle \tau_{D_{r}}^{>0}}}}^{\tau_{E_{r}}^{>0}}\mathbf{1}_{D_{r}}(X)_{n}>\frac{1}{2}\log L_{r}-1|\sigma^{\tau_{D_{r}}^{>0}}X\notin E_{r}].\label{eq:}
\end{multline}
and since $\{X\in D_{r}\}\subset\{X_{0}=1\}$ and 1 is a renewal state,
the first expression in the right hand product can be written as:
\begin{align*}
\mathbf{P}[\sigma^{\tau_{D_{r}}^{>0}}X\notin E_{r}] & =\sum_{n}\mathbf{P}[\sigma^{n}X\notin E_{r}|\tau_{D_{r}}^{>0}=n]\mathbf{P}[\tau_{D_{r}}^{>0}=n]\\
 & =\sum_{n}\mathbf{P}[X\notin E_{r}|X\in D_{r},\forall-n<i<0,\sigma^{i}X\notin D_{r}]\mathbf{P}[\tau_{D_{r}}^{>0}=n]\\
 & =\sum_{n}\mathbf{P}[X\notin E_{r}|X\in D_{r}]\mathbf{P}[\tau_{D_{r}}^{>0}=n]\\
 & =\mathbf{P}[X\notin E_{r}|X\in D_{r}]
\end{align*}

By the same argument, we can write the second expression in (\ref{eq:})
as
\[
\mathbf{P}[\sum_{n=\tau_{D_{r}}^{>{\scriptscriptstyle \tau_{D_{r}}^{>0}}}}^{\tau_{E_{r}}^{>0}}\mathbf{1}_{D_{r}}(X)_{n}>\frac{1}{2}\log L_{r}-1|\sigma^{\tau_{D_{r}}^{>0}}X\notin E_{r}]=\mathbf{P}[\sum_{n=\tau_{D_{r}}^{>0}}^{\tau_{E_{r}}^{>0}}\mathbf{1}_{D_{r}}(X)_{n}>\frac{1}{2}\log L_{r}-1|X\in D_{r}]
\]
and by induction we get 
\[
\mathbf{P}[\sum_{n=\tau_{D_{r}}^{>0}}^{\tau_{E_{r}}^{>0}}\mathbf{1}_{D_{r}}(X)_{n}>\frac{1}{2}\log L_{r}]=\mathbf{P}[X\notin E_{r}|X\in D_{r}]^{\frac{1}{2}\log L_{r}}.
\]
A similar idea can be applied to the second summand in the right hand
side of (\ref{eq:-24}) which gives us
\begin{align}
\mathbf{P}[t>\log L_{r}] & =2\mathbf{P}[X\notin E_{r}|X\in D_{r}]^{\frac{1}{2}\log L_{r}}\label{eq:-15}
\end{align}
Conditioned on the event that $X\in D_{r}$, the probability that
$X\notin E_{r}$ can be estimated by
\begin{align*}
\mathbf{P}[X\notin E_{r}|X\in D_{r}] & =\mathbf{P}[\exists1\leq n\leq\Lambda_{-1}L_{r}\,\mathbf{1}_{D_{r}}(X)_{n}=1|\mathbf{1}_{D_{r}}(X)_{0}=1]\\
 & \leq\sum_{n=1}^{\Lambda_{-1}L_{r}}\mathbf{P}[\sigma^{n}X\in D_{r}|X_{(N_{r}-1)|w|}=1]\\
 & \frac{1}{\mathbf{P}[X_{(N_{r}-1)|w|}=1]}\sum_{n=1}^{\Lambda_{-1}L_{r}}\mathbf{P}[\sigma^{n}X\in D_{r}]\\
 & =O(\Lambda_{-1}L_{r}\mathbf{P}[X\in D_{r}])\\
 & =O(\log^{-1}L_{r})
\end{align*}
and plugging that into (\ref{eq:-15}) yields
\[
\mathbf{P}[t>\log L_{r}]\leq L_{r}^{-\Omega(\log\log L_{r})}.
\]
$\hfill\square$

We now consider the second event in the right hand side union of
(\ref{eq:-16}).
\begin{claim}
\label{claim:In-the-same}In the same circumstance as in (\ref{eq:-16}),
one has
\[
\mathbf{P}[\{k_{1}-k_{0}>\Lambda_{d}L_{r}\}\backslash\{t>\log L_{r}\}]=\Delta(L_{r}).
\]
\end{claim}

\begin{proof}
Notice that on this event, 
\begin{equation}
[-L_{r},0]\cap[k_{0},k_{1}]\neq\emptyset\label{eq:-28}
\end{equation}
This follows from the fact that both $k_{0}\leq0$ and $-L_{r}<k_{1}$
hold: The first inequality is immediate from the definition of $k_{0}=\tau_{E_{r}}^{\leq0}$.
To see why the second inequality holds, observe that
\[
0<\tau_{E_{r}}^{>0}\leq k_{1}+t\Lambda_{-1}L_{r}\leq k_{1}+L_{r}
\]
(we used here the assumption $t\leq\log L_{r}$). Any index $i$ in
the intersection (\ref{eq:-28}) must precede or succeed a sequence
of $\lceil\frac{1}{2}\Lambda_{d}L_{r}\rceil$ zeros in the $\mathbf{1}_{D_{r}}(X)$-name,
according to whether $i$ sits in the first half or the second half
of the interval $[k_{0},k_{1}]$ respectively:
\begin{align}
\{k_{1}-k_{0}>\Lambda_{d}L_{r}\}\backslash\{t>\log L_{r}\}\subset & \bigcup_{i=-L_{r}}^{0}\{\mathbf{1}_{D_{r}}(X)_{[i+1,i+\lceil\frac{1}{2}\Lambda_{d}L_{r}\rceil]}=0^{\lceil\frac{1}{2}\Lambda_{d}L_{r}\rceil}\}\label{eq:return=000020times=000020concentration=0000201}\\
 & \cup\bigcup_{i=-L_{r}}^{0}\{\mathbf{1}_{D_{r}}(X)_{[i-\lceil\frac{1}{2}\Lambda_{d}L_{r}\rceil,i-1]}=0^{\lceil\frac{1}{2}\Lambda_{d}L_{r}\rceil}\}\nonumber 
\end{align}

Let $U\subset A^{[0,|w|N_{r}]}$ be the collection of words $u$ satisfying
$\hat{u}=0w^{N_{r}}$ and $\mathbf{P}[X_{[0,|w|N_{r}]}=u]>0$ (here
$\hat{u}$ refers to the $a$-distribution name of $u$, see Def.
\ref{def:The-distribution-of}). As $X$ is aperiodic, there are some
$n_{0}>0$ and $\eta>0$ s.t. for any $a',a''\in A$
\[
\mathbf{P}[X_{0}=a''|X_{-n_{0}}=a']\geq\eta.
\]
Since $X$ is of memory-one, we have
\[
\begin{array}{l}
\min_{a'\in A}\mathbf{P}[X\in D_{r}|X_{-|w|-n_{0}}=a']=\min_{a'\in A}\mathbf{P}[X\in U|X_{-|w|-n_{0}}=a']\\
=\min_{a'\in A}\sum_{u\in U}\mathbf{P}[X_{[1,|w|N_{r}]}=u_{[1,|w|N_{r}]}|X_{0}=u_{0}]\cdot\mathbf{P}[X_{0}=u_{0}|X_{-n_{0}}=a']\\
\geq\eta\sum_{u\in U}\mathbf{P}[X_{[1,|w|N_{r}]}=u_{[1,|w|N_{r}]}|X_{0}=u_{0}]\\
\geq\eta\mathbf{P}[X\in D_{r}]\\
=\eta L_{r}^{-1}
\end{array}
\]
 hence for all $n$,
\begin{align}
\mathbf{P}[\mathbf{1}_{D_{r}}(X)_{[1,n]}=0^{n}] & \leq\mathbf{P}[\forall0\leq m\leq\lfloor\frac{n}{|w|N_{r}+n_{0}}\rfloor,\,\mathbf{1}_{D_{r}}(X)_{m(|w|N_{r}+n_{0})}=0]\label{eq:-17}\\
 & \leq\left(1-\eta L_{r}^{-1}\right)^{\lfloor\frac{n}{|w|N_{r}+n_{0}}\rfloor}.\nonumber 
\end{align}
Notice that $L_{r}$ (\ref{eq:Def.=000020of=000020L_r}) can be written
as
\[
L_{r}=\mathbf{P}\left[D_{r}\right]^{-1}=\left(\mathbf{P}[\widehat{X}_{[-|w|,0]}=0w]\cdot\mathbf{P}[\widehat{X}_{[1,|w|]}=w|\widehat{X}_{0}=1]^{N_{r}-1}\right)^{-1}
\]
and so, by taking $d'\in(1,d-1)$ and $n=\Lambda_{d'+1}L_{r}\leq\lceil\frac{1}{2}\Lambda_{d}L_{r}\rceil$,
we have for some $C>0$
\begin{align}
\lfloor\frac{\lceil\frac{1}{2}\Lambda_{d}L_{r}\rceil}{|w|N_{r}+n_{0}}\rfloor & \geq\frac{\Lambda_{d'+1}L_{r}}{|w|N_{r}+n_{0}}\geq C\Lambda_{d'}L_{r}\log L_{r}\label{eq:-29}\\
\end{align}
Replacing $d'$ again by a smaller constant $d''\in(1,d')$ will yield
\begin{align*}
\mathbf{P}[\mathbf{1}_{D_{r}}(X)_{[1,\lceil\frac{1}{2}\Lambda_{d}L_{r}\rceil]}=0^{\lceil\frac{1}{2}\Lambda_{d}L_{r}\rceil}] & \leq\left(1-\eta L_{r}^{-1}\right)^{CL_{r}\log^{d''}L_{r}}\\
 & \leq L_{r}^{-\frac{C}{\eta}\log^{d''}L_{r}}\\
 & =\Delta(L_{r})
\end{align*}
 Thus the probability of the left side event in (\ref{eq:return=000020times=000020concentration=0000201})
is bounded by
\begin{align*}
\mathbf{P}[(k_{1}-k_{0}>\Lambda_{d}L_{r})\,\wedge\,(t\leq\log L_{r})] & \leq2(L_{r}+1)\mathbf{P}[\mathbf{1}_{D_{r}}(X)_{[1,\lceil\frac{1}{2}\Lambda_{d}L_{r}\rceil]}=0^{\lceil\frac{1}{2}\Lambda_{d}L_{r}\rceil}]\\
 & =\Delta(L_{r})
\end{align*}
 and the claim of the proposition follows.
\end{proof}
Following the last proposition, we fix $d>2$, so that 
\begin{equation}
\mathbf{P}[\tau_{E_{r}}^{>0}-\tau_{E_{r}}^{\leq0}>\Lambda_{d}L_{r}]=\Delta(L_{r}).\label{eq:-4}
\end{equation}
The marker process $\xi$ will be constructed so that its random $r$-skeleton
$S_{r}$ is likely to consist of all coordinates $[\tau_{E_{r}}^{\leq0},\tau_{E_{r}}^{>0}]$
up to addition or subtraction of the $(r-1)$-skeletons that overlap
its boundary. By taking $L_{r}\gg L_{r-1}$ we guarantee that this
addition/subtraction does not change significantly the size of the
$r$-skeleton. On the other hand, for $\xi$ to provide for each $n$
an $r=r(n)$ s.t. $|S_{r}(X)|=\tilde{\Theta}(n)$ with high probability,
the sequence $(L_{r})_{r\geq0}$ should satisfy $L_{r}=\tilde{O}(L_{r-1})$
for all $r$.

Accordingly, we fix $N_{1}$ to be some large positive integer\textcolor{blue}{{}
}(so that (\ref{eq:-2}) and (\ref{eq:-6}) will take place), and
define $\left(N_{r}\right)_{r\geq0}$ inductively by
\begin{equation}
N_{r+1}=N_{r}+\lfloor(d+2)\frac{\log N_{r}}{\log(\mathbf{P}[\widehat{X}_{[1,|w|]}=w|\widehat{X}_{0}=1]^{-1})}\rfloor.\label{eq:-3}
\end{equation}

This together with the definition of $L_{r}$ (\ref{eq:Def.=000020of=000020L_r})
gives
\begin{equation}
L_{r+1}\sim\frac{\Lambda_{d+2}L_{r}}{\log^{d+2}(\mathbf{P}[\widehat{X}_{[1,|w|]}=w|\widehat{X}_{0}=1]^{-1})}\label{eq:L_r+1=000020similar=000020to=000020L_r}
\end{equation}
Assuming $N_{1}$ is large enough, the following holds for all $r$:
\begin{prop}
\label{prop:(Properties-of-L_r)} For any $r$,
\begin{equation}
\Lambda_{-2}L_{r+1}<\Lambda_{d}L_{r}<\Lambda_{-1.9}L_{r+1}.\label{eq:-2}
\end{equation}
\end{prop}

We are now ready to define the marker process $\xi$. Define maps
$\left\{ \xi_{r}:A^{\mathbb{Z}}\rightarrow\mathbb{N}_{0}^{\mathbb{Z}}\right\} _{r\geq1}$
inductively on $r$ as follows: For $r=1$, define $\left(\xi_{1}(x)_{n}\right)_{n\in\mathbb{Z}}$
by
\[
\xi_{1}(x)_{n}=\mathbf{1}_{D_{1}}(x)_{n}
\]
thus $1$-skeletons are blocks of $x$ that begin with the first appearance
of $D_{1}$ and end right before the next appearance. Now let $r>1$,
and suppose we already defined $\xi_{r-1}(x)$ with symbols in$\{0,...,r-1\}$.
Fix a coordinate $n\in\mathbb{Z},$ and let 
\begin{align*}
\beta & =\min\left\{ i>n:\xi_{r-1}(x)_{i}=r-1\right\} \\
\alpha & =\max\left\{ i\leq\tau_{E_{r}}^{<\beta}(x):\xi_{r-1}(x)_{i}=r-1\right\} 
\end{align*}
First one should notice that
\[
\alpha\leq n<\beta
\]
This follows easily from the fact that $\alpha<\beta$, together with
the definition of $\beta$ that guarantees 
\[
\forall n<i<\beta,\xi_{r-1}(x)_{i}\neq r-1,
\]
whilst $\xi_{r-1}(x)_{\alpha}=r-1$.

Define $\xi_{r}(x)_{n}$ as follows: If
\begin{equation}
\tau_{E_{r}}^{\geq\beta}(x)-\tau_{E_{r}}^{<\beta}(x)\leq\Lambda_{d}L_{r}\label{eq:r-skeletonCond1}
\end{equation}
then define
\[
\xi_{r}(x)_{n}=\begin{cases}
r & n=\alpha\\
\xi_{r-1}(x)_{n} & n>\alpha
\end{cases}
\]
Otherwise, define
\[
\xi_{r}(x)_{n}=\begin{cases}
r & \xi_{r-1}(x)_{n}=r-1\\
\xi_{r-1}(x)_{n} & \xi_{r-1}(x)_{n}<r-1
\end{cases}.
\]

Below is a scheme of using the samples of $x\in E_{r}$ and $\xi_{r-1}(x)$
in the first and second rows, to determine the sample $\xi_{r}(x)$
in the third row, where having the symbol $\mathtt{r}$ or $\mathtt{x}$
in the $n$th coordinate stands for $r$ or $\xi_{r-1}(x)_{n}$ respectively:
\[
\begin{array}{rl}
\mathbf{1}_{E_{r}} & \ldots\overbrace{\mathtt{10000000000\,0}}^{\text{length}\leq\Lambda_{d}L_{r}}\overbrace{\mathtt{10000\,00000\,00}}^{\text{length}>\Lambda_{d}L_{r}}\overbrace{\mathtt{1\,00000000000}}^{\text{irelevant}}\overbrace{\mathtt{1000000\ldots}}^{\text{length}\leq\Lambda_{d}L_{r}}\\
\mathbf{1}\{\xi_{r-1}\!=\!r\!-\!1\} & \ldots\underbrace{\mathtt{10010010001}}\underbrace{\mathtt{1\,00000\,10000\,10\,0}}\underbrace{\mathtt{100000000000\,000100\ldots}}\\
\xi_{r} & \mathtt{\ldots\underbrace{\mathtt{rxxxxxxxxxx}}\underbrace{\mathtt{r\,xxxxx}}\underbrace{\mathtt{rxxxx}}\underbrace{\mathtt{rxx}}\underbrace{\mathtt{\mathtt{rxxxxxxx}xxxx\,xxxxxx\ldots}}}
\end{array}
\]

For any $n$ and $x$, the sequence $\left(\xi_{r}(x)_{n}\right)_{r=1}^{\infty}$
is monotonically increasing, hence the limit map $\xi:A^{\mathbb{Z}}\rightarrow\left(\mathbb{N}_{0}\cup\{\infty\}\right)^{\mathbb{Z}}$
defined by
\[
\xi(x)_{n}=\lim_{r}\xi_{r}(x)_{n}
\]
 is well defined a.e.

\subsection{\protect\label{subsec:Proof-of-lemma}Proof of lemma \ref{lem:Marker=000020process} }

Fixing $n\in\mathbb{Z},x\in A^{\mathbb{Z}}$, we omit the evaluated
point and index for brevity and denote by $S_{r}=S_{r}(n,x)$ the
$r$-skeleton that covers $x$ at $n$ according to $\xi$ (which
is a.s. well defined). We remark that for any $k\geq0$, $\xi(x)_{n}\neq\xi_{k}(x)_{n}$
iff $\xi_{k}(x)_{n}=k$ and $\xi(x)_{n}>k$. Thus, $\xi$ and $\xi_{k}$
share the same $r$-skeletons for all $r\leq k$, i.e. taking $S_{r}(n,x)$
according to $\xi_{k}$ or $\xi$ gives the same skeleton, for any
$n$ and $x$.

We first point out a deterministic bound on the size of any skeleton
with $\min\mathrm{rank}=r>0$, which is completely based on the construction
of $\xi$ and has nothing to do with the properties of $X$ as a process.
This will imply the upper bound of lemma \ref{lem:Marker=000020process}
on the likely $r$-skeleton (corollary \ref{cor:.Upper=000020bound=000020on=000020|S_r|}
below), and will help in establishing the lower bound too (Proposition
\ref{prop:Lower=000020bound=000020on=000020|S_r|}).
\begin{prop}
\label{prop:|S_r|=000020Deterministic=000020upper=000020bound}For
any $r\geq0$ and any $x$ for which $S_{r}=S_{r}(0,x)$ is well defined,
at least one of the following holds:

Either
\[
|S_{r}|\leq2\Lambda_{d}L_{r}
\]
or
\[
S_{r}=S_{1}
\]
\end{prop}

\begin{proof}
This is shown inductively on $r$: for $r=1$ the second condition
trivially holds. Let $r>0$ and assume the conclusion takes place
for all $(r-1)$-skeletons of $\xi$.

Suppose first that $S_{r}=S_{r-1}$. By the induction assumption,
either we have 
\[
S_{r}=S_{r-1}=S_{1}
\]
or
\[
|S_{r}|=|S_{r-1}|\leq2\Lambda_{d}L_{r-1}\leq2\Lambda_{d}L_{r}
\]
where the last inequality follows from $\left(L_{n}\right)_{n\geq0}$
being increasing.

Suppose now that $S_{r}\neq S_{r-1}$. By the remark preceding the
present proposition, we can view $S_{r}$ and $S_{r-1}$ as the $r$
and $r-1$ skeletons that cover $x$ at $n$ according to $\xi_{r}$
(instead of taking it according to $\xi$). From that point of view,
the inequality $S_{r}\neq S_{r-1}$ indicates that the condition (\ref{eq:r-skeletonCond1})
takes place, i.e.
\[
\tau_{E_{r}}^{\geq\beta}-\tau_{E_{r}}^{<\beta}\leq\Lambda_{d}L_{r}.
\]
Using the same notation as in the construction of $\xi$ (subsection
\ref{subsec:Construction-of-alpha}) with $n=0$, and letting 
\[
\gamma=\max\left\{ i\leq\tau_{E_{r}}^{\geq\beta}(x):\xi_{r-1}(x)_{i}=r-1\right\} \geq\beta>0
\]
the interval at which $S_{r}$ covers $x$ is $[\alpha,\gamma-1]$,
thus
\begin{equation}
\left|S_{r}\right|=\gamma-\alpha\leq(\tau_{E_{r}}^{\geq\beta}-\tau_{E_{r}}^{<\beta})+(\tau_{E_{r}}^{<\beta}-\alpha)\label{eq:deterministic=000020bound=0000201}
\end{equation}
(the evaluated point $x$ has been omitted for brevity). Notice that
$\alpha$ is the smallest coordinate of $S_{r-1}(\tau_{E_{r}}^{<\beta},x)$.
Since $E_{r}\subset D_{r}\subset D_{1}$, one has $\mathbf{1}_{D_{1}}(x)_{\tau_{E_{r}}^{<\beta}}=1$,
so that $\tau_{E_{r}}^{<\beta}$ is the first coordinate of $S_{1}(\tau_{E_{r}}^{<\beta},x)$.
Thus, if $S_{r-1}(\tau_{E_{r}}^{<\beta},x)=S_{1}(\tau_{E_{r}}^{<\beta},x)$
is a $1$-skeleton, then $\alpha=\tau_{E_{r}}^{<\beta}$, and the
second summand on the right hand side of (\ref{eq:deterministic=000020bound=0000201})
is 0. Otherwise, by the induction assumption,
\[
\tau_{E_{r}}^{<\beta}-\alpha\leq|S_{r-1}(\tau_{E_{r}}^{<\beta},x)|\leq2\Lambda_{d}L_{r-1}\leq L_{r}
\]
(The last inequality uses proposition \ref{prop:(Properties-of-L_r)}).
Either way, one gets in (\ref{eq:deterministic=000020bound=0000201})
that
\[
|S_{r}|\leq L_{r}+\Lambda_{d}L_{r}\leq2\Lambda_{d}L_{r}
\]
as claimed.
\end{proof}
\begin{cor}
\label{cor:.Upper=000020bound=000020on=000020|S_r|}For some $c>0$
and any $r$,
\[
\mathbf{P}\left[\left|S_{r}\right|\leq2\Lambda_{d}L_{r}\right]=1-O(e^{-cL_{r}}).
\]
\end{cor}

\begin{proof}
By proposition \ref{prop:|S_r|=000020Deterministic=000020upper=000020bound},
\[
\mathbf{P}\left[\left|S_{r}\right|>2\Lambda_{d}L_{r}\right]=\mathbf{P}\left[S_{r}=S_{1}\,\wedge\,\left|S_{r}\right|>2\Lambda_{d}L_{r}\right]=\mathbf{P}\left[\left|S_{1}\right|>2\Lambda_{d}L_{r}\right]
\]
But since $\xi_{1}=\mathbf{1}_{D_{r}}$, we can use inequality (\ref{eq:-17})
to get
\[
\mathbf{P}\left[\left|S_{1}\right|>2\Lambda_{d}L_{r}\right]\leq\left(1-\eta L_{1}^{-1}\right)^{\lfloor\frac{2\Lambda_{d}L_{r}}{|w|N_{1}+n_{0}}\rfloor}
\]
and as $\eta,L_{1},N_{1},n_{0},|w|$ are all constants, the claim
follows.
\end{proof}
The next proposition gives a lower bound for the size of the likely
$r$-skeleton $S_{r}(X)$. This is done by showing that the interval
$[\tau_{E_{r}}^{<\beta}(x),\tau_{E_{r}}^{\geq\beta}(x)]$ is likely
to be a good approximation for the interval covered by $S_{r}(X)$,
which by the definition of $E_{r}$ must be of size at least $\Lambda_{-1}(L_{r})$.
\begin{prop}
\label{prop:Lower=000020bound=000020on=000020|S_r|}$\mathbf{P}\left[|S_{r}|<\frac{1}{2}\Lambda_{-1}(L_{r})\right]=\Delta(L_{r})$.
\end{prop}

\begin{proof}
On the event that (\ref{eq:r-skeletonCond1}) holds, that is, 
\[
\tau_{E_{r}}^{\geq\beta}(x)-\tau_{E_{r}}^{<\beta}(x)\leq\Lambda_{d}L_{r}
\]
 one has, similarly to (\ref{eq:deterministic=000020bound=0000201}),
that
\begin{equation}
\left|S_{r}\right|=\gamma-\alpha\geq(\tau_{E_{r}}^{\geq\beta}-\tau_{E_{r}}^{<\beta})-(\tau_{E_{r}}^{\geq\beta}-\gamma)\label{eq:Lower=000020bound=000020on=000020S_r=000020(1)}
\end{equation}
where $\gamma$ is the first coordinate of $S_{r-1}(\tau_{E_{r}}^{\geq\beta},X)$.
As in the proof of proposition \ref{prop:|S_r|=000020Deterministic=000020upper=000020bound},
the inclusion $E_{r}\subset D_{r}\subset D_{1}$ shows that the first
coordinate of $S_{1}(\tau_{E_{r}}^{\geq\beta},X)$ is $\tau_{E_{r}}^{\geq\beta}$.
Thus if $S_{r-1}(\tau_{E_{r}}^{\geq\beta},X)=S_{1}(\tau_{E_{r}}^{\geq\beta},X)$,
then
\begin{equation}
\tau_{E_{r}}^{\geq\beta}-\gamma=0.\label{eq:eq:Lower=000020bound=000020on=000020S_r=000020(2)}
\end{equation}
On the other hand, if $S_{r-1}(\tau_{E_{r}}^{\geq\beta},X)\neq S_{1}(\tau_{E_{r}}^{\geq\beta},X)$
then using Proposition \ref{prop:|S_r|=000020Deterministic=000020upper=000020bound}
one has
\begin{equation}
\tau_{E_{r}}^{\geq\beta}-\gamma\leq|S_{r-1}(\tau_{E_{r}}^{\geq\beta},X)|\leq2\Lambda_{d}L_{r-1}.\label{eq:Lower=000020bound=000020on=000020S_r=000020(3)}
\end{equation}
Either way, since by the definition of $E_{r}$ the size $\tau_{E_{r}}^{\geq\beta}-\tau_{E_{r}}^{<\beta}$
must be at least $\Lambda_{-1}L_{r}$, we get from (\ref{eq:Lower=000020bound=000020on=000020S_r=000020(1)}),
(\ref{eq:eq:Lower=000020bound=000020on=000020S_r=000020(2)}) and
(\ref{eq:Lower=000020bound=000020on=000020S_r=000020(3)}) that:
\[
\left|S_{r}\right|\geq\Lambda_{-1}L_{r}-2\Lambda_{d}L_{r-1}\geq\frac{1}{2}\Lambda_{-1}L_{r}.
\]
(for the last inequality we used Proposition \ref{prop:(Properties-of-L_r)}).
Thus
\[
\mathbf{P}\left[|S_{r}|\geq\frac{1}{2}\Lambda_{-1}L_{r}\right]\geq\mathbf{P}\left[\tau_{E_{r}}^{\geq\beta}-\tau_{E_{r}}^{<\beta}\leq\Lambda_{d}L_{r}\right]
\]
and it remains to give a lower bound for the quantity on the right
hand side. On the complement event $B=\left\{ \tau_{E_{r}}^{\geq\beta}-\tau_{E_{r}}^{<\beta}>\Lambda_{d}L_{r}\right\} $,
the coordinate $\beta$ is the right end coordinate of $S_{r}(0)$.
So either $|S_{r}(n)|>2\Lambda_{d}L_{r}$ (call that event $A$) or
$\beta$ must be in the range $n<\beta\leq n+2\Lambda_{d}L_{r}$,
and in particular, for some $n<i\leq n+2\Lambda_{d}L_{r}$ one has
that
\[
\tau_{E_{r}}^{\geq i}-\tau_{E_{r}}^{<i}>\Lambda_{d}L_{r}
\]
(call each of these events $(B_{i})_{1\leq i\leq2\Lambda_{d}L_{r}}$).
As $B\subset A\cup\bigcup_{i=1}^{2\Lambda_{d}L_{r}}B_{i}$, we have
\begin{align*}
\mathbf{P}[B] & \leq\mathbf{P}[A]+\sum_{i=1}^{2\Lambda_{d}L_{r}}\mathbf{P}[B_{i}]\\
 & =\mathbf{P}[|S_{r}|>2\Lambda_{d}L_{r}]+\sum_{i=1}^{2\Lambda_{d}L_{r}}\mathbf{P}[\tau_{E_{r}}^{\geq i}-\tau_{E_{r}}^{<i}>\Lambda_{d}L_{r}]\\
 & =(2\Lambda_{d}L_{r}+2)\Delta(L_{r})\\
 & =\Delta(L_{r})
\end{align*}
where the one before the last inequality follows from Corollary \ref{cor:.Upper=000020bound=000020on=000020|S_r|}
and Equation (\ref{eq:-4}).
\end{proof}
We now turn to show that $S_{r}(X)$ is likely to be determined by
$X_{[-3\Lambda_{d}L_{r},3\Lambda_{d}L_{r}]}$ (Corollary \ref{cor:S_r=000020determined=000020by}).
Recall that for a factor map $\phi(X)=Y$ between two processes, its
coding radius $R_{\phi}$ is defined by
\[
R_{\phi}(x)=\min\left\{ r\geq0:\,(X_{i}=x_{i})_{|i|\leq r}\quad\Longrightarrow\quad\phi(X)_{0}=\phi(x)_{0}\;a.s.\right\} 
\]
Viewing the marker processes $\xi_{r}$ as factors of $X$, we make
the following observation on their coding radii :
\begin{prop}
\label{prop:radius=000020coding=000020of=000020alpha_r}Assuming $N_{1}$
is large enough, for any $r\geq1$, 
\[
R_{\xi_{r}}\leq2\Lambda_{d}L_{r}\qquad a.s.
\]
\end{prop}

\begin{proof}
For $r=1$, $\xi_{1}(\omega)_{0}$ is determined by the coordinates
$[-|w|,|w|(N_{1}-1)]$. Thus, assuming $N_{1}$ is large enough, we
have
\begin{equation}
R_{\xi_{1}}\leq|w|(N_{1}-1)<2\Lambda_{d}L_{1}.\label{eq:-6}
\end{equation}

Let $r>1$, and suppose the claim holds for all $k<r$. From the definition
of $\xi_{r}$ together with Proposition \ref{prop:|S_r|=000020Deterministic=000020upper=000020bound},
one can observe that $\xi_{r}(X)_{0}$ is completely determined by
the samples $X_{[-M,M]}$ and $\xi_{r-1}(X)_{[-M,M]}$, where
\[
M=\underbrace{\Lambda_{d}L_{r}}_{{\scriptscriptstyle M_{1}}}+\underbrace{2\Lambda_{d}L_{r-1}}_{{\scriptscriptstyle M_{2}}}+\underbrace{\Lambda_{-1}L_{r}+|w|(N_{r}-1)}_{{\scriptscriptstyle M_{3}}}.
\]
We briefly explain this observation: To determine $\xi_{r}(X)_{0}$,
first check the value of $\xi_{r-1}(X)_{0}$: if it is lesser than
$r-1$ then put $\xi_{r}(X)_{0}:=\xi_{r-1}(X)_{0}$ (this was remarked
already at the beginning of the present subsection). Otherwise, if
$\xi_{r-1}(X)_{0}=r-1$, denoting by $\beta$ the first positive coordinate
with $\xi_{r-1}(X)_{\beta}=r-1$, we need to check whether 
\begin{equation}
\tau_{E_{r}}^{\geq\beta}-\tau_{E_{r}}^{<\beta}\leq\Lambda_{d}L_{r}=M_{1}\label{eq:-8}
\end{equation}
or not, and accordingly put $\xi_{r}(X)_{0}:=r-1$ (if (\ref{eq:-8})
takes place but $0\neq\max\{i\leq\tau_{E_{r}}^{<\beta}:\xi_{r-1}(x)_{i}=r-1\}$)
or $r$ (otherwise). For that, notice first that $\tau_{E_{r}}^{<\beta}<\beta$.
As $\beta$ is the last coordinate of $S_{r-1}$, by Proposition \ref{prop:|S_r|=000020Deterministic=000020upper=000020bound}
we have that either $\beta\leq2\Lambda_{d}L_{r-1}$, which implies
$\tau_{E_{r}}^{<\beta}\leq2\Lambda_{d}L_{r-1}$, or $S_{r-1}=S_{0}$,
which in that case $\tau_{E_{r}}^{<\beta}\leq\tau_{D_{1}}^{<\beta}\leq0$.
In both cases we get that
\[
\tau_{E_{r}}^{<\beta}\leq2\Lambda_{d}L_{r-1}=M_{2}
\]
and so it is enough to know $\mathbf{1}\{\sigma^{n}X\in E_{r}\}$
and $\xi_{r-1}(X)_{n}$ within $n\in[-M_{1},M_{1}+M_{2}]\subset B_{0}(M_{1}+M_{2})$
in order to verify whether (\ref{eq:-8}) holds or not. By the definition
of $E_{r}$ (\ref{eq:-1}), for $M_{3}=\Lambda_{-1}L_{r}+|w|(N_{r}-1)$
we have for any $n\in\mathbb{Z}$ that
\[
X_{B_{n}(M_{3})}\Rrightarrow\mathbf{1}\{\sigma^{n}X\in E_{r}\}
\]
thus taking $M=M_{1}+M_{2}+M_{3}$ we have
\begin{align*}
\left(X_{B_{0}(M)},\xi_{r-1}(X)_{B_{0}(M)}\right) & \Rrightarrow\left(\mathbf{1}\{\sigma^{n}X\in E_{r}\}_{B_{0}(M_{1}+M_{2})},\xi_{r-1}(X)_{B_{0}(M_{1}+M_{2})}\right)\\
 & \Rrightarrow\xi_{r}(X)_{0}
\end{align*}
which is just as stated. 

By the induction assumption, for any $n\in[-M,M]$, the value $\xi_{r-1}(X)_{n}$
is determined by
\[
X_{[u-2\Lambda_{d}L_{r-1},u+2\Lambda_{d}L_{r-1}]}
\]
thus the whole sample $\xi_{r-1}(X)_{[-M,M]}$ is determined by
\[
X_{[-M-2\Lambda_{d}L_{r-1},M+2\Lambda_{d}L_{r-1}]}
\]
which of course includes the sample $X_{[-M,M]}$, so $\xi_{r}(X)$
is determined by it too. Thus by Proposition \ref{prop:(Properties-of-L_r)}
(while assuming $N_{1}$ is large enough),
\begin{align*}
R_{\xi_{r}} & \leq M+2\Lambda_{d}L_{r-1}\leq2\Lambda_{d}L_{r}
\end{align*}
and the proof is complete.
\end{proof}
A direct corollary of it is the following:
\begin{cor}
\label{cor:S_r=000020determined=000020by}For any $r$, any $r$-skeleton
$S$ that covers $X$ at $[a,b]$ is determined by $X_{[a-2\Lambda_{d}L_{r},b+2\Lambda_{d}L_{r}]}$. 
\end{cor}

Putting together all the results of this section, the marker process
lemma follows easily:

\textit{Proof of Lemma \ref{lem:Marker=000020process} (Marker process
Lemma).} We show how the marker process $\xi$ satisfies all required
properties of the lemma. Fix $n$ and let 
\[
r=\max\{r':6\Lambda_{d}L_{r'}\leq n\}.
\]

By Proposition \ref{prop:(Properties-of-L_r)} one has
\begin{equation}
n\leq6\Lambda_{d}L_{r+1}\leq\Lambda_{3d}L_{r}.\label{eq:n=000020is=000020bonded=000020by=000020L_r}
\end{equation}
By Corollaries \ref{cor:.Upper=000020bound=000020on=000020|S_r|}
and \ref{cor:S_r=000020determined=000020by}, 
\begin{align*}
\mathbf{P}\left[X_{[-n,n]}\Rrightarrow S_{r}\,,\,|S_{r}|\leq n\right] & =1-\Delta(L_{r})\\
 & =1-\Delta(n)
\end{align*}
(the second equality is by (\ref{eq:n=000020is=000020bonded=000020by=000020L_r}).)
By Proposition \ref{prop:Lower=000020bound=000020on=000020|S_r|}
together with (\ref{eq:n=000020is=000020bonded=000020by=000020L_r}),
\[
\mathbf{P}[|S_{r}(X)|<\Lambda_{-3d-1}(n)\leq\Lambda_{-1}L_{r}]=\Delta(n)
\]
and the proof of the Marker process Lemma is complete.$\hfill\square$

\section{Two processes with a common state\protect\label{sec:Common=000020entry}}

In this section we develop some of the main ingredients for constructing
the isomorphism of our main theorem, theorem \ref{thm:Main=000020Thm}.
At the end of this section, we will have enough tools to establish
Theorem \ref{thm:common=000020entry=000020version}, which is a weak
version of the main theorem. To state the theorem, we use the notion
of marked processes: 
\begin{defn}
\label{def:A-marked-process}\textit{(marked process)} A marked process
$X\sim\mathbf{A}=(A^{\mathbb{Z}},\mu,a)$ is an aperiodic Markov process
$(A^{\mathbb{Z}},\mu)$ endowed with a renewal state $a\in A$ that
satisfies the following ``probability reduction'' property:
\end{defn}

\begin{itemize}
\item \textit{(Probability reduction property)} There exists $c<1$ s.t.
for any $n\geq1$ and any $w\in(A\backslash\{a\})^{n}$,
\begin{equation}
\mathbf{P}[X_{[0,n+1]}=awa|\widehat{X}_{[0,n+1]}^{a}=10^{n}1]\in[0,c]\cup\{1\}.\label{eq:Marked=000020state=000020property}
\end{equation}
\end{itemize}
For i.i.d. processes, this property holds for all states. The same
is true for any $k$-stringing of an i.i.d. process (i.e. $X=(Y_{n}^{(k)})_{n}:=(Y_{n}\cdots Y_{n+k-1})_{n}$
for some i.i.d. process $Y$ and some fixed $k$). The property becomes
less trivial when dealing with Markov processes, and is not valid
in general. However, as it is shown in lemma \ref{lem:marked=000020process=000020lemma},
the $k$-stringing $X^{(k)}$ of any aperiodic Markov process $X$
admits a renewal state satisfying the above property, provided $k$
is large enough. 

For a marked process $X\sim\mathbf{A}$, the process $\widehat{X}$
will always refer to the distribution of the distinct state $a$,
and a marker process of $X$ will always refer to a marker factor
of $\widehat{X}$.
\begin{thm}
\label{thm:common=000020entry=000020version}Let $X,Y$ be two marked
processes of equal entropy $h(X)=h(Y)$ and equal marked-state distribution
$\widehat{X}\overset{{\scriptscriptstyle dist.}}{=}\hat{Y}$.

There exists a finitary isomorphism $\phi(X)\overset{{\scriptscriptstyle dist.}}{=}Y,\phi^{-1}(Y)\overset{{\scriptscriptstyle dist.}}{=}X$
with coding radii $R_{\phi}$ and $R_{\phi^{-1}}$ satisfying:
\[
\mathbf{P}[R_{\phi}>n],\mathbf{P}[R_{\phi^{-1}}>n]=\tilde{O}(n^{-\frac{1}{2}}).
\]
\end{thm}

The main theorem, Theorem \ref{thm:Main=000020Thm}, will be established
in section \ref{sec:The-general-case} by showing that between any
two aperiodic Markov processes one can compose a chain of isomorphisms
as above, while controlling the coding radius tails of the resulted
isomorphism.

\subsection{\protect\label{subsec:Definitions-and-preliminear}Fillers and systems
of fillers (SOF)}

Given a skeleton $S=(\widehat{x}_{[m,n]},\xi(x)_{[m,n]})$ for a marked
process $X\sim(A^{\mathbb{Z}},\mu,a)$ and a marker factor $\xi$,
a \textit{filler} of $S$ is any word $w\in A^{[m,n]}$ with $\hat{w}=\widehat{x}_{[m,n]}$.
The \textit{filler space} of $S$, denoted by $\left(A_{S},\mu_{S}\right)$,
is the set of fillers 
\[
A_{S}=\left\{ w\in A^{[m,n]}:\widehat{w}=\widehat{x}_{[m,n]}\right\} 
\]
endowed with the probability measure
\begin{align*}
\mu_{S}(w) & =\mathbf{P}[X_{[m,n]}=w|X_{[m,n+1]}=\widehat{x}_{[m,n+1]}]\\
 & =\mathbf{P}[X_{[m,n]}=w|\widehat{X}=\widehat{x}]
\end{align*}
where the last inequality follows from the fact that a skeleton obtained
from a marker processes has in its first coordinate the specified
renewal state of the process (and the same applies to the consecutive
skeleton of $x$). Consequently, if a skeleton $S$ has a decomposition
into sub-skeletons
\[
S=S^{1}S^{2}\cdots S^{m}
\]
then $A_{S}=A_{S^{1}}\times\cdots\times A_{S^{m}}$ and
\[
\mu_{S}=\mu_{S^{1}}\times\cdots\times\mu_{S^{m}}.
\]

Consider the disintegration of $\mu$ by the law $\widehat{\mu}$
of $\widehat{X}$,
\[
\mu=\int\mu_{\omega}d\widehat{\mu}(\omega)
\]
where for each $\omega\in\{0,1\}^{\mathbb{Z}}$, the measure $\mu_{\omega}$
is the conditional measure of $\mu$ on the fiber $\{x\in A^{\mathbb{Z}}:\widehat{x}=\omega\}$.
By the marker process lemma (Lemma \ref{lem:Marker=000020process}),
there is a nice marker factor $\xi$ on $(\{0,1\}^{\mathbb{Z}},\widehat{\mu})$
(see Def. \ref{def:nice=000020marker} for the definition of nice
marker factor). Given $r\in\mathbb{N}$ and $\omega\in\{0,1\}^{\mathbb{Z}}$,
let 
\[
...S^{-1}(\omega)S^{0}(\omega)S^{1}(\omega)...
\]
 be the decomposition of $\omega$ into $r$-skeletons according to
$\xi$ (say, $S^{0}(\omega)$ denotes the $r$-skeleton that covers
$\omega$ at $0$). We have that a.s. all the $r$-skeletons in the
decomposition are finite and begin with the renewal state $a$. It
can be easily observed that for $\widehat{\mu}$-a.e. $\omega$ and
any $r$, the conditional measure $\mu_{\omega}$ can be written as
the product measure
\[
\mu_{\omega}=\times_{i\in\mathbb{Z}}\mu_{S^{i}(\omega)}.
\]

A filler space $\left(A_{S},\mu_{S}\right)$ can be endowed with an
equivalence relation $\sim_{S}$ on the fillers. In that case, an
equivalence class will be denoted by $[w]_{\sim_{S}}$ or just $[w]$,
and the collection of classes is denoted by
\[
\tilde{A}_{S}=\left\{ [w]:w\in A_{S}\right\} 
\]

\begin{defn}
\label{def:SOF}\textit{(SOF) Given a marked process }$\mathbf{A}$
and its marker factor\textit{ $\xi$,} A \textit{system of fillers
}(SOF) $\mathscr{A}=\left(A_{S},\sim_{S}\right)_{S}$ for $\mathbf{A}$
and $\xi$ is a collection of equivalence relations $\sim_{S}$ on
all filler spaces $A_{S}$ of $\mathbf{A}$ that satisfies the following
properties: 
\end{defn}

\begin{enumerate}
\item \label{enu:The-collection-is}The collection is translation-invariant,
in the sense that for any skeleton $S=S(o,x)$ and its translation
$\sigma S:=S(o-1,\sigma x)$, the relation $\sim_{\sigma S}$ is the
translation $\sim_{\sigma S}=\sigma\otimes\sigma(\sim_{S})$ of the
relation $\sim_{S}$.
\item \label{enu:For-any-skeleton}For any skeleton decomposition $S=S^{1}\cdots S^{m}$,
the relation $\sim_{S}$ is a refinement of the product relation $\times_{i=1}^{m}\sim_{S^{i}}$,
defined by
\[
[w]_{\times_{i=1}^{m}\sim_{S^{i}}}=[w_{1}]_{\sim_{S^{1}}}\times\cdots\times[w_{m}]_{\sim_{S^{m}}}
\]
where $w_{i}$ is the sub-word of $w$ that fills $S^{i}$.
\item \label{enu:For-,-the}For $X\sim\mathbf{A}$ and $S_{r}=S_{r}(X)$,
its random $0$-coordinate skeleton of minimal rank $r$, the expected
fraction of coordinates $i\in S_{r}$ that are determined by the class
of $X_{S_{r}}$ tends to $1$ as $r\rightarrow\infty$ (see Def. \ref{def:Let--be}):
\[
\mathbf{E}\left[\frac{1}{|S_{r}|}\sum_{i}\mathbf{1}\{[X_{S_{r}}]\Rrightarrow X_{i}\}\right]\overset{{\scriptscriptstyle r\rightarrow\infty}}{\rightarrow}1.
\]
\end{enumerate}
Suppose we are given an SOF for $\mathbf{A}$. Identifying each class
$E\in\tilde{A}_{S}$ with the set $\left\{ x\in A^{\mathbb{Z}}:x_{S}\in E\right\} $,
define for each $r\geq0$ the algebra $\mathcal{A}_{r}\subset Borel\left(A^{\mathbb{Z}}\right)$
generated by all classes for all $r$-skeletons that cover zero:
\[
\mathcal{A}_{r}=\text{algebra}\left\{ E\in\tilde{A}_{S_{r}(x)}:x\in A^{\mathbb{Z}}\right\} 
\]

By property \ref{enu:For-any-skeleton} of the definition of SOF it
follows that $\left(\mathcal{A}_{r}\right)_{r\geq1}$ is increasing
(w.r.t. inclusion).
\begin{prop}
\label{prop:System=000020of=000020fillers=000020generates=000020the=000020algebra}
Suppose $\mathbf{A}$ has an SOF $\mathscr{A}$, and let $X\sim\mathbf{A}$.
Then a.s. for any $k$, the event $[X_{S_{r}(X)}]\Rrightarrow X_{[-k,k]}$
holds for large enough $r$. 

Consequently, the sequence of algebras $\left(\mathcal{A}_{r}\right)_{r\geq1}$
is an increasing sequence converges to the Borel $\sigma$-algebra
of $A^{\mathbb{Z}}$ (modulo $\mu$-null sets).
\end{prop}

\begin{proof}
Using a mass-transport argument, one has
\begin{align*}
\mathbf{P}\left[[X_{S_{r}}]\not\Rrightarrow X_{[-k,k]}\right] & =\mathbf{E}\left[\frac{1}{|S_{r}|}\sum_{i\in S_{r}}\mathbf{1}\{[X_{S_{r}}]\not\Rrightarrow X_{[i-k,i+k]}\}\right]\\
 & \leq(2k+1)\mathbf{E}\left[\frac{1}{|S_{r}|}\sum_{i\in S_{r}}\mathbf{1}\{[X_{S_{r}}]\not\Rrightarrow X_{i}\}\right]+2k\mathbf{E}\left[\frac{1}{|S_{r}|}\right]
\end{align*}
By property \ref{enu:For-,-the} in the definition of SOF, the first
summand of the last expression tends to 0 as $r\rightarrow\infty$.
By the second property of a marker process, $|S_{r}|\rightarrow\infty$
as $r\rightarrow\infty$ a.s. , thus the second summand of the last
expression tends to 0 as $r\rightarrow\infty$ (notice that $|S_{r}|\geq1$
a.s.). Thus a.s. there exists some $r'$ for which $[X_{S_{r'}}]\Rrightarrow X_{[-k,k]}$.
On the other hand, property \ref{enu:For-any-skeleton} in the definition
of SOF implies that for any $r\geq r'$,
\[
[X_{S_{r}}]\subset[X_{S_{r'}}]
\]
(viewed as sets in $A^{\mathbb{Z}}$) hence $[X_{S_{r}}]\Rrightarrow[X_{S_{r'}}]\Rrightarrow X_{[-k,k]}$.
This proves the first claim of the proposition.

The above claim implies that for any $k$, the algebra $\sigma\left(X_{[-k,k]}\right)$
can be approximated arbitrarily well by $\mathcal{A}_{r}$ with large
enough $r$. Thus 
\[
\sigma\left((X_{n})_{n\in\mathbb{Z}}\right)\subset\sigma\left(\bigvee_{k\in\mathbb{Z}}\left(X_{[-k,k]}\right)\right)\subset\sigma\left(\bigvee_{r\in\mathbb{Z}}\mathcal{A}_{r}\right)
\]
 and the conclusion follows.
\end{proof}
Let $X\in\mathcal{X},Y\in\mathcal{Y}$ be two random variables defined
on the same probability space. The information function is defined
for any $x\in\mathcal{X},y\in\mathcal{Y}$ by
\begin{align*}
I(x) & =-\log\mathbf{P}[X=x]\\
I(x,y) & =-\log\mathbf{P}[(X,Y)=(x,y)]\\
I(x|y) & =I(x,y)-I(y)
\end{align*}
Taking an independent copy of the random variables $(X',Y')\overset{D}{=}(X,Y)$,
we define the random information
\begin{align*}
I(X) & =-\log\mathbf{E}[\mathbf{1}\{X'=X\}|X]\\
I(X,Y) & =-\log\mathbf{E}[\mathbf{1}\{(X',Y')=(X,Y)\}|(X,Y)]\\
I(X|Y) & =I(X,Y)-I(Y)
\end{align*}
Occasionally we will consider the information function with respect
to a probability measure $\mu$ that differs from the universal one
$\mathbf{P}$. In that case we will add the subscript $I_{\mu}(x)=-\log\mu(\{x\})$.
\begin{defn}
\label{def:Good=000020SOF} \textit{(Good SOF)} Suppose that the marker
factor $\xi$ of $\mathbf{A}$ is nice (see Def. \ref{def:nice=000020marker})
and admits a lengths sequence $(L_{r})_{r}$. An SOF for $\mathbf{A}$
and $\xi$ is said to be \textit{good} if there are some $K\geq k$
such that:
\end{defn}

\begin{enumerate}
\item \label{enu:For-any-large}\textit{(Information Upper-Bound) }For any
skeleton $S$ of large enough size $|S|$, any class $[w]\in\tilde{A}_{S}$
satisfies
\[
I_{\mu_{S}}([w])\leq|S|h(X|\widehat{X})-\Lambda_{k,\frac{1}{2}}|S|.
\]
\item \label{enu:Good=000020SOF=0000202}\textit{(Information Lower-Bound)
}With probability $1-\Delta(L_{r})$, a lower bound for $I([X_{S_{r}}]|\widehat{X}_{S_{r}})=I_{\mu_{S_{r}}}([X_{S_{r}}])$
also holds:
\[
\mathbf{P}\left[I([X_{S_{r}}]|\widehat{X})\geq|S_{r}|h(X|\widehat{X})-\Lambda_{K,\frac{1}{2}}|S_{r}|\right]=1-\Delta(L_{r}).
\]
\end{enumerate}
Given a skeleton $S$ and $n\in\mathbb{N}$, write $N_{\mathbf{A}}(S,n)$
for the union's conditional measure of the $n$ largest classes in
$\tilde{A}_{S}$:
\[
N_{\mathbf{A}}(S,n)=\max\{\mu_{S}(\bigcup_{i=1}^{n}U_{i}):\,U_{i}\in\tilde{A}_{S}\}.
\]
What so good in a good SOF is that it satisfies the formula of the
following proposition. We use the notation 
\[
\bar{A}_{S}:=\tilde{A}_{S^{1}}\times\cdots\times\tilde{A}_{S^{m}}
\]
where $S=S^{1}\cdots S^{m}$ is the decomposition of an $r$-skeleton
$S$ into its $(r-1)$-sub-skeletons.
\begin{prop}
\label{prop:key=000020property=000020s.o.f.}Let $X\sim\mathbf{A}$
and $Y\sim\mathbf{B}$ be two marked processes of equal entropy, equal
marked-state distribution, and a common nice marker factor $\xi$
exhibiting lengths sequence $(L_{r})_{r\geq1}$ (see Def. \ref{def:nice=000020marker}).
Suppose both $X$ and $Y$ have good SOFs with respect to $\xi$.
By thinning out $\xi$ if needed (while keeping it nice), we have:
\[
\mathbf{E}[N_{\mathbf{A}}(S_{r},|\bar{B}_{S_{r}}|)]=\Delta(L_{r}).
\]
\end{prop}

The need for this specific formula will be clearer in the isomorphism
construction of Subsection \ref{subsec:The-isomorphism-construction}
below. There, $\mathbf{E}[N_{\mathbf{A}}(S_{r},|\bar{B}_{S_{r}}|)]$
serves as a bound for the measure of the set of classes in $\tilde{A}_{S_{r}}$
that are not encoded into any $\bar{B}_{S_{r}}$-class (See also the
marriage lemma in Proposition \ref{prop:Properties=000020of=000020societies}
below).

In order to prove Proposition \ref{prop:key=000020property=000020s.o.f.},
we need the following fact about concave functions that vanish at
zero:
\begin{prop}
Let $f:[0,\infty)\rightarrow[0,\infty)$ be a concave function that
vanishes at $0$. For any non-negative numbers $x_{1},...,x_{n}$
one has
\[
f(\sum_{i=1}^{n}x_{i})\leq\sum_{i=1}^{n}f(x_{i}).
\]
More generally, for any $\lambda\geq\max_{i}x_{i}$,
\[
\lambda^{-1}(\sum_{i=1}^{n}x_{i})f(\lambda)\leq\sum_{i=1}^{n}f(x_{i}).
\]
\end{prop}

\begin{proof}
Applying Jensen's inequality, one has
\[
\lambda^{-1}(\sum_{i=1}^{n}x_{i})f(\lambda)=\sum_{i=1}^{n}\left(\frac{x_{i}}{\lambda}f(\lambda)+\frac{\lambda-x_{i}}{\lambda}f(0)\right)\leq\sum_{i=1}^{n}f(x_{i})
\]
\end{proof}
We shall use this proposition for the functions $\Lambda_{k,\frac{1}{2}}$
($k\geq0$) in the following way (recall that $\Lambda_{a,b}(x):=x^{b}\log^{a}x$):
fixing $k$, although $\Lambda_{k,\frac{1}{2}}$ is not necessarily
concave on all $[0,\infty)$, and even not defined on $0$, for some
large enough $t_{0}$ we have that
\begin{itemize}
\item $\Lambda_{k,\frac{1}{2}}''(t)<0$ for all $t\in[t_{0},\infty)$.
\item $\Lambda_{k,\frac{1}{2}}'(t_{0})\leq\Lambda_{k,\frac{1}{2}}(t_{0})/t_{0}$.
\end{itemize}
Fixing such a $t_{0}$, the function 
\[
f(t)=\begin{cases}
\frac{t}{t_{0}}\Lambda_{k,\frac{1}{2}}(t_{0}) & 0\leq t<t_{0}\\
\Lambda_{k,\frac{1}{2}}(t) & t\geq t_{0}
\end{cases}
\]
is concave on $[0,\infty)$ and vanishes on $0$. Thus we get:
\begin{cor}
\label{cor:For-any-,}For any $k\geq0$, if $t_{0}$ is large enough
then for any $x_{1},...,x_{n}\geq t_{0}$ and any $\lambda\geq\max_{i}x_{i}$,
\[
\lambda^{-1}(\sum_{i=1}^{n}x_{i})\Lambda_{k,\frac{1}{2}}(\lambda)\leq\sum_{i=1}^{n}\Lambda_{k,\frac{1}{2}}(x_{i}).
\]
\end{cor}

\textit{Proof of Proposition \ref{prop:key=000020property=000020s.o.f.}:}
Let $Y\sim\mathbf{B}.$ As the SOF for $\mathbf{B}$ is good, Property
\ref{enu:For-any-large} of good SOF (Def. \ref{def:Good=000020SOF})
ensures the existence of some $k$, s.t. for any $S$ of large enough
size $|S|$, the conditional measure of any class $V\in\tilde{B}_{S}$
is at least
\[
\nu_{S}(V)=2^{-I_{\nu_{S}}(V)}\geq2^{-|S|h(Y|\widehat{Y})+\Lambda_{k,\frac{1}{2}}|S|}.
\]
Thus, letting $S_{r}=S^{1}S^{2}\cdots S^{m}$ being the decomposition
of $S_{r}$ into its $(r-1)$-sub-skeletons, 
\[
|\bar{B}_{S_{r}}|=\prod_{i=1}^{m}|\tilde{B}_{S^{i}}|\leq\prod_{i=1}^{m}\left(\min_{V\in\tilde{B}_{S^{i}}}\nu_{S^{i}}(V)\right)^{-1}\leq2^{|S|h(Y|\widehat{Y})-\sum_{i=1}^{m}\Lambda_{k,\frac{1}{2}}|S^{i}|}.
\]
In addition, having a good SOF for $\mathbf{A}$, Property \ref{enu:Good=000020SOF=0000202}
of good SOF provides some $K\geq k$ s.t. with probability $1-\Delta(L_{r})$,
\begin{equation}
\mu_{S_{r}}([X_{S_{r}}])=2^{-I([X_{S_{r}}]|\widehat{X}_{S_{r}})}\leq2^{-|S_{r}|h(X|\widehat{X})+\Lambda_{K,\frac{1}{2}}|S_{r}|}.\label{eq:-7}
\end{equation}
For a skeleton $S$, write 
\[
E(S)=\left\{ U\in\tilde{A}_{S}:\mu_{S}(U)>2^{-|S|h(X|\widehat{X})+\Lambda_{K,\frac{1}{2}}|S|}\right\} .
\]
 We have
\begin{align}
N_{\mathbf{A}}(S_{r},|\bar{B}_{S_{r}}|) & \leq\mu_{S_{r}}(E(S_{r}))+2^{-|S|h(X|\widehat{X})+\Lambda_{K,\frac{1}{2}}|S_{r}|}\cdot|\bar{B}_{S_{r}}|\nonumber \\
 & \leq\mu_{S_{r}}(E(S_{r}))+2^{-\left(\sum_{i=1}^{m}\Lambda_{k,\frac{1}{2}}|S^{i}|-\Lambda_{K,\frac{1}{2}}|S_{r}|\right)}\label{eq:key=000020property=000020s.o.f.=0000203}
\end{align}

and consequently, the expectation of the left hand side is bounded
by the expectation of the right hand side. Notice that
\begin{align}
\mathbf{E}[\mu_{S_{r}}(E(S_{r}))] & =\mathbf{P}\left[I([X_{S_{r}}]|\widehat{X}_{S_{r}})<|S_{r}|h(X|\widehat{X})-\Lambda_{K,\frac{1}{2}}|S_{r}|\right]=\Delta(L_{r}),\label{eq:key=000020property=000020s.o.f.=0000204}
\end{align}
so it only remains to show that
\begin{equation}
\mathbf{E}\left[2^{-\left(\sum_{i=1}^{m}\Lambda_{k,\frac{1}{2}}|S^{i}|-\Lambda_{K,\frac{1}{2}}|S_{r}|\right)}\right]=\Delta(L_{r}).\label{eq:-66}
\end{equation}
As $\xi$ is nice, independently of $r$ there is some $c>0$ so that
with probability $1-\Delta(L_{r})$, all $(r-1)$-sub-skeletons $S^{i}$
of $S_{r}$ satisfy 
\begin{equation}
|S^{i}|\in[\Lambda_{-c}L_{r-1},L_{r-1}]\label{eq:key=000020property=000020s.o.f.=0000201}
\end{equation}
and 
\begin{equation}
|S_{r}|\in[\Lambda_{-c}L_{r},L_{r}].\label{eq:key=000020property=000020s.o.f.=0000202}
\end{equation}
Denote by $B_{r}$ the event that both (\ref{eq:key=000020property=000020s.o.f.=0000201})
and (\ref{eq:key=000020property=000020s.o.f.=0000202}) hold for $r$.
Let 
\[
D_{r}=\bigcap_{|i|\leq L_{r}}\sigma^{i}B_{r}
\]
which has also probability $1-\Delta(L_{r})$. For $r$ large enough,
letting $t_{0}$ be as in Corollary \ref{cor:For-any-,}, one has
that
\begin{equation}
\forall1\leq i\leq m\qquad|S^{i}|>t_{0}.\label{eq:-67}
\end{equation}
Thus, on $D_{r}$, for any $\lambda>|S_{r}|$, we can use Corollary
to bound the first summand in the exponent of (\ref{eq:-66}) by
\[
\sum_{i=1}^{m}\Lambda_{k,\frac{1}{2}}|S^{i}|\geq\lambda^{-1}(\sum_{i=1}^{n}|S^{i}|)\Lambda_{k,\frac{1}{2}}(\lambda)=\frac{|S_{r}|}{\lambda}\Lambda_{k,\frac{1}{2}}(\lambda).
\]
Taking $M>2(K-k+1)$ and $\lambda:=\Lambda_{-M}|S_{r}|$ in the last
inequality, one gets that
\[
\sum_{i=1}^{m}\Lambda_{k,\frac{1}{2}}|S^{i}|=\Lambda_{K+1,\frac{1}{2}}|S_{r}|-o(\log^{M+1}|S_{r}|)
\]
which gives a bound on the exponent in (\ref{eq:-66}):
\begin{align*}
\sum_{i=1}^{m}\Lambda_{k,\frac{1}{2}}|S^{i}|-\Lambda_{K,\frac{1}{2}}|S_{r}| & \geq\Lambda_{K+1,\frac{1}{2}}|S_{r}|-\Lambda_{K,\frac{1}{2}}|S_{r}|-o(\log^{M+1}|S_{r}|)\\
 & \geq|S_{r}|^{\frac{1}{2}}\\
 & \geq\Lambda_{-c}L_{r}
\end{align*}
This, together with the fact that $D_{r}$ has probability $1-\Delta(L_{r})$,
implies (\ref{eq:-67}). 

Overall, taking expectations of both sides of (\ref{eq:key=000020property=000020s.o.f.=0000203}),
one gets
\begin{align*}
\mathbf{E}[N_{\mathbf{A}}(S_{r},|\bar{B}_{S_{r}}|)] & \leq\Delta(L_{r})+\Delta(L_{r})
\end{align*}
and the proof is complete.$\hfill\square$

\subsection{\protect\label{subsec:Azuma-Hoeffding-inequality-for}Azuma-Hoeffding
inequality for Markov chains and renewal processes}

We shall use the following concentration result, commonly referred
to as Azuma's inequality (see \cite{key-14} and \cite[Chapter 4]{key-15})
\begin{thm}
\label{thm:Azuma-Hoeffding-inequality}(Azuma-Hoeffding inequality)
Consider a probability space $(\Omega,\mathscr{F},\mathbf{P})$ with
filtration
\[
\{\emptyset,\Omega\}=\mathscr{F}_{0}\leq\mathscr{F}_{1}\leq\cdots\leq\mathscr{F}_{n}=\mathscr{F}
\]
Let $C>0$, and suppose that a random variable $X$ has the property
that 
\begin{equation}
\Vert\mathbf{E}[X|\mathscr{F}_{i+1}]-\mathbf{E}[X|\mathscr{F}_{i}]\Vert_{\infty}\leq C\qquad\forall i=0,...,n-1.\label{eq:-5}
\end{equation}
Then for any $\epsilon>0$
\[
\mathbf{P}[|X-\mathbf{E}[X]|>\epsilon]\leq e^{-\epsilon^{2}/2nC^{2}}.
\]
\end{thm}

Consider a stationary process $Y\sim(A^{\mathbb{Z}},\mu)$ on finite
state space $A$, with a renewal state $a\in A$ of exponential hitting
times, i.e. the time $\tau_{0}=\min\{r>0:Y_{r}=a\}$ that takes to
hit $a$ satisfies for some $c_{1},c_{2}>0$ and any $n$,
\begin{equation}
||\mathbf{\mathbf{P}}[\tau_{0}>n|\varSigma(Y_{(-\infty,0]})]||_{\infty}\leq c_{1}e^{-c_{2}n}.\label{eq:-64}
\end{equation}
Important examples for such a process are any aperiodic Markov process
$X$ on a finite state space $A$, and the distribution $\widehat{X}^{a}$
of any of its elements.
\begin{prop}
\label{prop:Let--be-1}Let $Y\sim(A^{\mathbb{Z}},\mu)$ be a stationary
process with finite state space $A$ and a renewal state $a\in A$
of exponential hitting times (as in (\ref{eq:-64})). There exists
some $C>0$ such that for any $n$,
\[
|\frac{1}{n}\mathbf{E}[I(Y_{[0,n-1]})]-h(Y)|<\frac{C}{\sqrt{n}}.
\]
\end{prop}

\selectlanguage{american}%
\begin{proof}
Given $k>0$, we have
\begin{align}
\left|h(Y)-H(Y_{0}|Y_{[1,k]})\right| & =\left|H(Y_{0}|Y_{[1,\infty)})-H(Y_{0}|Y_{[1,k]})\right|\nonumber \\
 & =\left|\mathbf{E}[\mathbf{E}[I(Y_{0}|Y_{[1,\tau_{0}]})-I(Y_{0}|Y_{[1,k]})|\tau_{0}]]\right|\nonumber \\
 & \leq\mathbf{E}[(I(Y_{0}|Y_{[1,k]})+I(Y_{0}|Y_{[1,\tau_{0}]}))\cdot\mathbf{1}\{\tau_{0}>k\}]\nonumber \\
 & \leq\mathbf{P}[\tau_{0}>k]\cdot k\log|A|+\mathbf{E}[\tau_{0}\cdot\mathbf{1}\{\tau_{0}>k\}]\nonumber \\
 & \leq\mathbf{P}[\tau_{0}>k]\cdot k\log|A|+\sum_{n=k+1}^{\infty}n\mathbf{P}[\tau_{0}=n]\label{eq:-1-1}
\end{align}
and the exponential hitting times property (\foreignlanguage{english}{\ref{eq:-64}})
implies that for some $c_{1},c_{2}$, the last expression is bounded
by $c_{1}e^{-c_{2}k}$.

Note that
\begin{align*}
\mathbf{E}[I(Y_{[0,n-1]})] & =\mathbf{E}[I(Y_{[0,\sqrt{n}]})]+\sum_{k=\sqrt{n}+1}^{n}H(Y_{0}|Y_{[1,k]})\\
 & \leq\sqrt{n}H(Y_{1})+\sum_{k=\sqrt{n}+1}^{n}H(Y_{0}|Y_{[1,k]}).
\end{align*}
Thus

\begin{align*}
\left|\frac{1}{n}\mathbf{E}[I(Y_{[0,n-1]})]-h(Y)\right| & \leq\frac{1}{n}\left(\sqrt{n}(H(Y_{1})+h(Y))+\sum_{k=\sqrt{n}+1}^{n}\left|\mathbf{E}[I(Y_{0}|Y_{[1,k]})]-h(Y)\right|\right)\\
 & \leq(H(Y_{1})+h(Y))n^{-\frac{1}{2}}+c_{1}e^{-c_{2}n^{\nicefrac{1}{2}}}.
\end{align*}
which proves our claim.
\end{proof}
\selectlanguage{english}%
\begin{prop}
Let $Y\sim(A^{\mathbb{Z}},\mu)$ be a stationary process with finite
state space $A$ and a renewal state $a\in A$ of exponential hitting
times (as in (\ref{eq:-64})). With respect to the filtration $\left\{ \mathscr{F}_{i}\right\} _{i=1}^{\infty}$
defined by
\[
\mathscr{F}_{i}=\varSigma\left(Y_{k}:1\leq k\leq i\right)
\]
there exists $C>0$, such that for any $n$, the random information
$X:=I(Y_{[1,n]})$ satisfies inequality (\ref{eq:-5}).
\end{prop}

\begin{proof}
Fix $0\leq i<n$, and let $\tau_{i}=\min\{r>0:Y_{i+r}=a\}$. For each
$j\in\{i,i+1\}$, we have
\begin{align}
\mathbf{E}[I(Y_{[1,n]})|\mathscr{F}_{j}] & =\mathbf{E}[I(Y_{[1,i]})|\mathscr{F}_{j}]+\mathbf{E}[I(Y_{[1,i+\tau_{i}]}|Y_{[1,i]})|\mathscr{F}_{j}]+\mathbf{E}[Y_{[1,n]}|Y_{[i+1,i+\tau_{i}]})|\mathscr{F}_{j}]\nonumber \\
 & =\mathbf{E}[I(Y_{[1,i]})|\mathscr{F}_{j}]+\mathbf{E}[I(Y_{[i+1,i+\tau_{i}]})|\mathscr{F}_{j}]+\mathbf{E}[I(Y_{[i+\tau_{i}+1,n]})|\mathscr{F}_{j}]\label{eq:-65}
\end{align}
where the change in the second summand is derived from the fact that
$Y_{[1,1+i]}\in\mathscr{F}_{j}$, and the change in the third summand
is derived from the renewal property of $a$.

The first summand $I(Y_{[1,i]})$ is $\mathscr{F}_{i}$-measurable,
thus doesn't contribute to (\ref{eq:-5}):
\begin{equation}
\Vert\mathbf{E}[I(Y_{[1,i]})|\mathscr{F}_{i+1}]-\mathbf{E}[I(Y_{[1,i]})|\mathscr{F}_{i}]\Vert_{\infty}=0.\label{eq:-18}
\end{equation}
The second summand of (\ref{eq:-65}) is bounded by
\[
\mathbf{E}[I(Y_{[i+1,i+\tau_{i}]})|\mathscr{F}_{j}]\leq\mathbf{E}[\tau_{i}\log|A|\,|\mathscr{F}_{j}],
\]
which by the exponential return times property (\ref{eq:-64}) is
bounded by some $C>0$ that is independent of $i$ or $n$. Thus its
contribution to (\ref{eq:-5}) is at most
\[
\Vert\mathbf{E}[I(Y_{[i+1,i+\tau_{i}]})|\mathscr{F}_{i+1}]-\mathbf{E}[I(Y_{[i+1,i+\tau_{i}]})|\mathscr{F}_{i}]\Vert_{\infty}\leq2C.
\]

The third summand of (\ref{eq:-65}) can be estimated by
\begin{align*}
\mathbf{E}[I(Y_{[i+\tau_{i}+1,n]})|\mathscr{F}_{j}] & =\mathbf{E}[I(Y_{[i+\tau_{i}+1,n+\tau_{i}]})-I(Y_{[i+\tau_{i}+1,n+\tau_{i}]}|Y_{[i+\tau_{i}+1,n]})|\mathscr{F}_{j}]\\
 & =H(Y_{[1,n-i]}|Y_{0}=a)-\mathbf{E}[I(Y_{[i+\tau_{i}+1,n+\tau_{i}]}|Y_{[i+\tau_{i}+1,n]})|\mathscr{F}_{j}]\\
 & =H(Y_{[1,n-i]}|Y_{0}=a)-O(\mathbf{E}[\tau_{i}|\mathscr{F}_{j}])
\end{align*}
thus its contribution to (\ref{eq:-5}) boils down to
\begin{align*}
\Vert\mathbf{E}[I(Y_{[i+\tau_{i}+1,n]})|\mathscr{F}_{i+1}]-\mathbf{E}[I(Y_{[i+\tau_{i}+1,n]})|\mathscr{F}_{i}]\Vert_{\infty} & \leq O(\mathbf{E}[\tau_{i}|\mathscr{F}_{i}]+\mathbf{E}[\tau_{i}|\mathscr{F}_{i+1}])
\end{align*}
Another use of the exponential return times as in (\ref{eq:-64})
implies that the right hand side of the last inequality is of size
$O(1)$. This proves the claimed inequality.
\end{proof}
The last proposition enables us to apply Theorem \ref{thm:Azuma-Hoeffding-inequality},
which together with Proposition \ref{prop:Let--be-1} gives us the
following useful inequality:
\begin{cor}
\label{cor:Let--beunconditional=000020probability=000020Azuma}Let
$Y\sim(A^{\mathbb{Z}},\mu)$ be a stationary process on a finite state
space $A$, with a renewal state $a\in A$ of exponential hitting
times (as in (\ref{eq:-64})). There exist some $c,C\geq0$ s.t. for
any $n$ and any $\epsilon>0$,
\[
\mathbf{P}[|I(Y_{[1,n]})-\mathbf{E}[I(Y_{[1,n]})]|>\epsilon]\leq ce^{-\epsilon^{2}/2nC^{2}}.
\]
In particular, for any $\delta>\frac{1}{2}$, taking $\epsilon=\Lambda_{\delta,\frac{1}{2}}n$,
the following holds:
\[
\mathbf{P}[|I(Y_{[1,n]})-nh(Y)|>\Lambda_{\delta,\frac{1}{2}}n]=\Delta(n).
\]
\end{cor}

We use this result to estimate the conditional measure $\mu_{S}([X_{S}])$
of a likely filler in a filler space $A_{S}$:
\begin{cor}
\label{cor:Azuma=000020information}Suppose $X\sim(A^{\mathbb{Z}},\mu)$
is an aperiodic Markov process, $a,b\in A$. Taking $Z=X$ or $Z=\widehat{X}^{a}$
or $Z=\widehat{X}^{a,b}$, for any $\delta>\frac{1}{2}$ and any segment
$\mathcal{J}\subset\mathbb{Z}$, the following events hold with probability
$1-\Delta(|\mathcal{J}|)$: 
\begin{equation}
I(Z_{\mathcal{J}})\in|\mathcal{J}|h(Z)\pm\Lambda_{\delta,\frac{1}{2}}|\mathcal{J}|\label{eq:-11}
\end{equation}
and

\begin{equation}
I(Z_{\mathcal{J}}|\widehat{X}_{[\tau^{\leq\mathcal{J}},\tau^{\geq\mathcal{J}}]}^{a})\in|\mathcal{J}|h(Z|\widehat{X}^{a})\pm\Lambda_{\delta,\frac{1}{2}}|\mathcal{J}|\label{eq:-37}
\end{equation}
where $\tau^{\leq\mathcal{J}}:=\tau^{\leq\alpha},\tau^{\geq\mathcal{J}}:=\tau^{\geq\beta}$
for $[\alpha,\beta]=\mathcal{J}$ are the last (first) appearance
of $a$ before (after) $\mathcal{J}$.

In addition, given a nice marker factor $\xi=\xi(\widehat{X}^{a})$
with lengths sequence $\left(L_{r}\right)_{r\in\mathbb{N}}$, for
any $\epsilon>0$, with probability $1-\Delta(L_{r})$, the above
holds simultaneously for all intervals $\mathcal{J}\subset S_{r}$
of size $|\mathcal{J}|>|S_{r}|^{\epsilon}$.
\end{cor}

\begin{proof}
The estimate of the probability for the event of (\ref{eq:-11}) is
just Corollary \ref{cor:Let--beunconditional=000020probability=000020Azuma}.
To prove the statement for the event of (\ref{eq:-37}), fix $\frac{1}{2}>\epsilon>0$.
By the exponential hitting times property of $a$ in $X$, we have
that with probability $1-\Delta(|\mathcal{J}|)$, 
\begin{equation}
d(\mathcal{J},\tau^{\leq\mathcal{J}}),d(\mathcal{J},\tau^{\geq\mathcal{J}})<|\mathcal{J}|^{\epsilon}\label{eq:-38}
\end{equation}
We claim that on the event that (\ref{eq:-38}) holds, (\ref{eq:-37})
holds too, and consequently has also probability $1-\Delta(|\mathcal{J}|)$
to hold. To see that, assume (\ref{eq:-38}) holds, then we have
\[
I(Z_{\mathcal{J}}|\widehat{X}_{[\tau^{\leq\mathcal{J}},\tau^{\geq\mathcal{J}}]}^{a})=I(Z_{\mathcal{J}}|\widehat{X}_{\mathcal{J}\pm|\mathcal{J}|^{\epsilon}}^{a}).
\]

This gives a lower bound:
\[
I(Z_{\mathcal{J}})-I(\widehat{X}_{\mathcal{J}\pm|\mathcal{J}|^{\epsilon}}^{a})\leq I(Z_{\mathcal{J}}|\widehat{X}_{[\tau^{\leq\mathcal{J}},\tau^{\geq\mathcal{J}}]}^{a})
\]
and an upper bound:
\begin{align*}
I(Z_{\mathcal{J}}|\widehat{X}_{[\tau^{\leq\mathcal{J}},\tau^{\geq\mathcal{J}}]}^{a}) & =I(Z_{\mathcal{J}}|\widehat{X}_{\mathcal{J}\pm|\mathcal{J}|^{\epsilon}}^{a})\\
 & \leq I(Z_{\mathcal{J}\pm|\mathcal{J}|^{\epsilon}}|\widehat{X}_{\mathcal{J}\pm|\mathcal{J}|^{\epsilon}})\\
 & =I(Z_{\mathcal{J}\pm|\mathcal{J}|^{\epsilon}})-I(\widehat{X}_{\mathcal{J}\pm|\mathcal{J}|^{\epsilon}})
\end{align*}
(here $\mathcal{J}\pm|\mathcal{J}|^{\epsilon}=\{i\in\mathbb{Z}:d(\mathcal{J},i)\leq|\mathcal{J}|^{\epsilon}\}$).
Given $\delta>\frac{1}{2}$, take $\delta'\in(\frac{1}{2},\delta)$,
and use (\ref{eq:-11}) to get
\begin{align*}
I(Z_{\mathcal{J}}) & \in|\mathcal{J}|h(Z)\pm\Lambda_{\delta',\frac{1}{2}}(|\mathcal{J}|)\\
I(Z_{\mathcal{J}\pm|\mathcal{J}|^{\epsilon}}) & \in(|\mathcal{J}|+2|\mathcal{J}|^{\epsilon})h(Z)\pm\Lambda_{\delta',\frac{1}{2}}(|\mathcal{J}|+2|\mathcal{J}|^{\epsilon})\\
I(\widehat{X}_{\mathcal{J}\pm|\mathcal{J}|^{\epsilon}}) & \in(|\mathcal{J}|+2|\mathcal{J}|^{\epsilon})h(\widehat{X})\pm\Lambda_{\delta',\frac{1}{2}}(|\mathcal{J}|+2|\mathcal{J}|^{\epsilon})
\end{align*}
so that for large enough $|\mathcal{J}|$ (depending on $\delta,\delta'$
and $\epsilon$), we get
\begin{align*}
I(Z_{\mathcal{J}}|\widehat{X}_{[\tau^{\leq\mathcal{J}},\tau^{\geq\mathcal{J}}]}^{a}) & \in|\mathcal{J}|h(Z|\widehat{X})\pm2\Lambda_{\delta',\frac{1}{2}}(|\mathcal{J}|+2|\mathcal{J}|^{\epsilon})\pm O(|\mathcal{J}|^{\epsilon})\\
 & \subset|\mathcal{J}|h(Z|\widehat{X})\pm\Lambda_{\delta,\frac{1}{2}}(|\mathcal{J}|)
\end{align*}
which is just (\ref{eq:-37}).

Finally, the additional claim follows from the fact that with probability
$1-\Delta(L_{r})$, $S_{r}\subset[-L_{r},L_{r}]$, hence any of the
segments $\mathcal{J}\subset S_{r}$ in consideration must lie in
$[-L_{r},L_{r}]$. As there are no more than $4L_{r}^{2}$ such segments,
ranging through all of them and collecting the error sets associated
with each segment, which is of measure $\Delta(|S_{r}|^{\epsilon})=\Delta(L_{r})$,
implies the claim.
\end{proof}

\subsection{\protect\label{subsec:The-isomorphism-construction}The isomorphism
construction}

\subsubsection{\protect\label{subsec:Societies}Societies}

The next definition of a society and the proposition that follows
are fully taken from \cite{key-2}.

Let $U$ and $V$ be finite sets with probability measures $\rho$
and $\sigma$. A \textit{society} from $U$ to $V$ is a map
\[
S:U\rightarrow\text{subsets of }V
\]
 such that for any $B\subset U$,
\[
\rho(B)\leq\sigma(S(B))
\]
where
\[
S(B)=\bigcup_{b\in B}S(b)
\]

If $R$ and $S$ are societies, then $R\prec S$ if $R(b)\subset S(b)$
for all $b\in U$. If $S$ is a society, then we define 
\[
S^{*}:V\rightarrow\text{subsets of }U
\]
by the duality
\[
b\in S^{*}(g)\text{ if and only if }g\in S(b).
\]

\begin{prop}
\label{prop:Properties=000020of=000020societies}1. $S^{*}$ is a
society from $V$ to $U$.

2. (Marriage lemma) For any society $S$, there exists a society $R\prec S$
such that
\[
\#\left\{ g\in V:\,\exists b_{1}\neq b_{2}\text{ with }g\in R(b_{1})\cap R(b_{2})\right\} <|U|.
\]

3. Products of societies are societies.

4. $S^{**}=S$, and if $R\prec S$ then $R^{*}\prec S^{*}$
\end{prop}

For a proof of 1-3, see \cite[Lemma 7]{key-2}.

\textit{Proof of property 4.} Let $b\in U$. Then 
\[
g\in S^{**}(b)\iff b\in S^{*}(g)\iff g\in S(b)
\]
which shows that $S^{**}(g)=S(g)$. If $R\prec S$ and $g\in V$,
then
\[
b\in R^{*}(g)\Rightarrow g\in R(b)\Rightarrow g\in S(b)\Rightarrow b\in S^{*}(g)
\]
 so $R^{*}(g)\subset S^{*}(g)$.$\hfill\square$ 

\subsubsection{Keane-Smorodinsky coding}
\begin{lem}
\label{lem:good=000020systems=000020of=000020fillers=000020imply=000020good=000020isomorphism}Let
$X$ and $Y$ be two marked processes of equal entropy, equal marked-state
distribution, and a common nice marker factor $\xi$ exhibiting lengths
sequence $(L_{r})_{r\geq1}$ (see Def. \ref{def:nice=000020marker}).
Suppose both $X$ and $Y$ have good SOFs with respect to $\xi$.

There is then an isomorphism $X\overset{\phi}{\rightarrow}Y$ s.t.
for any $n$, there are $\{r_{i}\}_{i=1}^{4}$ with $L_{r_{i}}=\tilde{\Theta}(n)$,
satisfying with probability $1-\Delta(n)$
\[
[X_{S_{r_{1}}}]\Rrightarrow[\phi(X)_{S_{r_{2}}}]
\]
and
\[
[Y_{S_{r_{3}}}]\Rrightarrow[\phi^{-1}(Y)_{S_{r_{4}}}].
\]
\end{lem}

We give the construction of the map $\phi$, followed by the proof
that this coding indeed satisfies the assertion of Lemma \ref{lem:good=000020systems=000020of=000020fillers=000020imply=000020good=000020isomorphism}.
The construction, which is just an adaptation of the original Kean-Smorodinsky's
coding \cite{key-2} to our context, uses the notion of society and
the related marriage lemma (Proposition \ref{prop:Properties=000020of=000020societies}),
as defined in Section \ref{subsec:Societies}. 

\subsubsection*{Step 1: Defining societies}

We provide societies $\left\{ G_{S}\right\} _{S}$ between $\tilde{A}_{S}$
and $\tilde{B}_{S}$ for all $r$-skeletons $S$ and for all $r$.
For even rank skeletons, $G_{S}$ will be defined from $\tilde{B}_{S}$
to $\tilde{A}_{S}$, while for odd rank skeletons, $G_{S}$ will be
defined from $\tilde{A}_{S}$ to $\tilde{B}_{S}$. In the present
construction, we shall treat a skeleton $S$ as representing the class
$\{\sigma^{n}S:n\in\mathbb{Z}\}$ of all translates of the skeleton.
Accordingly, for a given skeleton $S$, the societies $G_{\sigma^{n}S}$
are all defined simultaneously as translates of each other:
\begin{equation}
G_{\sigma^{n}S}(U)=\sigma^{n}G_{S}(\sigma^{-n}U)\label{eq:-39}
\end{equation}

Start with thinning out $\xi$ if needed so that Proposition \ref{prop:key=000020property=000020s.o.f.}
will hold. The societies $\left\{ G_{S}\right\} _{S}$ are defined
inductively on the rank of $S$: For a $1$-skeleton $S$, let $F_{S}:\tilde{A}_{S}\rightarrow\tilde{B}_{S}$
be the trivial society defined for all $U\in\tilde{A}_{S}$ by
\[
F_{S}(U)=\tilde{B}_{S}
\]
and let $G_{S}\prec F_{S}$ be as in (2) of Proposition \ref{prop:Properties=000020of=000020societies}.
Now suppose that $S$ is an $r$-skeleton with $r$ even, where $S=S^{1}\cdots S^{m}$
is its rank $r-1$ decomposition, and that for each $1\leq i\leq m$,
a society $G_{S^{i}}:\tilde{A}_{S^{i}}\rightarrow\tilde{B}_{S^{i}}$
is already defined. Recall the notation
\begin{align*}
\bar{A}_{S}: & =\tilde{A}_{S^{1}}\times\cdots\times\tilde{A}_{S^{m}}\\
\bar{B}_{S}: & =\tilde{B}_{S^{1}}\times\cdots\times\tilde{B}_{S^{m}}.
\end{align*}
where by Property \ref{enu:For-any-skeleton} in Definition \ref{def:SOF}
of SOF, each element of $\bar{A}_{S}$ is a union of classes in $\tilde{A}_{S}$
(and similarly for $\bar{B}_{S}$ and $\tilde{B}_{S}$). We obtain
a society $F_{S}:\bar{B}_{S}\rightarrow\tilde{A}_{S}$ defined for
any $\times_{i=1}^{m}V_{i}$ with $V_{i}\in\tilde{B}_{S^{i}}$ by
\[
F_{S}(\times_{i=1}^{m}V_{i})=\times_{i=1}^{m}G_{S^{i}}^{*}(V_{i}).
\]
Let $G_{S}:\bar{B}_{S}\rightarrow\tilde{A}_{S}$ be a society with
$G_{S}\prec F_{S}$ as in (2) of Proposition \ref{prop:Properties=000020of=000020societies}.
Finally, as $\tilde{B}_{S}$ refines $\bar{B}_{S}$, we define $G_{S}$
for any $V\in\tilde{B}_{S}$ by
\[
G_{S}(V)=G_{S}(V')
\]
where $V'$ is the unique element of $\bar{B}_{S}$ that contains
$V$.

Similarly, for odd $r$, suppose we have defined for each $i$ a society
$G_{S_{i}}:\tilde{B}_{S^{i}}\rightarrow\tilde{A}_{S^{i}}$. We obtain
a society $F_{S}:\bar{A}_{S}\rightarrow\tilde{B}_{S}$ defined for
any $\times_{i=1}^{m}U_{i}\in\bar{A}_{S}$ by
\[
F_{S}(\times_{i=1}^{m}U_{i})=\times_{i=1}^{m}G_{S^{i}}^{*}(U_{i}).
\]
Let $G_{S}:\bar{A}_{S}\rightarrow\tilde{B}_{S}$ be a society with
$G_{S}\prec F_{S}$ as in (2) of Proposition \ref{prop:Properties=000020of=000020societies},
and define $G_{S}$ on $\tilde{A}_{S}$ by
\[
G_{S}(U)=G_{S}(U')
\]
where $U'$ is the unique element of $\bar{A}_{S}$ that contains
$U$.

\subsubsection*{Step 2: Defining the isomorphism. }

Using the societies provided in step 1, we define the maps $\phi:\mathbf{A}\rightarrow\mathbf{B}$
and $\psi:\mathbf{B}\rightarrow\mathbf{A}$ (which will be shown to
satisfy $\psi=\phi^{-1}$). Let $X\sim\mathbf{A}$, fix some even
$r$, and let $S=S_{r}(X)$ with random rank $r-1$ decomposition
$S=S^{1}\ast\cdots\ast S^{m}$. By the construction of $G_{S}:\bar{B}_{S}\rightarrow\tilde{A}_{S}$,
the number of classes $U\in\tilde{A}_{S}$ with no unique $V\in\bar{B}_{S}$
for whom $U\subset G_{S}(V)$, is at most
\[
\#\left\{ U\in\tilde{A}_{S}:\exists V_{1},V_{2}\in\bar{B}_{S}\text{ s.t. }U\in G_{S}(V_{1})\cap G_{S}(V_{2})\right\} \leq|\bar{B}_{S}|
\]
Conditioned on $S$, the probability of $[X_{S}]$ not being one of
these classes is at most
\begin{align}
\mathbf{P}[\exists!V\in\bar{B}_{S}:\,[X_{S}]\in G_{S}(V)|S] & \geq1-\max\{\mu_{S}(\bigcup_{i=1}^{|\bar{B}_{S_{r}}|}U_{i}):\,U_{i}\in\tilde{A}_{S}\}\label{eq:Un=000020uniquely=000020associated=000020classes,=000020estimate}\\
 & =1-N_{\mathbf{A}}(S,|\bar{B}_{S}|)\nonumber 
\end{align}
Taking expectations on both sides, one gets by Proposition \ref{prop:key=000020property=000020s.o.f.}:
\begin{align}
\mathbf{P}[\exists!V\in\bar{B}_{S}:\,[X_{S}]\in G_{S}(V)] & \geq1-\mathbf{E}[N_{\mathbf{A}}(S,|\bar{B}_{S}|)]\label{eq:Prob.=000020for=000020unique=000020asignment}\\
 & =1-\Delta(L_{r})\nonumber 
\end{align}
and in particular this converges to 1 as $r\rightarrow\infty$. Thus
a.s. there is an infinite increasing sequence of (even) integers $\left(r_{n}\right)_{n}$
s.t. for each $n$, there is a unique element $V\in\bar{B}_{S_{r_{n}}}$
with $[X_{S_{r_{n}}}]\in G_{S_{r_{n}}}(V)$. Equivalently, $[X_{S_{r_{n}}}]$
is sent under the dual of $G_{S_{r_{n}}}$ to a single element $G_{S_{r_{n}}}^{*}([X_{S_{r_{n}}}])\in\bar{B}_{S_{r_{n}}}$
(as opposed to a union of elements).

For each $n$, identify the elements $V\in\bar{B}_{S_{r_{n}}}$ with
the sets $\{y\in B^{\mathbb{Z}}:y_{S_{r_{n}}(x)}\in V\}$. With this
identification, $\left(G_{S_{r_{n}}}^{*}([X_{S_{r_{n}}}])\right)_{n}$
is an inclusion-decreasing sequence of non-empty sets, compact w.r.t.
the product topology of $B^{\mathbb{Z}}$ (each one of them being
a finite union of cylinder sets), thus their intersection $\bigcap_{n=1}^{\infty}G_{S_{r_{n}}}^{*}([X_{S_{r_{n}}}])$
is non-empty. Pick some $y\in\bigcap_{n=1}^{\infty}G_{S_{r_{n}}}^{*}([X_{S_{r_{n}}}])$
and define $\phi(x):=y$. 

Similarly, there is a.s. an infinite increasing sequence of (odd)
integers $\left(r_{n}\right)_{n}$ s.t. for each $n$, $G_{S_{r_{n}}(Y)}^{*}([Y_{S_{r_{n}}(Y)}])$
is a single element of $\bar{A}_{S_{r_{n}}(Y)}$, and viewing these
elements as sets in $A^{\mathbb{Z}}$, their intersection $\bigcap_{n=1}^{\infty}G_{S_{r_{n}}(Y)}^{*}([Y_{S_{r_{n}}(Y)}])$
is non-empty. We pick some $x$ from this intersection and define
$\psi(y):=x$. 

Note that a.s. $S_{r}(X)=S_{r}(\phi(X))$ for infinitely many $r$,
thus $\widehat{X}=\widehat{\phi(X)}$ a.s. Similarly, we have $\hat{Y}=\widehat{\psi(Y)}$
a.s.

We claim that the map $\phi$ satisfies the conclusion of Lemma \textit{\ref{lem:good=000020systems=000020of=000020fillers=000020imply=000020good=000020isomorphism}.}

\textit{Proof of Lemma \ref{lem:good=000020systems=000020of=000020fillers=000020imply=000020good=000020isomorphism}:}
To prove the lemma, we claim that the map $\phi$ is measure preserving,
stationary, with inverse $\phi^{-1}=\psi$, and satisfies the determination
property asserted in the lemma. We divide the proof accordingly.

\subsubsection*{$\phi$ is measure preserving:}

Fix $\omega\in\{0,1\}^{\mathbb{Z}}$ generic for $\widehat{X}$ and
condition $X$ on the event $\widehat{X}=\omega$ (thus $S_{r}=S_{r}(X)=S_{r}(\omega)$
is fixed for all $r$). Let $r_{0}\in\mathbb{N}$ and $V_{0}\in\tilde{B}_{S_{r_{0}}}$,
then $\phi(X)\in V_{0}$ only if for some even $r>r_{0}$, $G_{S_{r}}^{*}([X_{S_{r}}])$
is a single element of $\bar{B}_{S_{r}}$, and as a set in $B^{\mathbb{Z}}$
is contained in $V_{0}$. Thus for any even $r\geq r_{0}$, $\phi^{-1}V_{0}$
is included in the union $D_{r}\cup E_{r}$, where
\begin{align*}
D_{r} & =\cup\{U\in\tilde{A}_{S_{r}}:G_{S_{r}}^{*}(U)\subset V_{0})\}\\
E_{r} & =\cup\{U\in\tilde{A}_{S_{r}}:G_{S_{r}}^{*}(U)\text{ is not a single element in }\bar{B}_{S_{r}}\}
\end{align*}

Thus 
\begin{align*}
\mu_{\omega}(\phi^{-1}V_{0}) & \leq\liminf_{r\text{ even}}(\mu_{\omega}(D_{r})+\mu_{\omega}(E_{r}))\\
 & \leq\sup_{r\text{ even}}\nu_{\omega}(G_{S_{r}}^{*}(D_{r}))+\liminf_{r\text{ even}}N_{\mathbf{A}}(S_{r}(\omega),|\bar{B}_{S_{r}(\omega)}|)\\
 & \overset{{\scriptscriptstyle \omega-a.s.}}{=}\nu_{\omega}(V_{0})
\end{align*}
As $V_{0}$ was chosen arbitrarily, we have
\[
\phi_{*}\mu_{\omega}|_{\tilde{B}_{S_{r_{0}}}}\leq\nu_{\omega}|_{\tilde{B}_{S_{r_{0}}}}
\]
Having in both sides probability measures, the inequality becomes
equality. Since $r$ and $\omega$ are arbitrary, $\phi$ preserves
measure on a family of sets that approximates the whole Borel $\sigma$-algebra
of $B^{\mathbb{Z}}$ (Proposition \ref{prop:System=000020of=000020fillers=000020generates=000020the=000020algebra}),
thus $\phi$ is measure preserving.

\subsubsection*{$\phi$ is stationary:}

By the definition of a marker process, a.s. for large enough $r$,
$S_{r}$ covers $X$ at $[0,1]$. This implies that for large enough
$r$, $S_{r}(\sigma X)=\sigma S_{r}$. That, together with (\ref{eq:-39}),
implies that if $\left(r_{n}\right)$ is an increasing sequence of
integers such that for all $n$, $G_{S_{r_{n}}}^{*}([X_{S_{r_{n}}}])$
is a single element in $\bar{B}_{S_{r_{n}}}$, then the same holds
for $G_{S_{r_{n}}(\sigma X)}^{*}([(\sigma X)_{S_{r_{n}}(\sigma X)}])$
in $\bar{B}_{S_{r_{n}}(\sigma X)}$, for all $n$ large enough, and
\begin{align*}
\phi(\sigma X) & \in\bigcap_{n=1}^{\infty}G_{S_{r_{n}}(\sigma X)}^{*}([(\sigma X)_{S_{r_{n}}(\sigma X)}])=\sigma\left(\bigcap_{n=1}^{\infty}G_{S_{r_{n}}}^{*}([X_{S_{r_{n}}}])\right)
\end{align*}
Hence, in order to show that $\sigma^{-1}\phi(\sigma x)=\phi(x)$,
it is enough to show that a.s. $\bigcap_{n=1}^{\infty}G_{S_{r_{n}}}^{*}([X_{S_{r_{n}}}])$
is a singleton. Notice that for any $n$, 
\[
G_{S_{r_{n}}}^{*}([X_{S_{r_{n}}}])=[\phi(X)_{S_{r_{n}}}]_{\bar{B}_{S_{r_{n}}}}\subset[\phi(X)_{S_{r_{n}-1}}].
\]
By Proposition \ref{prop:System=000020of=000020fillers=000020generates=000020the=000020algebra},
$\left(\mathcal{B}_{r_{n}-1}\right)_{n}$ generates the Borel $\sigma$-algebra
of $B^{\mathbb{Z}}$ modulo $\nu$. We already proved that $\phi$
preserves measure, thus $X$-a.s. 
\[
\bigcap_{n=1}^{\infty}G_{S_{r_{n}}}^{*}([X_{S_{r_{n}}}])\subset\bigcap_{n=1}^{\infty}[\phi(X)_{S_{r_{n}}}]=\{\phi(X)\}
\]
and stationarity of $\phi$ follows.

\subsubsection*{$\psi$ is the inverse of $\phi$:}

By the measure preserving property of $\phi$, a.s. both $\phi(X)$
and $\psi(\phi(X))$ are defined. Thus we have sequences $\left(r_{n}\right)_{n}$
and $\left(r'_{n}\right)_{n}$ so that 
\begin{align*}
\{\phi(X)\} & =\bigcap_{n=1}^{\infty}G_{S_{r_{n}}}^{*}([X_{S_{r_{n}}}])\\
\{\psi(\phi(X))\} & =\bigcap_{n=1}^{\infty}G_{S_{r'_{n}}}^{*}([\phi(X)_{S_{r'_{n}}}])
\end{align*}
(notice that $S_{r}(\phi(X))=S_{r}$ for all $r$). By passing to
a subsequence, we can assume that $r_{n}\geq r'_{n}$. Writing $V_{n}=G_{S_{r_{n}}}^{*}([X_{S_{r_{n}}}])$,
we have $[X_{S_{r_{n}}}]\subset G_{S_{r_{n}}}(V_{n})$, and $V_{n}\subset[\phi(X)_{S_{r_{n}-1}}]$.
Thus
\begin{equation}
[X_{S_{r_{n}}}]\subset G_{S_{r_{n}}}(V_{n})\subset G_{S_{r_{n}}}([\phi(X)_{S_{r_{n}-1}}])\label{eq:-40}
\end{equation}
Writing the rank $r_{n}-1$ decomposition of $S_{r_{n}}$ as $S_{r_{n}}=S^{1}*\cdots*S_{r_{n}-1}*\cdots*S^{m}$
where $S_{r_{n}-1}$ is the $(r_{n}-1)$-skeleton that covers $X$
at $0$, we have 
\[
G_{S_{r_{n}}}\prec\times_{i=1}^{m}G_{S^{i}}^{*}\prec G_{S_{r_{n}-1}}^{*}
\]
thus 
\begin{equation}
G_{S_{r_{n}}}([\phi(X)_{S_{r_{n}-1}}])\subset G_{S_{r_{n}-1}}^{*}([\phi(X)_{S_{r_{n}-1}}])\label{eq:-41}
\end{equation}
Using property 4 of Proposition \ref{prop:Properties=000020of=000020societies},
it can be shown by induction that $r_{n}-1\geq r'_{n}$ implies $G_{S_{r_{n}-1}}^{*}\prec G_{S_{r'_{n}}}^{*}$.
Thus
\begin{equation}
G_{S_{r_{n}-1}}^{*}([\phi(X)_{S_{r_{n}-1}}])\subset G_{S_{r_{n}-1}}^{*}([\phi(X)_{S_{r'_{n}}}])\subset G_{S_{r'_{n}}}^{*}([\phi(X)_{S_{r'_{n}}}])\label{eq:-42}
\end{equation}
and overall we get from (\ref{eq:-40}), (\ref{eq:-41}) and (\ref{eq:-42})
that
\[
[X_{S_{r_{n}}}]\subset G_{S_{r'_{n}}}^{*}([\phi(X)_{S_{r'_{n}}}])
\]
and so,
\[
\{\psi(\phi(X))\}=\bigcap_{n=1}^{\infty}G_{S_{r'_{n}}}^{*}([\phi(X)_{S_{r'_{n}}}])\subset\bigcap_{n=1}^{\infty}[X_{S_{r_{n}}}]=\{x\}
\]
Which shows that $\psi\circ\phi=Id_{\mathbf{A}}$. In the same manner
one shows that $\phi\circ\psi=Id_{\mathbf{B}}$, thus $\psi=\phi^{-1}$.

\subsubsection*{Determination property: }

Finally, to show the last assertion of the lemma, let $n\in\mathbb{N}$,
and let $r$ be an even integer for which $L_{r}=\tilde{\Theta}(n)$.
$\xi$ being a nice marker process,  such an $r$ indeed exists, and
we further have $L_{r-1}=\tilde{\Theta}(n)$. We have
\[
[X_{S_{r}}]\Rrightarrow G_{S_{r}}^{*}([X_{S_{r}}])
\]
and by (\ref{eq:Prob.=000020for=000020unique=000020asignment}), $G_{S_{r}}^{*}([X_{S_{r}}])\in\bar{B}_{S_{r}}$
is a single element with probability $1-\Delta(L_{r})$, and in that
case $G_{S_{r}}^{*}([X_{S_{r}}])\subset[\phi(X)_{S_{r-1}}]$, thus
\[
\mathbf{P}\left[[X_{S_{r}}]\Rrightarrow G_{S_{r}}^{*}([X_{S_{r}}])\Rrightarrow[\phi(X)_{S_{r-1}}]\right]=1-\Delta(L_{r}).
\]
using the same reasoning to $Y$ and $\phi^{-1}(Y)$, and the fact
that $\Delta(L_{r})=\Delta(n)$ yields the claimed determination property,
and completes the proof of the lemma.$\hfill\square$

For the above lemma to be applied, one needs to provide good SOFs
for the processes in consideration. The next proposition provides
such good SOFs for any marked process. These good SOFs satisfy an
additional property, under which the ``determination property''
asserted at the end of Lemma \ref{lem:good=000020systems=000020of=000020fillers=000020imply=000020good=000020isomorphism}
interprets to a bound on the coding radius of the isomorphism (for
more details, see the proof of Theorem \textit{\ref{thm:common=000020entry=000020version}}
at the end of this section).

\subsection{Construction of Good SOF}
\begin{prop}
\label{prop:Existence=000020of=000020good=000020filler=000020systems}For
any marked process $X\sim\mathbf{A}$ with nice marker process $\xi$
and lengths sequence $\left(L_{r}\right)_{r}$, there exists a good
SOF that additionally satisfies for some $\delta>0$:
\begin{equation}
\mathbf{P}\left[\frac{1}{|S_{r}|}\sum_{i\in S_{r}}\mathbf{1}\{[X_{S_{r}}]\Rrightarrow X_{i}\}\geq1-\Lambda_{\delta,\frac{1}{2}}|S_{r}|\right]=1-\Delta(L_{r}).\label{eq:good=000020expected=000020number=000020of=000020determined=000020coordinates}
\end{equation}
\end{prop}

In Sec. \ref{subsec:Proof-outline-of} we give an outline of the good-SOF
construction and the proof of the above proposition. The detailed
construction is given in Sec. \ref{subsec:Construction-of-a}, followed
by the proof of the proposition in Sec. \ref{subsec:Proof-of-Proposition}.

\subsubsection{\protect\label{subsec:Proof-outline-of}Proof outline of Proposition
\ref{prop:Existence=000020of=000020good=000020filler=000020systems}.}

For simplicity, let us consider first $X$ as an i.i.d. process, with
marginal measure $p$ on $A$. Given a skeleton $S$, the class of
a word $w\in A_{S}$ will be defined by its name along a certain subset
of indices $\mathcal{J}(w)\subset S$, so that
\[
w'\in[w]\iff w'_{\mathcal{J}(w)}=w_{\mathcal{J}(w)}.
\]
The subset $\mathcal{J}(w)$ will be determined by the values of $w$
along $\mathcal{J}(w)$, so that if $w'\in A_{S}$ and $w'_{\mathcal{J}(w)}=w{}_{\mathcal{J}(w)}$
then $\mathcal{J}(w')=\mathcal{J}(w)$. 

In more detail, we begin by defining $\mathcal{J}(w)=\emptyset$ for
all $w\in A_{S}$ with $\text{mr}(S)=1$. For $w\in A_{S}$ with $\text{mr}(S)>1$,
we define $\mathcal{J}(w)$ inductively: 

Fix some $\delta>\frac{1}{2}$, and let $S=S^{1}\ast\cdots\ast S^{m}$
be the $(\text{mr}(S)-1)$-decomposition of $S$ and accordingly $w=w^{1}*\cdots*w^{m}$.
Assuming $\mathcal{J}(w^{i})\subset S^{i}$ are already defined for
all $1\leq i\leq m$, let 
\begin{equation}
\mathcal{J}(w)=[\alpha,\gamma]\cup\bigcup_{i=1}^{m}\mathcal{J}(w^{i})\label{eq:-59}
\end{equation}
where $\alpha$ is the smallest index covered by $S$, and $\gamma$
is the smallest index for which the class $[w]$ will satisfy
\begin{equation}
I([w]|\widehat{wa})\geq|S|h(X|\widehat{X})-\Lambda_{\delta,\frac{1}{2}}|S|\label{eq:-12-1}
\end{equation}
where in case there is no such $\gamma$, we put $\mathcal{J}(w)=S$. 

Corollary \ref{cor:Azuma=000020information} guarantees that with
probability $1-\Delta(|S|)$, there is such $\gamma$ so that (\ref{eq:-12-1})
will take place, which is just property \ref{enu:Good=000020SOF=0000202}
of good SOF (Definition \ref{def:Good=000020SOF}). On the other hand,
the fact that $\gamma-1$ does not satisfy (\ref{eq:-12-1}) implies
property \ref{enu:For-any-large} of good SOF. To obtain (\ref{eq:good=000020expected=000020number=000020of=000020determined=000020coordinates}),
notice that 
\[
|S\backslash\mathcal{J}(w)|\cdot\min_{a'\in A}\log(p(a')^{-1})\leq I(w|\widehat{wa})-I([w]|\widehat{wa}).
\]
The upper bound in Corollary \ref{cor:Azuma=000020information} can
be used to imply that with probability $1-\Delta(|S|)$,
\[
I(w|\widehat{wa})-I([w]|\widehat{wa})\leq\Lambda_{\delta,\frac{1}{2}}|S|
\]
thus the additional property (\ref{eq:good=000020expected=000020number=000020of=000020determined=000020coordinates})
holds and the proof is completed for the case that $X$ is an i.i.d.
process.

When applying the same method to the general case of any marked process
$X$, we require in (\ref{eq:-59}) that $\gamma$ will be the smallest
index that satisfy both (\ref{eq:-12-1}) and $w_{\gamma+1}=a$. In
other words, letting $\alpha=i_{0}<\cdots<i_{n}=\beta+1$ be all coordinates
of $S$ with the symbol $a$, and denoting $N_{j}=[i_{j-1},i_{j}-1]$
for $j=1...n$, we obtain $\mathcal{J}(w)$ by adding to $\bigcup_{i=1}^{m}\mathcal{J}(w^{i})$
intervals $N_{j}$ one-by-one until $\mathcal{J}(w)$ finally satisfies
(\ref{eq:-12-1}). Using the reduction property of $X$ and the fact
that the state space is finite, the probability of any such interval
can be estimated by
\begin{equation}
\min_{(a',a'')}\mathbf{P}[X_{0}=a'|X_{-1}=a'']^{i_{j}-i_{j-1}}\leq\mathbf{P}[X_{N_{j}}=w_{N_{j}}|\widehat{X}_{[\alpha,\beta]}=S]\leq c\label{eq:-60}
\end{equation}
where $c<1$ is the probability reduction constant (\ref{eq:Marked=000020state=000020property}),
which is independent of the size of $N_{j}$. Letting $\gamma=i_{j}-1,$the
fact that (\ref{eq:-12-1}) doesn't hold for $[\alpha,i_{j-1}]\cup\bigcup_{i=1}^{m}\mathcal{J}(w^{i})$
leads to the following upper bound
\[
I([w]|\widehat{w})\leq|S|h(X|\widehat{X})-\Lambda_{\delta,\frac{1}{2}}|S|+(i_{j}-i_{j-1})\cdot\log\left(\min_{(a',a'')}\mathbf{P}[X_{0}=a'|X_{-1}=a'']^{-1}\right).
\]
To obtain property \ref{enu:For-any-large} of good SOF from that
upper bound, we need to assume that $(i_{j}-i_{j-1})=O(\Lambda_{\delta,\frac{1}{2}}|S|)$.
The exponential hitting times property guarantees this is true with
probability $1-\Delta(|S|)$, thus the definition of $\mathcal{J}$
as above applies only for that part of the space (on the rest of the
probability space we define $\mathcal{J}(w)$ as the union of $\mathcal{J}(w_{i})$
only). The constant $c$ in (\ref{eq:-60}) is then used to estimate
the number of symbols that are determined by the class $[X_{S}]$,
yielding the validity of (\ref{eq:good=000020expected=000020number=000020of=000020determined=000020coordinates}).$\hfill\blacksquare$

\subsubsection{\protect\label{subsec:Construction-of-a}Construction of a good SOF
$\mathcal{A}$}

Let 
\begin{align*}
\theta & =\min_{(a',a'')}\mathbf{P}[X_{0}=a'|X_{-1}=a'']
\end{align*}
where $(a',a'')$ ranges over all tuples with $\mathbf{P}[X_{0}=a'|X_{-1}=a'']>0$.
Fix $L_{0}$ so that for any $L\geq L_{0}$
\[
\Lambda_{1,\frac{1}{2}}L-\log^{2}L\cdot\log\theta^{-1}\geq\sqrt{L}.
\]
Given a skeleton $S=\left(\widehat{x}_{[\alpha,\beta]},\xi(x)_{[\alpha,\beta]}\right)$,
define for each word $w\in A_{S}$ a subset of coordinates $\mathcal{J}(w)\subset S$
as follows: 

\textbf{Case 0 -} $|S|<L_{0}$ or $r=1$: Put
\[
\mathcal{J}(w)=\emptyset.
\]
Suppose now that $|S|\geq L_{0}$, write $S=S^{1}\cdots S^{m}$ and
accordingly $w=w^{1}\cdots w^{m}$ for its $(\text{mr}(S)-1)$-decomposition,
and assume $\mathcal{J}(w^{i})\subset S^{i}$ are all defined (or
alternatively, assuming $\text{mr}(S)=1$ and letting $m=0$). Write
\[
\mathcal{J}'(w)=\bigcup_{i=1}^{m}\mathcal{J}(w^{i})
\]
(which are empty when $\text{mr}(S)=1$) and consider two cases:

\textbf{Case 1 -} $|S|\geq L_{0}$ and $S$ contains an interval of
size at least $\log^{2}|S|$ with no $a$ (the ``rare'' case, see
(\ref{eq:-47}) below): Put
\[
\mathcal{J}(w)=\mathcal{J}'(w).
\]

\textbf{Case 2 -} $|S|\geq L_{0}$ and $S$ does not contain an $a$-free
interval of size $>\log^{2}|S|$: Define $\mathcal{J}(w)$ to be
\[
\mathcal{J}(w)=[\alpha,\gamma]\cup\mathcal{J}'(w)
\]
where $[\alpha,\gamma]$ is the smallest prefix of $S$ with $z_{\gamma+1}=a$
for which
\begin{equation}
I(w_{\mathcal{J}(w)}|\widehat{wa})\geq|S|h(X|\widehat{X})-\Lambda_{1,\frac{1}{2}}|S|\label{eq:-21}
\end{equation}
(where $w_{\mathcal{J}(w)}$ and $\widehat{wa}$ are considered as
values of $X_{\mathcal{J}(w)}$ and $\widehat{X}_{[\alpha,\beta+1]}$
respectively), and if there is no such a prefix, put $\mathcal{J}(w)=S$.
Now define for any skeleton $S$ the equivalence relation $\sim$
on $A_{S}$ by
\[
w\sim w'\iff w_{\mathcal{J}(w)}=w'_{\mathcal{J}(w)}.
\]
It can be easily observed that if $w\sim w'$ then $\mathcal{J}(w)=\mathcal{J}(w')$,
thus $\sim$ is indeed an equivalence relation. Finally, we let $\tilde{A}_{S}$
be the set of equivalence classes of $\sim$ in $A_{S}$. Denote the
collection $\{(\tilde{A}_{S},\mu_{S})_{S}:S\text{ a \ensuremath{\xi}-skeleton}\}$
by $\mathscr{A}.$ To establish Proposition \ref{prop:Existence=000020of=000020good=000020filler=000020systems},
we show that $\mathscr{A}$ is a good SOF that satisfies the additional
property asserted there.

\subsubsection{\textit{\protect\label{subsec:Proof-of-Proposition}Proof of Proposition
\ref{prop:Existence=000020of=000020good=000020filler=000020systems}:}}

It is clear from the definition that $\mathscr{A}$ satisfies properties
\ref{enu:The-collection-is} and \ref{enu:For-any-skeleton} of SOF
(Definition \ref{def:SOF}). The validity of the third property will
follow once we prove that Property \ref{enu:For-,-the} of good SOF
is valid (Definition \ref{def:Good=000020SOF}). Thus we turn now
to show it is a good SOF.

We first show that Property \ref{enu:For-any-large} of good SOF (Definition
\ref{def:Good=000020SOF}) holds with $k=0$, i.e. that for any skeleton
$S$ of large enough size $|S|$, any class $[w]\in\tilde{A}_{S}$
satisfies
\begin{equation}
I_{\mu_{S}}([w])\leq|S|h(X|\widehat{X})-\sqrt{|S|}.\label{eq:-31}
\end{equation}

There are three cases:

\textbf{Case 0 -} $|S|<L_{0}$: then there is only one class of measure
one, so that $I_{\mu_{S}}([w])=0$ and (\ref{eq:-31}) trivially follows. 

\textbf{Case 1 -} $|S|\geq L_{0}$ and $S$ contains an $a$-free
interval of size $>\log^{2}|S|$: If $\text{mr}(S)=1$ then again
we get $I_{\mu_{S}}([w])=0$ for all fillers and (\ref{eq:-31}) trivially
holds. Otherwise, writing $S=S^{1}\ast\cdots\ast S^{m}$ as in the
construction of $\mathscr{A}$, then as $\mathcal{J}(w)=\mathcal{J}'(w)$,
we have
\[
[w]=[w^{1}]\times\cdots\times[w^{m}].
\]
Assuming (\ref{eq:-31}) holds for skeletons of minimal rank $\text{mr}(S)-1$,
then
\[
I([w]|\hat{w})=\sum_{i}I([w^{i}]|\widehat{w^{i}})\leq|S|h(X|\widehat{X})-\sum_{i}\sqrt{|S^{i}|}
\]
which is at most $|S|h(X|\widehat{X})-\sqrt{|S|}$ by Corollary \ref{cor:For-any-,}.
Thus (\ref{eq:-31}) holds for $\tilde{C}_{S}$.

\textbf{Case 2} - $|S|\geq L_{0}$ and $S$ does not contain an $a$-free
interval of size $>\log^{2}|S|$: and suppose that (\ref{eq:-21})
takes place for some prefix (otherwise the claim follows immediately).
Writing $\mathcal{J}(w)=[\alpha,\gamma]\cup\bigcup_{i}\mathcal{J}(w^{i})$
, one has
\[
|S|h(X|\widehat{X})-\Lambda_{1,\frac{1}{2}}|S|>I(w_{[\alpha,\tau^{<\gamma}]\cup\bigcup_{i}\mathcal{J}(w^{i})}|\hat{w}).
\]
 By the assumption that $S$ doesn't contain an interval of size $>\log^{2}|S|$
with no $a$, we have $\gamma-\tau^{<\gamma}\leq\log^{2}|S|$, and
so,
\begin{align*}
I(w_{\mathcal{J}(w)}|\widehat{wa}) & \leq I(w_{[\alpha,\gamma]\cup\bigcup_{i}\mathcal{J}(w^{i})}|\hat{w})\\
 & \leq I(w_{[\alpha,\tau^{<\gamma}]\cup\bigcup_{i}\mathcal{J}(w^{i})}|\hat{w})+I(w_{[\tau^{<\gamma},\gamma]}|\hat{w}_{[\tau^{<\gamma},\gamma]})\\
 & \leq|S|h(X|\widehat{X})-\Lambda_{1,\frac{1}{2}}|S|+\log^{2}|S|\cdot\log\theta^{-1}\\
 & \leq|S|h(X|\widehat{X})-\sqrt{|S|}
\end{align*}
where the last inequality holds by choosing large enough $L_{0}$
and assuming $|S|\geq L_{0}$. Overall we see that Property \ref{enu:For-any-large}
of good SOF (Definition \ref{def:Good=000020SOF}) holds with $k=0$.

We now show that Property \ref{enu:Good=000020SOF=0000202} of good
SOF holds: From the definition of $\mathcal{J}$, if a skeleton $S$
that appears in $X$ doesn't contain an interval free of $a$'s of
size $>\log^{2}|S|$, and in addition
\begin{equation}
I(X_{S}|\widehat{X})\geq|S|h(X|\widehat{X})-\Lambda_{1,\frac{1}{2}}|S|\label{eq:-13}
\end{equation}
then (\ref{eq:-21}) takes place for some prefix $[\alpha,\gamma]$,
and the conditional information of the whole class of $X_{S}$ also
satisfies
\begin{equation}
I([X_{S}]|\widehat{X})\geq|S|h(X|\widehat{X})-\Lambda_{1,\frac{1}{2}}|S|.\label{eq:-14}
\end{equation}
Writing
\begin{align*}
E_{r} & =\{\forall N=[i,i+\log^{2}|S_{r}|]\subset S_{r},\exists j\in N:X_{j}=a\}\\
F_{r} & =\{I(X_{S_{r}}|\widehat{X})\in|S_{r}|h(X|\widehat{X})\pm\Lambda_{1,\frac{1}{2}}|S_{r}|\}\\
G_{r} & =\{I([X_{S_{r}}]|\widehat{X})\geq|S_{r}|h(X|\widehat{X})-\Lambda_{1,\frac{1}{2}}|S_{r}|\}
\end{align*}
the exponential hitting times property of the state $a$ (together
with $\xi$ being nice) implies that 
\begin{equation}
\mathbf{P}[E_{r}]=1-\Delta(L_{r}).\label{eq:-47}
\end{equation}
Moreover, by Corollary \ref{cor:Azuma=000020information}, 
\[
\mathbf{P}[F_{r}]=1-\Delta(L_{r})
\]
On the event $E_{r}\cap F_{r}$, (\ref{eq:-13}) is valid and implies
the validity of (\ref{eq:-14}), thus
\[
\mathbf{P}[G_{r}]\geq\mathbf{P}[E_{r}\cap F_{r}]=1-\Delta(L_{r})
\]
 and $\mathscr{A}$ satisfies Property \ref{enu:Good=000020SOF=0000202}
of good SOF.

It remains to show that (\ref{eq:good=000020expected=000020number=000020of=000020determined=000020coordinates})
holds, i.e. that
\[
\mathbf{P}\left[\frac{1}{|S_{r}|}\sum_{i\in S_{r}}\mathbf{1}\{[X_{S_{r}}]\Rrightarrow X_{i}\}\geq1-\tilde{O}(|S_{r}|^{-\frac{1}{2}})\right]=1-\Delta(L_{r}).
\]
 Let $\alpha=i_{0}<\cdots<i_{n}=\beta+1$ be all coordinates of $X_{S_{r}}$
with the symbol $a$, where as before, $[\alpha,\beta]$ denotes the
segment covered by $S_{r}$, and let $N_{j}=[i_{j-1},i_{j}-1]$ for
$j=1...n$. From the construction above it is clear that $\mathcal{J}(X_{S_{r}})=\bigcup_{j\in J}N_{j}$
for some set $J\subset\{1,...,n\}$. By the probability reduction
property (\ref{eq:Marked=000020state=000020property}) of the marked
state $a$, we have for some fixed $c<1$
\[
\mu_{S_{r}}(w)=\mu_{S_{r}}([w])\cdot\prod_{j\in J^{c}}\mu_{S_{r}}(w_{N_{j}})\leq\mu_{S_{r}}([w])\cdot c^{|J^{c}|}
\]
thus
\[
|J^{c}|\leq\frac{I(w|\hat{w})-I([w]|\hat{w})}{\log c^{-1}}.
\]
Assuming $X\in E_{r}$, we have that
\[
\#\{i\in S_{r}:[X_{S_{r}}]\not\Rrightarrow X_{i}\}\leq|J^{c}|\log^{2}|S_{r}|
\]
assuming further $X\in E_{r}\cap F_{r}\cap G_{r}$, we get from the
last two inequalities:
\[
\#\{i\in S_{r}:[X_{S_{r}}]\not\Rrightarrow\not X_{i}\}\leq\frac{2}{\log c^{-1}}\Lambda_{3,\frac{1}{2}}|S_{r}|
\]
which holds with probability $\mathbf{P}[E_{r}\cap F_{r}\cap G_{r}]=1-\Delta(L_{r})$,
thus shows the fulfillment of Property (\ref{eq:good=000020expected=000020number=000020of=000020determined=000020coordinates}).
In particular, it implies Property \ref{enu:For-,-the} of SOF. We
have thus shown that $\mathscr{A}$ is a good SOF for $\mathbf{A}$
that satisfies the additional property (\ref{eq:good=000020expected=000020number=000020of=000020determined=000020coordinates}),
and the proof is complete.$\hfill\square$

We are now ready to prove theorem \ref{thm:common=000020entry=000020version}:

\subsubsection{\textit{Proof of Theorem \ref{thm:common=000020entry=000020version}:}}

\noindent Given two marked processes $X$ and $Y$ with equal marked-state
distribution $\widehat{X}\overset{D}{=}\widehat{Y}$, use Corollary
\ref{cor:Marker=000020process} to find a nice marker factor of their
common distribution $\xi=\xi(\widehat{X})=\xi(\widehat{Y})$ (with
lengths sequence $(L_{r})_{r}$). Use Proposition \ref{prop:Existence=000020of=000020good=000020filler=000020systems}
to find good SOFs with respect to $\xi$ for $X$ and $Y$. This allows
us to find an isomorphism $X\overset{\phi}{\rightarrow}Y$ as in lemma
\ref{lem:good=000020systems=000020of=000020fillers=000020imply=000020good=000020isomorphism}. 

Given $n\in\mathbb{N}$, lemma \ref{lem:good=000020systems=000020of=000020fillers=000020imply=000020good=000020isomorphism}
provides some $r_{1},r_{2}$ with 
\begin{equation}
L_{r_{1}},L_{r_{2}}=\tilde{\Theta}(n)\label{eq:weak=000020theorem=0000200}
\end{equation}
 that satisfy
\begin{equation}
\mathbf{P}[[X_{S_{r_{1}}}]\Rrightarrow[\phi(X)_{S_{r_{2}}}]]=1-\Delta(n).\label{eq:weak=000020theorem=0000201}
\end{equation}

To estimate the probability of $\{X_{[-n,n]}\Rrightarrow\phi(X)_{0}\}$,
notice that 
\begin{align}
\{X_{[-n,n]}\Rrightarrow\phi(X)_{0}\}\supset & \{X_{[-n,n]}\Rrightarrow X_{[-L_{r_{1}},L_{r_{1}}]}\}\label{eq:-23}\\
 & \cap\{X_{[-L_{r_{1}},L_{r_{1}}]}\Rrightarrow X_{S_{r_{1}}}\}\nonumber \\
 & \cap\{[X_{S_{r_{1}}}]\Rrightarrow[\phi(X)_{S_{r_{2}}}]\}\nonumber \\
 & \cap\{[\phi(X)_{S_{r_{2}}}]\Rrightarrow\phi(X)_{0}\}\nonumber 
\end{align}
Thus it is enough to show that each of the four events on the right
is of probability $=1-\tilde{O}(n^{-\frac{1}{2}})$. By (\ref{eq:weak=000020theorem=0000200}),
a small change in $n$ gives that the first event on the right of
(\ref{eq:-23}),
\[
\{X_{[-n,n]}\Rrightarrow X_{[-L_{r_{1}},L_{r_{1}}]}\}
\]
 always holds. By the definition of lengths sequence (right after
Corollary \ref{cor:Marker=000020process}), the second event on the
right hand side of (\ref{eq:-23}),
\[
\{X_{[-L_{r_{1}},L_{r_{1}}]}\Rrightarrow X_{S_{r_{1}}}\}
\]
holds with probability $1-\Delta(n)$. The third event on the right
of (\ref{eq:-23}) has a good estimate in (\ref{eq:weak=000020theorem=0000201}),
hence it remains to show that the fourth event, $\{[\phi(X)_{S_{r_{2}}}]\Rrightarrow\phi(X)_{0}\}$,
has high enough probability. 

As our SOFs are provided by Proposition \ref{prop:Existence=000020of=000020good=000020filler=000020systems}
and thus satisfying (\ref{eq:good=000020expected=000020number=000020of=000020determined=000020coordinates}),
a mass-transport argument yields:
\begin{align*}
\mathbf{P}[[\phi(X)_{S_{r_{2}}}]\Rrightarrow\phi(X)_{0}] & =\mathbf{E}[\frac{1}{|S_{r_{2}}|}\sum_{i}\mathbf{1}\{[X_{S_{r_{2}}}]\Rrightarrow X_{i}\}]\\
 & \geq\mathbf{E}[1-\tilde{O}(|S_{r_{2}}|^{-\frac{1}{2}}]-\Delta(L_{r_{2}})\\
 & =1-\tilde{O}(L_{r_{2}}^{-\frac{1}{2}})\\
 & =1-\tilde{O}(n^{-\frac{1}{2}})
\end{align*}
where the last two equalities follow from $\xi$ being nice and (\ref{eq:weak=000020theorem=0000200}).
Overall we get
\[
\mathbf{P}[X_{[-n,n]}\Rrightarrow\phi(X)_{0}]\geq1-\Delta(n)-\Delta(n)-\tilde{O}(n^{-\frac{1}{2}})=1-\tilde{O}(n^{-\frac{1}{2}})
\]

\noindent And the proof is complete.$\hfill\square$

\section{\protect\label{sec:The-general-case}The general case}

\subsection{Chaining two Markov processes by a sequence of bi-marked processes}

To lift the results of section \ref{sec:Common=000020entry} (Theorem
\ref{thm:common=000020entry=000020version}) to the general case of
any two aperiodic Markov processes of equal entropy, we get into the
next lemma, Lemma \ref{lem:marked=000020process=000020lemma}. In
this lemma we endow any two aperiodic Markov processes $X,Y$ of equal
entropy with a finite sequence of processes $\{X(i)\}_{i=0}^{n}$
that satisfies the following: First, $X$ and $Y$ are isomorphic
to $X(0)$ and $X(n)$ (respectively) by a finite isomorphism. Second,
any two consecutive processes $X(i-1),X(i)$ are isomorphic by an
isomorphism $\phi^{i}$ as in Theorem \ref{thm:common=000020entry=000020version}.
Furthermore, the processes $\{X(i)\}_{i=0}^{n}$ have some particular
properties that enable us to construct the isomorphisms $\phi^{i}$
in a way that makes the decomposition $\phi=\phi^{n}\circ\cdots\circ\phi^{1}$
satisfying the claimed tail law of our main theorem. To this end,
we introduce a very particular kind of a process, which plays an important
role in this section. 
\begin{defn}
\textit{(bi-marked process)} Suppose $M$ is a stationary process
on three states $\{a,b,c\}$, such that $a$ and $b$ are renewal
states, and moreover, each of the distributions $\widehat{M}^{a},\widehat{M}^{b}$
has law identical to the distribution of some aperiodic Markov process
state. Thus $a,b$ both have exponential hitting times, and each of
the distributions $\widehat{M}^{a},\widehat{M}^{b}$ has a nice marker
factor (see Def. \ref{def:nice=000020marker} and the marker process
lemma, Lemma \ref{lem:Marker=000020process}). We call such a process
a \textit{bi-marked process} (with marked states $a$ and $b$).
\begin{defn}
\textit{(Independent splitting of a bi-marked process)} Suppose further
that $B(p)$ is an i.i.d. process independent of $M$, with marginal
law $p$. Let $Z$ be the process defined by
\begin{equation}
Z_{n}=\begin{cases}
M_{n} & M_{n}\neq c\\
B(p)_{n} & M_{n}=c
\end{cases}\label{eq:-50}
\end{equation}
 In that case we say that $Z$ is an \textit{independent splitting
of the bi-marked process $M$}.
\end{defn}

\end{defn}

Writing $Z\sim(C^{\mathbb{Z}},\nu)$, one has that both $\mathbf{C}_{a}:=(C^{\mathbb{Z}},\nu,a)$
and $\mathbf{C}_{b}:=(C^{\mathbb{Z}},\nu,b)$ are marked processes
(i.e. the states $a,b\in C$ are renewal and satisfying the probability
reduction property (\ref{eq:Marked=000020state=000020property})). 

The lemma is based on a collection of well known observations that
were used in the same manner to derive finitary isomorphisms between
processes (\cite{key-2,key-17}, see also \cite{key-18} and \cite[lemma 2.5]{key-24}).
Recall that the $k$-stringing $X^{(k)}$ of a process $X$ is the
process defined by $\left(X_{n}^{(k)}\right)_{n\in\mathbb{Z}}=\left(X_{[n,n+k-1]}\right)_{n\in\mathbb{Z}}$.
\begin{lem}
\label{lem:marked=000020process=000020lemma}For any two finite state
aperiodic Markov processes $X$ and $Y$ of the same entropy, there
is a finite sequence $\left(X(i)\sim(A(i)^{\mathbb{Z}},\mu(i))\right)_{i=0}^{n}$
of stationary processes with equal entropy, satisfying the following:

1. For some $k\geq1$, $X(0)=X^{(k)}$ and $X(n)=Y^{(k)}$ are the
$k$-stringings of $X$ and $Y$.

2. For each $1\leq i<n$, $X(i)$ is an independent splitting of a
bi-marked process $M(i)$ with marked states $a(i),b(i)$, so that
for all $1\leq i\leq n$, 
\[
\widehat{X(i-1)}^{b(i-1)}\overset{{\scriptscriptstyle dist.}}{=}\widehat{X(i)}^{a(i)}
\]
where $b(0)\in A(0)$ and $a(n)\in A(n)$ are renewal states satisfying
the probability reduction property (\ref{eq:Marked=000020state=000020property}).
\end{lem}

This lemma is the main result of \cite{key-17}, although the states
$b(0),a(n)$ provided there do not satisfy the probability reduction
property in general. However, a slight modification in the choice
of these states handles this problem. We indicate the modification
needed for that property to hold, while the full construction of the
states is omitted. 

In more detail, given an aperiodic Markov process $X$ on a state
space $A=\{1,...,|A|\}$, it is shown in \cite{key-17} that for any
large enough $k$, there are words $\alpha^{i}\in A^{k}$ ($0\leq i\leq|A|$)
which are allowable (i.e. $\mathbf{P}[X_{[1,k]}=\alpha^{i}]>0$ for
all $i$), so that for some $g\in A$,
\begin{align*}
\alpha_{1}^{i}=g, & \alpha_{k}^{i}=i\qquad(\forall i=1...|A|)
\end{align*}
and for all $i\in A$ and $1<j\leq k$, 
\[
\mathbf{P}[X_{[j,j+k-1]}=\alpha^{i}|X_{[1,k]}=\alpha^{0}]=0.
\]
The word $\alpha^{0}$ serves as the marked state of $X^{(k)}$ (choosing
some large $k$), which in our notations is denoted by the state $b(0)$
of the process $X(0)$. In order to make $\alpha^{0}$ satisfy the
probability reduction property, the idea is to choose it in a way
that for any allowable word $w\in A^{n}$ with $\hat{w}^{\alpha^{0}}=10\cdots01$
(here $w$ is viewed as a word of $X^{(k)}$) of large enough size
(say, $n>n_{0}$), by changing a bounded number of symbols in $w$
(which does not depend on $|w|$), another allowable word $u\in A^{n}$
is obtained, satisfying $\hat{u}^{\alpha^{0}}=10\cdots01$. The bound
on the number of symbols that have been changed forces that for some
fixed $c>0$, 
\[
\mathbf{P}[X_{[0,n-k]}^{(k)}=u|\hat{u}^{\alpha^{0}}=10\cdots01]\geq c\mathbf{P}[X_{[0,n-k]}^{(k)}=w|\hat{w}^{\alpha^{0}}=10\cdots01]
\]
 This implies
\[
\mathbf{P}[X_{[0,n-k]}^{(k)}=w|\hat{w}^{\alpha^{0}}=10\cdots01]\leq\frac{1}{1+c}
\]
thus the probability reduction property holds for large words of size
greater than $n_{0}$. As there are finitely many words of size at
most $n_{0}$, the probability reduction property trivially holds
on that range, and the modification is complete.

\subsection{\protect\label{subsec:Good-SOFs-for}Good SOFs for bi-marked process
splittings}

Given $Z\sim\left(C^{\mathbb{Z}},\mu\right)$ an independent splitting
of a bi-marked process $M$ with marked states $a,b$, endow $\mathbf{C}_{a}=\left(C^{\mathbb{Z}},\mu,a\right)$
and $\mathbf{C}_{b}=\left(C^{\mathbb{Z}},\mu,b\right)$ with nice
marker factors of their distributions, $\varsigma$ and $\vartheta$
respectively, with associated lengths sequences \textit{$\left(L_{\varsigma,r}\right)_{r},\left(L_{\vartheta,r}\right)_{r}$}
. We refer to $\varsigma$-skeletons with the capital letter $S$
and to $\vartheta$-skeletons with the capital letter $T$. Thus $S_{r}=S_{r}(0,Z)$
is the random $r$-skeleton for $\varsigma$ that covers $Z$ at 0,
while $T_{r}=T_{r}(0,Z)$ is the random $r$-skeleton for $\vartheta$
that covers $Z$ at $0$. 
\begin{lem}
\label{Nested=000020fillers=000020in=000020r=000020property}Under
the above circumstance, there exist good SOFs, $\mathscr{C}_{a}$
for $(\mathbf{C}_{a},\varsigma)$ and $\mathscr{C}_{b}$ for $(\mathbf{C}_{b},\vartheta)$,
that satisfy the following:

For any $0<\epsilon<\frac{1}{2}$, with probability $1-\Delta(L_{\varsigma,r})$,
for any $\vartheta$-skeleton $T$ appearing in $Z_{S_{r}}$ and of
size $|T|\geq|S_{r}|^{\epsilon}$, if $Z_{S_{r}}\Rrightarrow T,[Z_{T}]_{\mathbf{C}_{b}}$,
then $[Z_{S_{r}}]_{\mathbf{C}_{a}}\Rrightarrow T,[Z_{T}]_{\mathbf{C}_{b}}$
too. 

Similarly, for any $0<\epsilon<\frac{1}{2}$, with probability $1-\Delta(L_{\vartheta,r})$,
for any $\varsigma$-skeleton $S$ appearing in $Z_{T_{r}}$ and of
size $|S|\geq|T_{r}|^{\epsilon}$, if $Z_{T_{r}}\Rrightarrow S,[Z_{S}]_{\mathbf{C}_{a}}$,
then $[Z_{T_{r}}]_{\mathbf{C}_{b}}\Rrightarrow S,[Z_{S}]_{\mathbf{C}_{a}}$
too. 
\end{lem}

\subsubsection*{Proof outline of Lemma \ref{Nested=000020fillers=000020in=000020r=000020property}.}

As in the proof outline for Proposition \ref{prop:Existence=000020of=000020good=000020filler=000020systems},
we will restrict our discussion to the simpler case where $Z$ is
an i.i.d. process, and remark the places where omitting this assumption
must be carefully handled. We describe here what are the modifications
needed to apply to the construction there, in order to get good-SOFs
that satisfy the additional properties listed in Lemma \ref{Nested=000020fillers=000020in=000020r=000020property}.

Given a $\varsigma$-skeleton $S$ and a word $w\in C_{S}$, for the
class $[w]$ to determine (with probability $1-\Delta(|S|)$) any
(large enough) $\vartheta$-skeleton $T$ appearing in $Z_{S}$, we
will require $[w]$ to reveal the whole sample $\widehat{w}^{b}$,
at least with high probability. For that, besides of defining $\mathcal{J}(w)$
as we did in Lemma \ref{Nested=000020fillers=000020in=000020r=000020property},
we attach to $w$ another index set, $\text{\ensuremath{\mathcal{I}(w)}}\ensuremath{\subset}S$,
with which the class of $w$ is defined as 
\begin{equation}
w'\in[w]\iff\left(\begin{array}{c}
\widehat{\ensuremath{w}}_{\mathcal{I}(w)}^{b}=\widehat{w'}_{\mathcal{I}(w)}^{b}\\
\&\\
w_{\mathcal{J}(w)}=w'{}_{\mathcal{J}(w)}
\end{array}\right).\label{eq:-61}
\end{equation}
We begin with the requirement that $\mathcal{I}(w)$ will contain
$\cup\mathcal{I}(w^{i})$ and $\mathcal{J}(w)$ will contain $\cup\mathcal{J}(w^{i})$.
Fix some $\delta>\frac{1}{2}$, and enlarge the set $\mathcal{I}(w)$,
and then maybe enlarge $\mathcal{J}(w)$ too, until the following
will hold:
\begin{equation}
I([w]|\widehat{wa})\geq|S|h(Z|\widehat{Z})-\Lambda_{\delta,\frac{1}{2}}|S|.\label{eq:-12}
\end{equation}
Using an Azuma-Hoeffding argument, it is shown in Proposition \ref{prop:-is-a}
that 
\begin{equation}
\mathbf{P}[\mathcal{I}(w)=S]=1-\Delta(|S|).\label{eq:-62}
\end{equation}
In particular, any $\vartheta$-skeleton $T'$ that is determined
by $\widehat{w}^{b}$, is very likely to be determined by $[w]$.
Up to this point, all definitions are reflected with no change for
any $\vartheta$-skeleton $T$, attaching to each word $w\in C_{T}$
the index sets $\mathcal{J}(w),\mathcal{I}(w)\subset T$, where $\mathcal{I}(w)=T$
with probability $1-\Delta(T)$, and with the relation respective
to (\ref{eq:-61}) gets the form of 
\begin{equation}
w'\in[w]\iff\left(\begin{array}{c}
\widehat{\ensuremath{w}}_{\mathcal{I}(w)}^{a}=\widehat{w'}_{\mathcal{I}(w)}^{a}\\
\&\\
w_{\mathcal{J}(w)}=w'{}_{\mathcal{J}(w)}
\end{array}\right).\label{eq:-61-1}
\end{equation}

Turning back to $\varsigma$-skeleton fillers $w\in C_{S}$, if $\mathcal{I}(w)=S$,
there's a room to enlarge $\mathcal{J}(w)$ as long as (\ref{eq:-12})
is yet to hold. Given a large enough $T'$ that is determined by $w$,
we note that for the class $[w]_{S}$ to determine $[w_{T'}]_{T'}$,
one shouldn't enlarge the set $\mathcal{J}(w)$ just by some prefix
of $S$, as this may not contain all the information about the class
$[w_{T'}]_{T'}$. Instead, we attach to $w$ an index set $\mathcal{J}_{0}(w)$
of indices covered by the skeleton that are ``allowed'' to be added
into $\mathcal{J}(w)$. The subset $\mathcal{J}_{0}(w)\subset S$
doesn't really on the full name of the word $w$, but rather on its
$(a,b)$-name $\widehat{w}^{ab}$, where
\[
(\widehat{w}^{ab})_{n}=\begin{cases}
a & w_{n}=a\\
b & w_{n}=b\\
0 & Otherwise
\end{cases}
\]
(the value $\widehat{w}^{ab}$ is determined from $[w]$, assuming
$\mathcal{I}(w)=S$). Having $\mathcal{J}_{0}(w)$ in hand, $\mathcal{J}(w)$
is defined as the union of $\cup\mathcal{J}(w^{i})$ together with
the smallest prefix of $\mathcal{J}_{0}(w)$ that satisfy inequality
(\ref{eq:-12}). $\mathcal{J}_{0}(w)$ is constructed as an enlargement
of $\cup\mathcal{J}_{0}(w^{i})$, thus inductively one gets that $\mathcal{J}(w)\subset\mathcal{J}_{0}(w)$.
According to the choice of the $\delta$ in (\ref{eq:-12}), we build
the $\mathcal{J}_{0}(w)$ to be small enough so that with high probability,
the inequality of (\ref{eq:-12}) will not hold, regardless of what
the subset $\mathcal{J}(w)\subset\mathcal{J}_{0}(w)$ is (however,
inequality (\ref{eq:-12}) will take place for some higher delta,
so that Property \ref{enu:For-any-large} of good SOF (Definition
\ref{def:Good=000020SOF}) will still hold). That is, 
\begin{equation}
\mathbf{P}[\mathcal{J}(w)=\mathcal{J}_{0}(w)]=1-\Delta(|S|).\label{eq:-63}
\end{equation}
with high probability, $\mathcal{J}(w)$ will be equal to $\mathcal{J}_{0}(w)$. 

For the definition of $\mathcal{J}(w)$ for fillers $w\in C_{T}$
of a $\vartheta$-skeleton $T$ (assuming as before that $\mathcal{I}(w)=T$),
we use the fact that $\widehat{w}^{ab}$ is completely known (by the
assumption of $\mathcal{I}(w)=T$) to determine from $[w]_{T}$ all
$\varsigma$-skeletons that are determined by $w$. Letting $S^{1},...,S^{m}$
be all maximal $\varsigma$-sub-skeletons appear at and determined
by $w$, We define $\mathcal{J}_{0}(w)$ as the union 
\[
\mathcal{J}_{0}(w)=\bigcup_{i=1}^{m}\mathcal{J}_{0}(w_{S^{i}})
\]
and $\mathcal{J}(w)$ is then defined to be the smallest subset of
it satisfying (\ref{eq:-12}). While the validity of (\ref{eq:-12})
for $\mathscr{C}_{a}$ for some fixed $\delta$ is quite straight
forward from the definitions of $\mathcal{J}_{0}$, it takes more
effort to show that it holds for $\mathscr{C}_{b}$ too. It is not
immediately evident that the determined sub-skeletons $S^{1},...,S^{m'}$
cover a sufficient number of indices of $T$. To demonstrate this,
we leverage the unique properties of the marker processes discussed
in Section \ref{sec:Marker-processes}.$\hfill\blacksquare$

Below we give a detailed proof of the lemma. We begin with the construction
of the SOF \textit{$\mathscr{C}_{a}$} (Sec. \ref{subsec:Construction-of})
and prove it is a good SOF and some key properties of it (Sec. \ref{subsec:Properties-of}).
Then we'll move to the construction of the SOF \textit{$\mathscr{C}_{b}$}
(Sec. \ref{subsec:Construction-of-1}) and prove it is a good SOF
and show its additional properties (Sec. \ref{subsec:Properties-of-1}).
Finally, we'll show how all those properties imply the above lemma
(Sec. \ref{subsec:Proof-of-Lemma}).

\subsubsection{\protect\label{subsec:Construction-of}Construction of \textit{$\mathscr{C}_{a}$}}

Fix $\delta>\frac{1}{2}$. Let 
\begin{align*}
\theta & =\min_{(a',a'')}\mathbf{P}[Z_{0}=a'|Z_{-1}=a'']
\end{align*}
where $(a',a'')\in C\times C$ ranges over all tuples with $\mathbf{P}[Z_{0}=a'|Z_{-1}=a'']>0$.
Fix $t_{0}>0$ to be so large that Corollary \ref{cor:For-any-,}
will hold for $k=2\delta$ with this $t_{0}$, and so large that for
any $t\geq t_{0}$
\begin{equation}
\Lambda_{\delta,\frac{1}{2}}t-\log^{2}t\cdot\log\theta^{-1}\geq\sqrt{t}.\label{eq:-45}
\end{equation}
Given a skeleton $S=\left(\hat{z}_{[\alpha,\beta]}^{a},\varsigma(z)_{[\alpha,\beta]}\right)$
and a filler $w\in C_{S}$, let
\[
\mathcal{L}(w)=\{i\in S:w_{i}\neq a,b\}\subset S
\]
and define the subsets $\mathcal{J}_{0}(w)\subset\mathcal{L}(w)$,
$\mathcal{J}(w)\subset\mathcal{J}_{0}(w)$ and $\mathcal{I}(w)\subset S$
as follows: 

\subparagraph*{Trivial definitions for small skeletons. }

If $|\mathcal{L}(w)|<t_{0}$, then we put $\mathcal{I}(w)=\mathcal{J}(w)=\mathcal{J}_{0}(w)=\emptyset$. 

\subparagraph{Define $\mathcal{J}_{0}(w)$ in general.}

Suppose now that $|\mathcal{L}(w)|\geq t_{0}$, write $S=S^{1}\ast\cdots\ast S^{m}$
and accordingly $w=w^{1}*\cdots w^{m}$ for its $(\text{mr}(S)-1)$-decomposition,
and assume $\mathcal{I}(w^{i})$, $\mathcal{J}_{0}(w^{i})\subset\mathcal{L}(w^{i})$
and $\mathcal{J}(w^{i})\subset\mathcal{J}_{0}(w^{i})$ are all defined
(or alternatively, assuming $\text{mr}(S)=1$ and letting $m=0$).
To define $\mathcal{\mathcal{J}}_{0}(w)$, we first assume that for
all $1\leq i\leq m$ with $|S^{i}|\geq t_{0}$,
\[
|\mathcal{J}_{0}(w^{i})|=\lfloor|\mathcal{L}(w^{i})|-\Lambda_{2\delta,\frac{1}{2}}|\mathcal{L}(w^{i})|\rfloor.
\]
Corollary \ref{cor:For-any-,} (together with the choice of $t_{0}$
above) implies
\begin{equation}
|\bigcup_{i}\mathcal{J}_{0}(w^{i})|\leq\lfloor\sum_{|\mathcal{L}(w^{i})|\geq t_{0}}\left(|\mathcal{L}(w^{i})|-\Lambda_{2\delta,\frac{1}{2}}|\mathcal{L}(w^{i})|\right)\rfloor\leq\lfloor|\mathcal{L}(w)|-\Lambda_{2\delta,\frac{1}{2}}|\mathcal{L}(w)|\rfloor\label{eq:-52}
\end{equation}
thus we define
\[
\mathcal{J}_{0}(w)=\left([\alpha,\gamma_{\mathcal{J}_{0}}]\cap\mathcal{L}(w)\right)\cup\bigcup_{i}\mathcal{J}_{0}(w^{i})
\]
where $[\alpha,\gamma_{\mathcal{J}_{0}}]$ is the largest prefix of
$S$ that satisfies
\begin{equation}
|\mathcal{J}_{0}(w)|=\lfloor|\mathcal{L}(w)|-\Lambda_{2\delta,\frac{1}{2}}|\mathcal{L}(w)|\rfloor\label{eq:-49}
\end{equation}
Such a $\gamma_{\mathcal{J}_{0}}$ certainly exists by (\ref{eq:-52}),
and the validity of (\ref{eq:-49}) enables to define $\mathcal{J}_{0}(w)$
inductively for all fillers of all $\varsigma$-skeletons. 

\subparagraph{Define $\mathcal{I}(w),\mathcal{J}(w)$ in general.}

To define $\mathcal{I}(w),\mathcal{J}(w)$, write 
\[
\begin{array}{cc}
\mathcal{I}'(w)=\bigcup_{i=1}^{m}\mathcal{I}(w_{S^{i}}), & \mathcal{J}'(w)=\bigcup_{i=1}^{m}\mathcal{J}(w^{i})\end{array}
\]
(which are empty when $\text{mr}(S)=1$) and consider two cases:

\textbf{Case 1 -} $S$ contains an interval of size at least $\log^{2}|S|$
with no $a$ (the ``rare'' case, see (\ref{eq:-47})): Put
\[
\begin{array}{lr}
\mathcal{I}(w)=\mathcal{I}'(w), & \mathcal{J}(w)=\mathcal{J}'(w).\end{array}
\]

\textbf{Case 2 -} $S$ does not contain an $a$-free interval of size
$>\log^{2}|S|$: Define $\mathcal{I}(w)$ to be
\[
\mathcal{I}(w)=[\alpha,\gamma_{a}]\cup\mathcal{I}'(w)
\]
where $[\alpha,\gamma_{a}]$ is the smallest prefix of $S$ with $z_{\gamma_{a}+1}=a$
for which
\begin{equation}
I(w_{\mathcal{J}'(w)},\widehat{w}^{b}{}_{\mathcal{I}(w)}|\widehat{wa}^{a})\geq|S|h(Z|\widehat{Z}^{a})-\Lambda_{\delta,\frac{1}{2}}|S|\label{eq:-19}
\end{equation}
(where $w_{\mathcal{J}'(w)},\hat{w}^{b}{}_{\mathcal{I}(w)},\widehat{wa}^{a}$
are considered as values of the random variables $X_{\mathcal{J}'(w)},\widehat{X}^{b}{}_{\mathcal{I}(w)},\widehat{X}_{[\alpha,\beta+1]}^{a}$
respectively). If there is no such a prefix, put $\mathcal{I}(w)=S$.
Now define $\mathcal{J}(w)$ to be
\[
\mathcal{J}(w)=\left([\alpha,\gamma_{\mathcal{J}}]\cap\mathcal{J}_{0}(w)\right)\cup\mathcal{J}'(w)
\]
where $[\alpha,\gamma_{\mathcal{J}}]$ is the smallest prefix of $S$
for which
\begin{equation}
I(w_{\mathcal{J}(w)},\hat{w}^{b}{}_{\mathcal{I}(w)}|\widehat{wa}^{a})\geq|S|h(Z|\hat{Z}^{a})-\Lambda_{\delta,\frac{1}{2}}|S|,\label{eq:-25}
\end{equation}
and if there is no such a prefix, put $\mathcal{J}(w)=\mathcal{J}_{0}(w)$.

Finally, define for any skeleton $S$ the equivalence relation $\sim$
on $C_{S}$ by
\[
w\sim w'\iff\left(\begin{array}{c}
\widehat{\ensuremath{w}}_{\mathcal{I}(w)}^{b}=\widehat{w'}_{\mathcal{I}(w)}^{b}\\
\&\\
w_{\mathcal{J}(w)}=w'{}_{\mathcal{J}(w)}
\end{array}\right)
\]
and let $\mathscr{C}_{a}=\{(\tilde{C}_{S},\mu_{S})_{S}:S\text{ a \ensuremath{\varsigma}-skeleton}\}$.

\subsubsection{\protect\label{subsec:Properties-of}Properties of $\mathscr{C}_{a}$}
\begin{prop}
\label{prop:-is-a}$\mathscr{C}_{a}$ is a good SOF that in addition
satisfies:
\begin{equation}
\mathbf{P}[\mathcal{I}(Z_{S_{r}})=S_{r}\wedge\mathcal{J}(Z_{S_{r}})=\mathcal{J}_{0}(S_{r})]=1-\Delta(L_{\varsigma,r}).\label{eq:-22}
\end{equation}
\end{prop}

\begin{proof}
From the definition of $\mathscr{C}_{a}$ it is clear that it satisfies
the first two properties of SOF (Definition \ref{def:SOF}). That
$\mathscr{C}_{a}$ satisfies Property \ref{enu:For-any-large} of
good SOF (Definition \ref{def:Good=000020SOF}) with $k=0$, i.e.
that for any skeleton $S$ of large enough size $|S|$, any class
$[w]\in\tilde{C}_{S}$ satisfies
\begin{equation}
I_{\mu_{S}}([w])\leq|S|h(Z|\widehat{Z}^{a})-\sqrt{|S|}\label{eq:-43}
\end{equation}
is proved in a similar manner to the proof of Proposition \ref{prop:Existence=000020of=000020good=000020filler=000020systems}:
the only case with a slight difference is Case 2, where $|S|>t_{0}$
and $S$ does not contain an $a$-free interval of size $>\log^{2}|S|$:
Assume (\ref{eq:-25}) holds (otherwise (\ref{eq:-43}) follows directly).
We either have 
\begin{equation}
\begin{array}{ccc}
\mathcal{I}(w)=[\alpha,\gamma_{a}]\cup\mathcal{I}'(w) & \text{and} & \mathcal{J}(w)=\mathcal{J}'(w)\end{array}\label{eq:-44}
\end{equation}
 or $\mathcal{J}(w)\varsupsetneq\mathcal{J}'(w)$, where
\begin{equation}
\begin{array}{ccc}
\mathcal{I}(w)=S & \text{and} & \mathcal{J}(w)=\left([\alpha,\gamma_{\mathcal{J}}]\cap\mathcal{J}_{0}(w)\right)\cup\mathcal{J}'(w)\end{array}\label{eq:-46}
\end{equation}
 If (\ref{eq:-44}) holds, we can assume that $\mathcal{I}(w)\supsetneq\mathcal{I}'(w)$
(or else the same argument as in Case 1 will yield (\ref{eq:-43})).
Thus
\begin{align*}
I([w]|\widehat{wa}^{a}) & \leq I(w_{\mathcal{J}(w)},\widehat{w}^{b}{}_{[\alpha,\tau_{a}^{<\gamma_{a}}]\cup\mathcal{I}'(w)}|\widehat{wa}^{a})+I[\widehat{w}^{b}{}_{[\tau_{a}^{<\gamma_{a}},\gamma_{a}]}|\widehat{wa}^{a})\\
 & \leq|S|h(Z|\widehat{Z}^{a})-\Lambda_{\delta,\frac{1}{2}}|S|+(\gamma_{a}-\tau_{a}^{<\gamma_{a}})\log\theta^{-1}
\end{align*}
(recall that $\tau_{a}^{<\gamma}$ is the last occurrence of $a$
before $\gamma$). Similarly, if (\ref{eq:-46}) holds instead of
(\ref{eq:-44}), then we have
\begin{align*}
I([w]|\widehat{wa}^{a}) & \leq I(w_{\left([\alpha,\gamma-1]\cap\mathcal{J}_{0}(w)\right)\cup\mathcal{J}'(w)},\widehat{w}^{b}{}_{\mathcal{I}(w)}|\widehat{wa}^{a})+I(w{}_{\gamma}|\widehat{w}_{\gamma}^{a,b})\\
 & \leq|S|h(Z|\widehat{Z}^{a})-\Lambda_{\delta,\frac{1}{2}}|S|+\log\theta^{-1}
\end{align*}
 Either way, since $(\gamma-\tau_{a}^{<\gamma})$ is assumed to be
at most $\log^{2}|S|$ and by the choice of $t_{0}$ (\ref{eq:-45}),
\begin{align*}
I([w]|\widehat{wa}^{a}) & \leq|S|h(Z|\widehat{Z}^{a})-\sqrt{|S|}
\end{align*}
and (\ref{eq:-43}) holds.

We now turn to prove the additional property (\ref{eq:-22}), and
then use it to show that $\mathscr{C}_{a}$ satisfies the remained
properties of SOF and good SOF. Notice the size difference between
the sets of indices that were left outside of $\mathcal{J}_{0}(w)$
and outside of $\mathcal{J}(w)$: While the one of $\mathcal{J}_{0}(w)$
was quantified in (\ref{eq:-49}) by terms of $2\delta$, the set
of indices that were left outside of $\mathcal{I}(w)$ in (\ref{eq:-19})
and $\mathcal{J}(w)$ in (\ref{eq:-25}) was quantified by terms of
$\delta$. This difference will cause, with high probability, the
set $\mathcal{J}(w)$ to contain the entire set $\mathcal{J}_{0}(w)$
before achieving the point where (\ref{eq:-25}) is satisfied.

We have:
\begin{align}
I(Z_{\mathcal{J}_{0}(Z_{S_{r}})},\widehat{Z}_{S_{r}}^{b}|\widehat{Z}^{a}) & =I(\widehat{Z}_{S_{r}}^{b}|\widehat{Z}^{a})+I(Z_{\mathcal{J}_{0}(Z_{S_{r}})}|\widehat{Z}^{a,b})\nonumber \\
 & =I(M_{S_{r}}|\widehat{Z}^{a})+I(Z_{\mathcal{J}_{0}(Z_{S_{r}})}|M)\label{eq:-57}
\end{align}
As $Z$ has $a$ and $b$ as renewal states of exponential hitting
times, one can use the results of Subsection \ref{subsec:Azuma-Hoeffding-inequality-for}
to estimate the last expression: The first summand above, $I(M_{S_{r}}|\widehat{Z}^{a})$,
can be estimated by Corollary \ref{cor:Azuma=000020information} of
Azuma's inequality, that with probability $1-\Delta(L_{\varsigma,r})$,
\begin{equation}
I(M_{S_{r}}|\widehat{Z}^{a})\in|S_{r}|h(M|\widehat{Z}^{a})\pm\Lambda_{\delta,\frac{1}{2}}|S_{r}|.\label{eq:-55}
\end{equation}
The second summand of (\ref{eq:-57}), $I(Z_{\mathcal{J}_{0}(Z_{S_{r}})}|M)$,
can be estimated using the fact that $Z$ is an independent splitting
of $M$ (\ref{eq:-50}): We denote by $p$ the marginal law of $B(p)$,
which is the marginal law of $Z$ given it is not $a$ or $b$. The
definition of $\mathcal{J}_{0}(Z_{S_{r}})$ depends only on $\hat{Z}_{S_{r}}^{a,b}=M_{S_{r}}$
and thus independent of $B(p)$, hence
\begin{align}
I(Z_{\mathcal{J}_{0}(Z_{S_{r}})}|\hat{Z}^{a,b}) & \overset{{\scriptscriptstyle dist.}}{=}I(B(p)_{\mathcal{J}_{0}(Z_{S_{r}})})\overset{{\scriptscriptstyle dist.}}{=}I(B(p)_{[1,|\mathcal{J}_{0}(Z_{S_{r}})|]}).\label{eq:-48}
\end{align}
Applying Azuma's inequality (Theorem \ref{thm:Azuma-Hoeffding-inequality})
to the random variable $X:=\sum_{i=1}^{n}\mathbf{1}\{M_{i}=c\}$,
one can conclude similarly to Corollary \ref{cor:Azuma=000020information},
that with probability $1-\Delta(L_{\varsigma,r})$
\[
|\mathcal{L}(Z_{S_{r}})|\in\mathbf{P}[M_{0}=c]|S_{r}|\pm\Lambda_{\delta,\frac{1}{2}}|S_{r}|
\]
and so, by (\ref{eq:-49}) we get with probability $1-\Delta(L_{\varsigma,r})$
\begin{equation}
\mathcal{J}_{0}(Z_{S_{r}})\in\mathbf{P}[M_{0}=c]|S_{r}|-\Lambda_{2\delta,\frac{1}{2}}(\mathbf{P}[M_{0}=c]|S_{r}|)\pm O(\Lambda_{\delta,\frac{1}{2}}|S_{r}|).\label{eq:-51}
\end{equation}

This, together with (\ref{eq:-48}) gives that with probability $1-\Delta(L_{\varsigma,r})$,
\begin{align}
I(Z_{\mathcal{J}_{0}(Z_{S_{r}})}|M) & \in H(p)\left(\mathbf{P}[M_{0}=c]|S_{r}|-\Lambda_{2\delta,\frac{1}{2}}(\mathbf{P}[M_{0}=c]|S_{r}|)\right)\pm O(\Lambda_{\delta,\frac{1}{2}}|S_{r}|)\label{eq:-54}
\end{align}

Notice that
\begin{align*}
h(Z|M) & =\lim_{n}\frac{1}{n}\mathbf{E}[\mathbf{E}[I(\{X_{i}:M_{i}=c\}|M_{[1,n]})|M_{[1,n]}]]\\
 & =\lim_{n}\frac{1}{n}\mathbf{E}[H(p)\sum_{i=1}^{n}\mathbf{1}\{M_{i}=c\}]\\
 & =H(p)\mathbf{P}[M_{0}=c]
\end{align*}
thus
\begin{align}
h(Z|\widehat{Z}^{a}) & =h(Z|M)+h(M|\widehat{Z}^{a})=H(p)\mathbf{P}[M_{0}=c]+h(M|\widehat{Z}^{a})\label{eq:-53}
\end{align}
and if both (\ref{eq:-55}) and (\ref{eq:-54}) hold, then summing
them and assuming $|S_{r}|$ is large enough, one has
\begin{align}
I(Z_{\mathcal{J}_{0}(Z_{S_{r}})},\hat{Z}_{S_{r}}^{b}|\hat{Z}^{a}) & \in h(Z|\hat{Z}^{a})|S_{r}|-H(p)\Lambda_{2\delta,\frac{1}{2}}(\mathbf{P}[M_{0}=c]|S_{r}|)\pm O(\Lambda_{\delta,\frac{1}{2}}|S_{r}|)\label{eq:-56}\\
 & <h(Z|\hat{Z}^{a})|S_{r}|-\Lambda_{\delta,\frac{1}{2}}|S_{r}|.\nonumber 
\end{align}

On the event that the above holds and $S_{r}$ does not contain an
$a$-free interval of size $>\log^{2}|S|$, by the construction of
$\mathscr{C}_{a}$ (see Case 2) we get that $\mathcal{I}(Z_{S_{r}})=S$
and $\mathcal{J}(Z_{S_{r}})=\mathcal{J}_{0}(Z_{S_{r}})$. As these
events hold with probability $1-\Delta(L_{\varsigma,r})$, their intersection
is of the same order, thus the additional property (\ref{eq:-22})
is proved.

The validity of (\ref{eq:-22}) together with the lower bound appears
in (\ref{eq:-56}) gives that with probability $1-\Delta(L_{\varsigma,r})$,
\[
I([Z_{S_{r}}]|\hat{Z}^{a})=I(Z_{\mathcal{J}_{0}(Z_{S_{r}})},\hat{Z}_{S_{r}}^{b}|\hat{Z}^{a})\geq h(Z|\hat{Z}^{a})|S_{r}|-c\Lambda_{2\delta,\frac{1}{2}}|S_{r}|
\]
for some fixed $c>0$. Thus Property \ref{enu:Good=000020SOF=0000202}
of good SOF holds for $\mathscr{C}_{a}$.

Finally, another consequence of (\ref{eq:-22}) is that 
\[
\mathbf{P}[[Z_{S_{r}}]\Rrightarrow\{Z_{i}:i\in S_{r}\backslash(\mathcal{L}(Z_{S_{r}})\backslash\mathcal{J}_{0}(Z_{S_{r}}))]=1-\Delta(L_{\varsigma,r})
\]
 From the construction of $\mathcal{J}_{0}(Z_{S_{r}})\subset\mathcal{L}(Z_{S_{r}})$
(\ref{eq:-49}), 
\begin{align*}
\frac{\left|S_{r}\backslash(\mathcal{L}(Z_{S_{r}})\backslash\mathcal{J}_{0}(Z_{S_{r}}))\right|}{|S_{r}|} & =\frac{|S_{r}|-\left(|\mathcal{L}(Z_{S_{r}})|-|\mathcal{J}_{0}(Z_{S_{r}})|\right)}{|S_{r}|}\\
 & =1-\frac{\Lambda_{2\delta,\frac{1}{2}}|\mathcal{L}(w_{S})|+O(1)}{|S_{r}|}
\end{align*}
which converges a.s. to 1. Hence Property \ref{enu:For-,-the} of
SOFs (Definition \ref{def:SOF}) follows, and the proof is complete.
\end{proof}
We turn to define $\mathscr{C}_{b}$: 

\subsubsection{\protect\label{subsec:Construction-of-1}Construction of \textit{$\mathscr{C}_{b}$}}

Given a $\vartheta$-skeleton $T$, Let $T^{1},...,T^{m}$ be all
proper $\vartheta$-sub-skeletons of $T$, and for a given $w\in C_{T}$,
let $S^{1},...,S^{m'}$ be all maximal $\varsigma$-sub-skeletons
appear at and determined by $w$ (maximality w.r.t. inclusion). Suppose
that index sets $\mathcal{I}(\cdot),\mathcal{J}(\cdot)$ are defined
for all fillers $\left(w_{S^{i}}\right)_{i=1}^{m'}$ and $\left(w_{T^{j}}\right)_{j=1}^{m}$.
We now write $\mathcal{J}_{0}(w)=\bigcup_{i=1}^{m'}\mathcal{J}_{0}(w_{S^{i}})$
$\mathcal{I}'(w)=\bigcup_{i=1}^{m}\mathcal{I}(w_{T^{i}})$ and $\mathcal{J}'(w)=\bigcup_{i=1}^{m}\mathcal{J}(w_{T^{i}})$.
The definition of $\mathcal{I}$ and $\mathcal{J}$ in Case 0 - where
$|T|\leq t_{0}$ - and in Case 1 - where $|T|>t_{0}$ but there is
a $b$-free interval in $T$ of size $>\log^{2}|T|$ - is just the
same as we defined it on fillers of $\varsigma$-skeletons $S$.

In Case 2, where $|T|>t_{0}$ and $T$ does not contain a $b$-free
interval of size $>\log^{2}|T|$, the definition of $\mathcal{I}(w)$
is just the same as for $\mathscr{C}_{a}$, switching the roles of
$a$ and $b$: we define 
\[
\mathcal{I}(w)=[\alpha,\gamma_{b}]\cup\mathcal{I}'(w)
\]
where $[\alpha,\gamma_{b}]$ is the smallest prefix of $T$ with $z_{\gamma_{b}+1}=b$
for which
\[
I(w_{\mathcal{J}'(w)},\hat{w}^{a}{}_{\mathcal{I}(w)}|\widehat{wb}^{b})\geq|T|h(Z|\hat{Z}^{b})-\Lambda_{\delta,\frac{1}{2}}|T|
\]
and if there is no such a prefix, put $\mathcal{I}(w)=|T|$.

The definition of $\mathcal{J}(w)$ is also similar to that of $\mathscr{C}_{a}$:
we let
\[
\mathcal{J}(w)=\left([\alpha,\gamma]\cap\mathcal{J}_{0}(w)\right)\cup\mathcal{J}'(w)
\]
where $[\alpha,\gamma]$ is the smallest prefix of $S$ for which
\[
I(w_{\mathcal{J}(w)},\hat{w}_{\mathcal{I}(w)}^{a}|\widehat{wb}^{b})\geq|T|h(Z|\hat{Z}^{b})-\Lambda_{\delta,\frac{1}{2}}|T|
\]
and if there is no such a prefix, set $\mathcal{J}(w)=\mathcal{J}_{0}(w)$.
Notice however that $\mathcal{J}_{0}(w)$ has different definition
than the one defined for $\varsigma$-skeleton fillers.

Defining for any skeleton $T$ the relation $\sim$ on $C_{T}$ by
\[
w\sim w'\iff\begin{array}{c}
\widehat{w}_{\mathcal{I}(w)}^{a}=\widehat{w'}_{\mathcal{I}(w)}^{a}\\
\&\\
w_{\mathcal{J}(w)}=w'{}_{\mathcal{J}(w)}
\end{array}
\]
gives an equivalence relation. Let $\mathscr{C}_{b}=\{(\tilde{C}_{T},\mu_{T})_{T}:S\text{ a \ensuremath{\vartheta}-skeleton}\}$.

\subsubsection{\protect\label{subsec:Properties-of-1}Properties of $\mathscr{C}_{b}$}

Before proving that this construction satisfies the claimed properties
of Lemma \ref{Nested=000020fillers=000020in=000020r=000020property},
we make some preparations. For a skeleton $S$ and a fixed $\epsilon>0$,
let $\mathring{S^{\epsilon}}$ denote the set of indices in $S$ which
are not in its $|S|^{\epsilon}$-boundary:
\[
\mathring{S^{\epsilon}}=\{i\in S:d(i,S^{c})>|S|^{\epsilon}\}.
\]
Considering a point $z\in C^{\mathbb{Z}}$ with a $\vartheta$-skeleton
$T$ appearing in it, we write $S(i,T,z)$ for the largest $\varsigma$-skeleton
that covers $z$ at $i$, appears at and determined by $z_{T}$:
\begin{equation}
S(i,T,z)=\arg\max\{|S|:i\in S\subset T\,\&\,z_{T}\Rrightarrow S\}.\label{eq:-20}
\end{equation}
Where there is no ambiguity about $T$ and $z$, we abbreviate $S(i)=S(i,T,z)$.
\begin{prop}
\label{prop:For-any-,}For some $c>0$ and any $\epsilon>0$,
\[
\mathbf{P}\left[\bigcap_{i\in\mathring{T_{r}}^{\epsilon}}\{\Lambda_{-c}d(i,T_{r}^{c})\leq|S(i,T_{r},Z)|\}\right]=1-\Delta(L_{\vartheta,r}).
\]
\end{prop}

\begin{proof}
The marker process lemma (Lemma \ref{lem:Marker=000020process}) implies
that for some $c>0$, for any choice of $i=i_{r}\in[-L_{\vartheta,r},L_{\vartheta,r}]$,
any $\epsilon>0$ and any choice of $d=d_{r}\in($ $L_{\vartheta,r}^{\epsilon},L_{\vartheta,r}]$,
there is some sequence of $r(d_{r})$ for which 
\begin{align*}
\mathbf{P}[\left(Z_{[i-d,i+d]}\Rrightarrow S_{r(d)}(i,Z)\right)\&\left(\Lambda_{-c}d\leq|S_{r(d)}(i,Z)|\leq d\right)] & =1-\Delta(d)\\
 & =1-\Delta(L_{\vartheta,r}^{\epsilon})\\
 & =1-\Delta(L_{\vartheta,r}).
\end{align*}
Intersecting the last event over all $i$'s and $d$'s gives an event,
call it $D_{r}$, of probability 
\begin{align}
\mathbf{P}[D_{r}]=1-2L_{\vartheta,r}^{2}\Delta(L_{\vartheta,r})=1-\Delta(L_{\vartheta,r}).\label{eq:-26}
\end{align}
Using the nice property of $\vartheta$ (see Def. \ref{def:nice=000020marker}),
the event $E_{r}:=\{|T_{r}|\leq L_{\vartheta,r}\}$ is of probability
$1-\Delta(L_{\vartheta,r})$, which yields $\mathbf{P}[D_{r}\cap E_{r}]=1-\Delta(L_{\vartheta,r})$.
On $D_{r}\cap E_{r}$, for all $i\in\mathring{T_{r}^{\epsilon}}$
we have both $i\in[-L_{\vartheta,r},L_{\vartheta,r}]$ and $L_{\vartheta,r}^{\epsilon}<d(i,T_{r}^{c})$
so that
\[
\Lambda_{-c}d(i,T_{r}^{c})\leq|S_{r(d(i,T_{r}^{c}))}(i,Z)|\leq d(i,T_{r}^{c}).
\]

The right hand side inequality implies $S_{r(d(i,T_{r}^{c}))}(i,Z)\subset T_{r}$,
so that $S_{r(d(i,T_{r}^{c}))}\subset S(i)$. Taking into account
the left hand side inequality gives now
\[
\Lambda_{-c}d(i,T_{r}^{c})\leq S(i)|
\]
which ends the proof.
\end{proof}
\begin{prop}
\label{prop:Let--be}Let $Q$ be a shift-invariant collection of finite
words with symbols in $C$. Suppose that $\mathbf{P}[Z_{S_{r}}\in Q]=1-\Delta(L_{\varsigma,r})$.
Then for any $\epsilon>0$,
\[
\mathbf{P}\left[\bigcap_{i\in\mathring{T_{r}^{\epsilon}}}\{Z_{S(i,T_{r},Z)}\in Q\}\right]=1-\Delta(L_{\vartheta,r}).
\]
\end{prop}

\begin{proof}
We have
\[
\bigcap_{i\in\mathring{T_{r}^{\epsilon}}}\{Z_{S(i)}\in Q\}\supset\{\Lambda_{-c}L_{\varsigma,r}\leq|T_{r}|\leq L_{\varsigma,r}\}\cap\bigcap_{{\scriptscriptstyle \begin{array}{c}
r'\in\mathbb{N},i\in\mathbb{Z}\\
\leq|T_{r'}(i,Z)|\leq d(i,S_{r}^{c})
\end{array}}}\{Z_{T_{r'}(i,Z)}\in Q\}.
\]
Thus 
\begin{align*}
\mathbf{P}\left[\bigcap_{i\in\mathring{S_{r}^{\epsilon}}}\{Z_{T(i)}\in Q\}\right]\geq & \mathbf{P}\left[\Lambda_{-c}L_{\varsigma,r}\leq|S_{r}|\leq L_{\varsigma,r}\right]\\
 & -\left|\{T_{r'}(i,Z)\subset[-L_{\varsigma,r},L_{\varsigma,r}]\}\right|\cdot\min_{T_{r'}(i,Z)}\mathbf{P}[Z_{T_{r'}(i,Z)}\notin Q]\\
\geq & 1-\Delta(L_{\varsigma,r})-L_{\varsigma,r}^{2}\Delta(L_{\varsigma,r}^{\epsilon})\\
= & 1-\Delta(L_{\varsigma,r})
\end{align*}
which is just as claimed.
\end{proof}
Turning back to the construction of $\mathscr{C}_{b}$, just as we
had $|\mathcal{J}_{0}(Z_{S_{r}})|=\mathbf{P}[M_{0}=c]|S_{r}|-\tilde{O}(|S_{r}|)$
with high probability (\ref{eq:-51}), the size of $\mathcal{J}_{0}(Z_{T_{r}})$
has the same estimate with the same probability order:
\begin{prop}
For some $c,c_{1},c_{2}>0$, the following holds with probability
$1-\Delta(L_{\vartheta,r})$:
\[
\mathbf{P}[M_{0}=c]|T_{r}|-c_{1}\Lambda_{2\delta+c,\frac{1}{2}}|T_{r}|\leq|\mathcal{J}_{0}(Z_{T_{r}})|\leq\mathbf{P}[M_{0}=c]|T_{r}|-c_{2}\Lambda_{2\delta,\frac{1}{2}}|T_{r}|
\]
\end{prop}

\begin{proof}
Assuming $|T_{r}|$ is large enough, the definition of $\mathcal{J}_{0}$
(\ref{eq:-49}) together with Corollary \ref{cor:For-any-,} bounds
$|\mathcal{J}_{0}(Z_{T_{r}})|$ from above by
\begin{align*}
|\mathcal{J}_{0}(Z_{T_{r}})| & =\sum_{i=1}^{m'}|\mathcal{J}_{0}(Z_{S^{i}})|\\
 & =\sum_{i=1}^{m'}\lfloor|\mathcal{L}(Z_{S^{i}})|-\Lambda_{2\delta,\frac{1}{2}}|\mathcal{L}(Z_{S^{i}})|\rfloor\\
 & \leq|\mathcal{L}(Z_{T_{r}})|-\Lambda_{2\delta,\frac{1}{2}}|\mathcal{L}(Z_{T_{r}})|.
\end{align*}
On the other hand, we have
\begin{align*}
|\mathcal{J}_{0}(Z_{T_{r}})| & =\left|\bigcup_{i\in T_{r}}\mathcal{J}_{0}(Z_{S(i)})\right|\\
 & =\sum_{\begin{array}{c}
i\in T_{r}\\
S(i)\neq\emptyset
\end{array}}\frac{|\mathcal{J}_{0}(Z_{S(i)})|}{|S(i)|}\\
 & =|\bigcup_{\begin{array}{c}
i\in T_{r}\\
S(i)\neq\emptyset
\end{array}}\mathcal{L}(Z_{S(i)})|-\sum_{{\scriptstyle \begin{array}{c}
i\in T_{r}\\
S(i)\neq\emptyset
\end{array}}}\left(\lceil\Lambda_{2\delta,\frac{1}{2}}|\mathcal{L}(Z_{S(i)})|\rceil/|S(i)|\right)\wedge1.
\end{align*}
 By Proposition \ref{prop:For-any-,}, for some fixed $c>0$, with
probability $1-\Delta(L_{\vartheta,r})$ we have
\begin{align*}
\sum_{i\in\mathring{T}_{r}^{\epsilon}}\lceil\Lambda_{2\delta,\frac{1}{2}}|\mathcal{L}(Z_{S(i)})|\rceil/|S(i)| & \leq\sum_{i\in\mathring{T}_{r}^{\epsilon}}\lceil\Lambda_{2\delta,\frac{1}{2}}|S(i)|\rceil/|S(i)|\\
 & \leq\sum_{i\in\mathring{T}_{r}^{\epsilon}}\Lambda_{2\delta+c,-\frac{1}{2}}d(i,T_{r}^{c})\\
 & =2\sum_{n=|T_{r}|^{\epsilon}}^{|T_{r}|/2}\Lambda_{2\delta+c,-\frac{1}{2}}n\\
 & =\Theta\left(\int_{|T_{r}|^{\epsilon}}^{|T_{r}|/2}\Lambda_{2\delta+c,-\frac{1}{2}}xdx\right)\\
 & =\Theta(\Lambda_{2\delta+c,\frac{1}{2}}|T_{r}|)
\end{align*}
By the same proposition, together with Azuma's inequality, we have
for some fixed $\delta_{0}>\delta>\frac{1}{2}$, that with probability
$1-\Delta(L_{\vartheta,r})$:
\[
|\bigcup_{\begin{array}{c}
i\in T_{r}\\
S(i)\neq\emptyset
\end{array}}\mathcal{L}(Z_{S(i)})|\geq|\mathcal{L}(Z_{\mathring{T}_{r}^{\epsilon}})|\geq\mathbf{P}[M_{0}=c]|\mathring{T}_{r}^{\epsilon}|-\Lambda_{\delta,\frac{1}{2}}|\mathring{T}_{r}^{\epsilon}|
\]
Thus
\begin{align*}
|\mathcal{J}_{0}(Z_{T_{r}})| & =\mathbf{P}[M_{0}=c]|T_{r}|-O(\Lambda_{\delta_{0}+c,\frac{1}{2}}|T_{r}|)-O(\Lambda_{\delta,\frac{1}{2}}|T_{r}|)-O(|T_{r}|^{\epsilon})\\
 & =\mathbf{P}[M_{0}=c]|T_{r}|-O(\Lambda_{\delta_{0}+c,\frac{1}{2}}|T_{r}|).
\end{align*}
\end{proof}
By the last estimate of $|\mathcal{J}_{0}(Z_{T_{r}})|$, following
the same lines as in the proof of Proposition \ref{prop:-is-a}, we
get:
\begin{prop}
\label{prop:-is-a-1}$\mathscr{C}_{b}$ is a good SOF that in addition
satisfies:
\begin{equation}
\mathbf{P}[\mathcal{I}(Z_{T_{r}})=T_{r},\mathcal{J}(Z_{T_{r}})=\mathcal{J}_{0}(T_{r})]=1-\Delta(L_{\vartheta,r}).\label{eq:-22-1}
\end{equation}
\end{prop}

\subsubsection{\protect\label{subsec:Proof-of-Lemma}Proof of Lemma \textit{\ref{Nested=000020fillers=000020in=000020r=000020property}}}

 We have shown that these systems are good SOFs in Propositions \ref{prop:-is-a}
and \ref{prop:-is-a-1}, thus it remains to show they satisfy the
additional property there. We will prove its first part, i.e. that
with probability $1-\Delta(L_{\varsigma,r})$, for any $\vartheta$-skeleton
$T$ appearing in $Z_{S_{r}}$ and of size $|T|\geq|S_{r}|^{\epsilon}$,
if $Z_{S_{r}}\Rrightarrow T,[Z_{T}]_{\mathbf{C}_{b}}$, then $[Z_{S_{r}}]_{\mathbf{C}_{a}}\Rrightarrow T,[Z_{T}]_{\mathbf{C}_{b}}$
too. the second part is proved following the same lines.

Let $\mathcal{T}$ be the set of all $\vartheta$-skeletons $T\subset S_{r}$
of size $|T|\geq|S_{r}|^{\epsilon}$ such that $Z_{S_{r}}\Rrightarrow T$.
Clearly, any $T\in\mathcal{T}$ is included at some maximal $T'\in\mathcal{T}$,
and $[Z_{T'}]\Rrightarrow[Z_{T}]$. Thus it is enough to prove the
claim for all maximal skeletons 
\[
\mathcal{T}'=\{T(i,S_{r},Z):i\in S_{r},|T|\geq|S_{r}|^{\epsilon}\}
\]
where as in (\ref{eq:-20}), $T(i,S,Z)$ stands for maximal skeleton:
\[
T(i,S,Z)=\arg\max\{|T|:i\in T\subset S\,\&\,z_{T}\Rrightarrow S\}.
\]
 By Proposition \ref{prop:-is-a-1} together with Proposition \ref{prop:Let--be},
with probability $1-\Delta(L_{\varsigma,r})$ we have for all $T\in\mathcal{T}'$:
\begin{align*}
[Z_{T}] & =\{w\in C_{T}:\hat{w}^{a,b}=\hat{Z}_{T}^{a,b},w_{\mathcal{J}_{0}(Z_{T})}=Z_{\mathcal{J}_{0}(Z_{T})}\}\\{}
[Z_{S_{r}}] & =\{w\in C_{S_{r}}:\hat{w}^{a,b}=\hat{Z}_{S_{r}}^{a,b},w_{\mathcal{J}_{0}(Z_{S_{r}})}=Z_{\mathcal{J}_{0}(Z_{S_{r}})}\}
\end{align*}
thus, as $\mathcal{J}_{0}(Z_{T})\subset\mathcal{J}_{0}(Z_{S_{r}})$
(by definition) and $T\subset S_{r}$, we get with the same probability
that 
\[
[Z_{S_{r}}]\Rrightarrow\hat{Z}_{S_{r}}^{a,b}\,\wedge\,\mathcal{J}_{0}(Z_{S_{r}})\Rrightarrow\hat{Z}_{T}^{a,b}\,\wedge\,\mathcal{J}_{0}(Z_{T})\Rrightarrow[Z_{T}].
\]
Which prove one direction of the claim. The other direction proved
in the same way.$\hfill\square$

\subsection{Proof of the main theorem}

As in Sec. \ref{subsec:Proof-of-Lemma}, $Z$ is assumed to be an
independent splitting of a bi-marked process with marked states $a,b$,
$\mathbf{C}_{a}=\left(C^{\mathbb{Z}},\mu,a\right)$ and $\mathbf{C}_{b}=\left(C^{\mathbb{Z}},\mu,b\right)$
are endowed with nice marker factors $\varsigma$ and $\vartheta$
respectively, with associated lengths sequences \textit{$\left(L_{\varsigma,r}\right)_{r},\left(L_{\vartheta,r}\right)_{r}$}
, referring to their skeletons by $S$ and $T$ respectively. $[Z_{S}]$
and $[Z_{T}]$ refer to their corresponding filler equivalence class
in the good SOFs that where guaranteed to exist in Lemma \ref{Nested=000020fillers=000020in=000020r=000020property}.
\begin{prop}
\label{prop:With-the-same}Suppose $\phi(Z)$ is a factor of $Z$
satisfying
\begin{equation}
\mathbf{P}\left[\frac{1}{|T_{r}|}\sum_{i}\mathbf{1}\{[Z_{T_{r}}]\Rrightarrow\phi(Z)_{i}\}\geq1-\tilde{O}(|T_{r}|^{-\frac{1}{2}})\right]=1-\Delta(L_{\vartheta,r}),\label{eq:-30}
\end{equation}
then one has also
\[
\mathbf{P}\left[\frac{1}{|S_{r}|}\sum_{i}\mathbf{1}\{[Z_{S_{r}}]\Rrightarrow\phi(Z)_{i}\}\geq1-\tilde{O}(|S_{r}|^{-\frac{1}{2}})\right]=1-\Delta(L_{\varsigma,r}).
\]
\end{prop}

\begin{proof}
Let $Q$ be the collection of all fillers $w$ of $\vartheta$-skeletons
$T$ that satisfy
\[
\frac{1}{|T|}\sum_{i\in T}\mathbf{1}\{[w]\Rrightarrow\phi(w)_{i}\}\geq1-\tilde{O}(|T|^{-\frac{1}{2}}).
\]
By (\ref{eq:-30}), 
\[
\mathbf{P}[Z_{T_{r}}\in Q]=1-\Delta(L_{\vartheta,r})
\]
thus Proposition \ref{prop:Let--be} implies that with probability
$1-\Delta(L_{\varsigma,r})$, for all $i\in\mathring{S}_{r}^{\epsilon}$,
$Z_{T(i,S_{r},Z)}\in Q$. That, together with the additional property
of Lemma \ref{Nested=000020fillers=000020in=000020r=000020property}
and Proposition \ref{prop:For-any-,} (that with probability $1-\Delta(L_{\varsigma,r})$,
for all $i\in\mathring{S}_{r}^{\epsilon}$, $[Z_{S_{r}}]\Rrightarrow[Z_{T(i)}]$
and $|T(i)|\geq\Lambda_{-c}d(i,S_{r}^{c})$ ), implies that with probability
$1-\Delta(L_{\varsigma,r})$,
\begin{align*}
\frac{1}{|S_{r}|}\sum_{i\in S_{r}}\mathbf{1}\{[Z_{S_{r}}]\Rrightarrow\phi(Z)_{i}\} & \geq\frac{1}{|S_{r}|}\sum_{i\in\mathring{S}_{r}^{\epsilon}}\frac{1}{|T(i)|}\sum_{j\in T(i)}\mathbf{1}\{[Z_{T(i)}]\Rrightarrow\phi(Z)_{j}\}\\
 & \geq1-2|S_{r}|^{\epsilon}-\frac{1}{|S_{r}|}\sum_{i\in\mathring{S}_{r}^{\epsilon}}\tilde{O}(|T(i)|^{-\frac{1}{2}})\\
 & \geq1-2|S_{r}|^{\epsilon}-\frac{1}{|S_{r}|}\sum_{i\in\mathring{S}_{r}^{\epsilon}}\tilde{O}(d(i,S_{r}^{c})^{-\frac{1}{2}})\\
 & \geq1-\tilde{O}(|S_{r}|^{-\frac{1}{2}})
\end{align*}
and the claim follows.
\end{proof}
Suppose $X$ and $Y$ are two marked processes with equal entropy
and equal marked state distribution. We recall that the marker process
lemma (Lemma \ref{lem:Marker=000020process} and Corollary \ref{cor:Marker=000020process})
provides a common nice marker factor $\xi$ with lengths sequence
$\left(L_{r}\right)_{r}$, while Lemma \ref{lem:good=000020systems=000020of=000020fillers=000020imply=000020good=000020isomorphism}
together with Proposition \ref{prop:Existence=000020of=000020good=000020filler=000020systems}
provide an isomorphism $X\overset{\phi}{\rightarrow}Y$ s.t. for any
$n$, there are $\{r_{i}\}_{i=1}^{4}$ with $L_{r_{i}}=\tilde{\Theta}(n)$,
satisfying with probability $1-\Delta(n)$
\begin{equation}
[X_{S_{r_{1}}}]\Rrightarrow[\phi(X)_{S_{r_{2}}}]\label{eq:-32}
\end{equation}
and
\begin{equation}
[Y_{S_{r_{3}}}]\Rrightarrow[\phi^{-1}(Y)_{S_{r_{4}}}].\label{eq:-33}
\end{equation}

\begin{prop}
\label{prop:In-the-above}In the above circumstance, suppose further
that $\psi(Y)$ is a factor of $Y$ with the property that
\begin{equation}
\frac{1}{|S_{r}|}\sum_{i}\mathbf{1}\{[Y_{S_{r}}]\Rrightarrow\psi(Y)_{i}\}\geq1-\tilde{O}(|S_{r}|^{-\frac{1}{2}})\label{eq:-34}
\end{equation}
 with probability $1-\Delta(L_{r})$. Then one has also
\[
\mathbf{P}\left[\frac{1}{|S_{r}|}\sum_{i}\mathbf{1}\{[X_{S_{r}}]\Rrightarrow\psi\circ\phi(X)_{i}\}\geq1-\tilde{O}(|S_{r}|^{-\frac{1}{2}})\right]=1-\Delta(L_{r}).
\]
\end{prop}

\begin{proof}
Given $r$, choose $n=\tilde{\Theta}(L_{r})$ but yet small enough
so that in (\ref{eq:-32}) we will have $r_{1}\leq r$. Writing $S_{r}(i)$
for the $r$-skeleton that covers $X$ at $i$, we have:
\begin{align*}
\frac{1}{|S_{r}|}\sum_{i}\mathbf{1}\{[X_{S_{r}}] & \Rrightarrow\psi\circ\phi(X)_{i}\}\geq\frac{1}{|S_{r}|}\sum_{i}\mathbf{1}\{[X_{S_{r_{1}}(i)}]\Rrightarrow\psi\circ\phi(X)_{i}\}\\
 & \geq\frac{1}{|S_{r}|}\sum_{i}\mathbf{1}\{[X_{S_{r_{1}}(i)}]\Rrightarrow[\phi(X)_{S_{r_{2}}(i)}]\}\mathbf{1}\{[\phi(X)_{S_{r_{2}}(i)}]\Rrightarrow\psi\circ\phi(X)_{i}\}\\
 & =\frac{1}{|S_{r}|}\sum_{i}\mathbf{1}\{[X_{S_{r_{1}}(i)}]\Rrightarrow[\phi(X)_{S_{r_{2}}(i)}]\}\left(\frac{1}{|S_{r_{2}}(i)|}\sum_{j\in S_{r_{2}}(i)}\mathbf{1}\{[\phi(X)_{S_{r_{2}}(i)}]\Rrightarrow\psi\circ\phi(X)_{j}\}\right)
\end{align*}
By (\ref{eq:-32}) and (\ref{eq:-34}), for any $i\in S_{r}$, the
event 
\[
\mathbf{1}\{[X_{S_{r_{1}}(i)}]\Rrightarrow[\phi(X)_{S_{r_{2}}(i)}]\}\left(\frac{1}{|S_{r_{2}}(i)|}\sum_{j\in S_{2}(i)}\mathbf{1}\{[\phi(X)_{S_{r_{2}}(i)}]\Rrightarrow\psi\circ\phi(X)_{j}\}\right)\geq1-\tilde{O}(|S_{r_{2}}|^{-\frac{1}{2}})
\]
is of probability $1-\Delta(L_{r})$. Intersecting on all these events
together with the event $|S_{r}|\leq L_{r}$ and the events $|S_{r_{2}}(i)|=\tilde{\Theta}(L_{r})$
(all of whose probability is $1-\Delta(L_{r})$ too) gives an event
of probability $1-(2L_{r}+1)\Delta(L_{r})=1-\Delta(L_{r})$, on which
\begin{align*}
\frac{1}{|S_{r}|}\sum_{i}\mathbf{1}\{[X_{S_{r}}]\Rrightarrow\psi\circ\phi(X)_{i}\} & \geq1-\frac{1}{|S_{r}|}\sum_{i}\tilde{O}(|S_{r_{2}}(i)|^{-\frac{1}{2}})\\
 & =1-\tilde{O}(L_{r}^{-\frac{1}{2}})
\end{align*}
just as claimed.
\end{proof}
Let $\left(X(i)\right)_{i=0}^{N}$ Be a sequence as in Lemma \ref{lem:marked=000020process=000020lemma},
and for each pair of consecutive elements of the sequence $X(i-1),X(i)$,
let $\xi^{i}$ and $\phi^{i}$ be the associated nice marker factor
and isomorphism provided by the marker process lemma, lemma \ref{lem:good=000020systems=000020of=000020fillers=000020imply=000020good=000020isomorphism}
and proposition \ref{prop:Existence=000020of=000020good=000020filler=000020systems}
as above.
\begin{cor}
\label{cor:-is-a}$\phi=\phi^{N}\circ\cdots\circ\phi^{1}$ is a finitary
isomorphism from $X(0)$ to $X(n)$ satisfying
\begin{align}
\mathbf{P}\left[\frac{1}{|S_{r}|}\sum_{i}\mathbf{1}\{[X(0)_{S_{r}}]\Rrightarrow\phi(X(0))_{i}\}\geq1-\tilde{O}(|S_{r}|^{-\frac{1}{2}})\right] & =1-\Delta(L_{\xi^{1},r})\label{eq:-35}\\
\mathbf{P}\left[\frac{1}{|T_{r}|}\sum_{i}\mathbf{1}\{[X(N)_{T_{r}}]\Rrightarrow\phi^{-1}(X(N))_{i}\}\geq1-\tilde{O}(|T_{r}|^{-\frac{1}{2}})\right] & =1-\Delta(L_{\xi^{N},r})\label{eq:-36}
\end{align}
where $S_{r}$ is taken with respect to $\xi^{1}$ and $T_{r}$ is
taken with respect to $\xi^{N}$.
\end{cor}

\begin{proof}
It is clear that $\phi$ is a finitary isomorphism. the validity of
(\ref{eq:-35}) is shown by induction: The case $N=1$ is just theorem
\ref{thm:common=000020entry=000020version} (or alternatively, use
(\ref{eq:good=000020expected=000020number=000020of=000020determined=000020coordinates})
of proposition \ref{eq:good=000020expected=000020number=000020of=000020determined=000020coordinates}
and then use proposition \ref{prop:In-the-above} with $\psi$ being
the identity). Suppose the claim holds for $N-1\geq1$. Then we have
that $\phi'=\phi^{N}\circ\cdots\circ\phi^{2}$ is a finitary isomorphism
satisfying 
\[
\mathbf{P}\left[\frac{1}{|T_{r}|}\sum_{i}\mathbf{1}\{[X(1)_{T_{r}}]\Rrightarrow\phi'(X(1))_{i}\}\geq1-\tilde{O}(|T_{r}|^{-\frac{1}{2}})\right]=1-\Delta(L_{\xi^{2},r})
\]
 with $T_{r}=T_{r}(X(1))$ being taken with respect to $\xi^{2}$.
This together with proposition \ref{prop:With-the-same} implies that
\[
\mathbf{P}\left[\frac{1}{|S_{r}|}\sum_{i}\mathbf{1}\{[X(1)_{S_{r}}]\Rrightarrow\phi'(X(1))_{i}\}\geq1-\tilde{O}(|S_{r}|^{-\frac{1}{2}})\right]=1-\Delta(L_{\xi^{1},r})
\]
 with $S_{r}=S_{r}(X(1))$ being taken with respect to $\xi^{1}$.
Now proposition \ref{prop:In-the-above} will give
\[
\mathbf{P}\left[\frac{1}{|S_{r}|}\sum_{i}\mathbf{1}\{[X(0)_{S_{r}}]\Rrightarrow\phi^{1}\circ\phi'(X(0))_{i}\}\geq1-\tilde{O}(|S_{r}|^{-\frac{1}{2}})\right]=1-\Delta(L_{\xi^{1},r})
\]
which proves the induction step. Hence equation (\ref{eq:-35}) follows.
equation (\ref{eq:-36}) is proved in the same manner.
\end{proof}
Our main theorem now follows from the results of this section:
\begin{thm}
(Main Theorem) For any two finite alphabet aperiodic Markov shifts
$X\sim(A^{\mathbb{Z}},\mu),Y\sim(B^{\mathbb{Z}},\nu)$ of the same
entropy, there exists a finitary isomorphism $\Phi:A^{\mathbb{Z}}\rightarrow B^{\mathbb{Z}}$
with coding radii $R_{\Phi},R_{\Phi^{-1}}$ satisfying:
\[
\mathbf{P}[R_{\Phi}>n]=\tilde{O}(n^{-\frac{1}{2}}),\qquad\mathbf{P}[R_{\Phi^{-1}}>n]=\tilde{O}(n^{-\frac{1}{2}}).
\]
\end{thm}

\begin{proof}
By lemma \ref{lem:marked=000020process=000020lemma} together with
corollary \ref{cor:-is-a}, there is some $k\geq1$ such that $X^{(k)}$
and $Y^{(k)}$ are isomorphic by a finitary isomorphism $\phi$, so
that for some nice marker factors $\xi^{1}=\xi^{1}(X^{(k)})$ and
$\xi^{N}=\xi^{N}(Y^{(k)})$ equations (\ref{eq:-35}) and (\ref{eq:-36})
take place with $X(0)=X^{(k)}$ and $X(N)=Y^{(k)}$. Given $n>k$,
as $\xi^{1}$ is nice, there is some $r=r(n-k)$ so that 
\[
\mathbf{P}[n-k\geq|S_{r}(X^{(k)})|=\tilde{\Theta}(n)]=1-\Delta(n)
\]
thus $\mathbf{P}[X_{[-n,n-k]}^{(k)}\Rrightarrow X_{S_{r}}^{(k)}]=1-\Delta(n)$.
We also have that $\mathbf{P}[X_{[-n,n]}\Rrightarrow X_{[-n,n-k]}^{(k)}]=1$,
$\mathbf{P}[X_{S_{r}}^{(k)}\Rrightarrow[X_{S_{r}}^{(k)}]]=1$, and
(\ref{eq:-35}) interprets to
\[
\mathbf{P}\text{[[\ensuremath{X_{S_{r}}^{(k)}}]}\Rrightarrow\phi(X_{S_{r}}^{(k)})_{0}]=1-\tilde{O}(n^{-\frac{1}{2}}).
\]
Thus, writing $\psi(Z)=Z^{(k)}$ for any of the processes $Z\in\{X,Y\}$,
we get that $\Phi:=\psi^{-1}\circ\phi\circ\psi$ is an isomorphism
from $X$ to $Y$, and
\begin{align*}
\mathbf{P}[X_{[-n,n]}\Rrightarrow\Phi(X)_{0}] & \geq\mathbf{P}[X_{[-n,n-k]}^{(k)}\Rrightarrow\phi(X^{(k)})_{0}]\\
 & =\mathbf{P}[[X_{S_{r}}^{(k)}]\Rrightarrow\phi(X^{(k)})_{0}]-\Delta(n)\\
 & =1-\tilde{O}(n^{-\frac{1}{2}})
\end{align*}
which shows that $\mathbf{P}[R_{\Phi}>n]=\tilde{O}(n^{-\frac{1}{2}})$.
The proof of the tail bound for $R_{\Phi^{-1}}$ follows the same
lines.
\end{proof}

\lyxaddress{uriel.gabor@gmail.com}

\lyxaddress{Einstein Institute of Mathematics}

\lyxaddress{The Hebrew University of Jerusalem}

\lyxaddress{Edmond J. Safra Campus, Jerusalem, 91904, Israel}

\begin{thebibliography}{10}
\bibitem{key-14}Tim Austin. Entropy and Ergodic Theory (lecture notes),
Lecture 15: A first look at concentration. https://www.math.ucla.edu/\textasciitilde tim/entropy\_15.pdf

\bibitem{key-18}L. Bowen, Every countably infinite group is almost
Ornstein, in Dynamical Systems and Group Actions, Contemp. Math.,
567, Amer. Math. Soc., Providence, RI, 2012, 67--78.

\bibitem{key-19}M. Denker, M. Keane, Almost topological dynamical
systems. Israel J. Math. 34, 139--160 (1979).

\bibitem{key-1}Nate Harvey, Alexander E. Holroyd, Yuval Peres, and
Dan Romik, Universal finitary codes with exponential tails, Proceedings
of the London Mathematical Society 94 (2006), no. 2, 475--496.

\bibitem{key-9}Nate Harvey and Yuval Peres, An invariant of finitary
codes with finite expected square root coding length. Ergodic Theory
and Dynamical Systems, 31 (2011), no. 1, 77-90.

\bibitem{key-3}Michael Keane and Meir Smorodinsky, A class of finitary
codes, Israel Journal of Mathematics 26 (1977), no. 3-4, 352--371.

\bibitem{key-2}Michael Keane and Meir Smorodinsky, Bernoulli schemes
of the same entropy are finitarily isomorphic, Annals of Mathematics
109 (1979), no. 2, 397--406.

\bibitem{key-17}Michael Keane and Meir Smorodinsky, Finitary isomorphisms
of irreducible Markov shifts Israel J. Math., 34 (4) (1980), pp. 281-286

\bibitem{key-10}A. N. Kolmogorov. A new metric invariant of transient
dynamical systems and automorphisms in Lebesgue spaces. Dokl. Akad.
Nauk SSSR (N.S.), 119:861--864, 1958.

\bibitem{key-8}Krieger, W.: On the finitary isomorphisms of Markov
shifts that have finite expected coding time. Z. Wahrscheinlichkeitstheorie
und verw. Gebiete 65, 323-328 (1983)

\bibitem{key-4}L. D. Meshalkin (1959), A case of isomorphism of Bernoulli
schemes. Dokl. Akad. Nauk SSSR, 128, 41--44.

\bibitem{key-5}D. S. Ornstein (1970), Bernoulli shifts of the same
entropy are isomorphic. Adv. in Math. 4, 337--352.

\bibitem{key-6}W. Parry (1979), Finitary isomorphisms with finite
expected code lengths. Bull. London Math. Soc. 11, 170--176. 

\bibitem{key-12}Parry, W., Schmidt, K.: Natural coefficients and
invariants for Markov shifts. Invent. math.76, 15--32 (1984)

\bibitem{key-7}K. Schmidt (1984), Invariants for finitary isomorphisms
with finite expected code lengths. Invent. Math. 76, 33--40.

\bibitem{key-11}J. Serafin (1996), The finitary coding of two Bernoulli
schemes with unequal entropies has finite expectation. Indag. Math.
(N.S.) 7, no. 4, 503--519.

\bibitem{key-24}Seward, Brandon. (2018). Bernoulli shifts with bases
of equal entropy are isomorphic. https://arxiv.org/abs/1805.08279

\bibitem{key-15}Michel Talagrand. A new look at independence. Ann.
Probab., 24(1):1-- 34, 1996.

\end{thebibliography}
\end{document}